\documentclass{amsart}
\usepackage[margin=1.5in]{geometry}
\usepackage{amsthm}
\usepackage{amssymb}
\usepackage{mathrsfs}
\usepackage[non-compressed-cites]{amsrefs}
\usepackage{esint}
\usepackage{caption}
\usepackage{enumitem}
\usepackage{mathtools}
\usepackage{xcolor}
\usepackage{microtype}
\usepackage{multirow}
\usepackage{array}
\usepackage{tikz}
\usepackage{threeparttable}
\usepackage{hyperref}

\usetikzlibrary{arrows.meta,decorations.markings,calc}

\allowdisplaybreaks

\newcommand{\R}{\mathbb{R}}
\newcommand{\Sph}{\mathbb{S}}
\newcommand{\del}{\partial}
\newcommand{\Ric}{\operatorname{Ric}}
\renewcommand{\div}{\operatorname{div}}
\newcommand{\Rm}{\operatorname{Rm}}
\newcommand{\tr}{\operatorname{tr}}
\newcommand{\vol}{\operatorname{vol}}
\newcommand{\dvol}{\mathrm{dvol}}
\newcommand\inner[2]{\langle #1, #2 \rangle}
\newcommand{\norm}[1]{\left\lVert #1 \right\rVert}
\renewcommand{\O}{\mathcal{O}}
\newcommand{\Z}{\mathbb{Z}}
\newcommand{\CP}{\mathbb{CP}}
\newcommand{\C}{\mathbb{C}}
\renewcommand{\epsilon}{\varepsilon}
\renewcommand{\Re}{\operatorname{Re}}
\newcommand{\diag}{\operatorname{diag}}
\DeclareMathOperator{\Span}{span}
\newcommand{\Hom}{\operatorname{Hom}}

\newcommand{\bfe}{\mathbf{e}}

\newcommand{\bfM}{\mathbf{M}}
\newcommand{\bfp}{\mathbf{p}}
\newcommand{\bfq}{\mathbf{q}}
\newcommand{\bfr}{\mathbf{r}}
\newcommand{\bfu}{\mathbf{u}}

\newcommand{\sfA}{\mathsf{A}}
\newcommand{\sfB}{\mathsf{B}}
\newcommand{\sfC}{\mathsf{C}}
\newcommand{\sfE}{\mathsf{E}}

\newcommand{\sfG}{\mathsf{G}}

\newcommand{\sfK}{\mathsf{K}}
\newcommand{\sfM}{\mathsf{M}}
\newcommand{\sfN}{\mathsf{N}}
\newcommand{\sfT}{\mathsf{T}}

\newcommand{\calB}{\mathcal{B}}
\newcommand{\calD}{\mathcal{D}}
\newcommand{\calE}{\mathcal{E}}
\newcommand{\calH}{\mathcal{H}}
\newcommand{\calK}{\mathcal{K}}
\newcommand{\calM}{\mathcal{M}}

\newcommand{\calQ}{\mathcal{Q}}
\newcommand{\calR}{\mathcal{R}}

\newcommand{\fka}{\mathfrak{a}}

\newcommand{\scrE}{\mathscr{E}}

\newcommand{\scrS}{\mathscr{S}}
\newcommand{\scrT}{\mathscr{T}}

\newtheorem{theorem}{Theorem}[section]
\newtheorem{corollary}[theorem]{Corollary}
\newtheorem{lemma}[theorem]{Lemma}

\newtheorem{proposition}[theorem]{Proposition}

\newtheorem{remark}[theorem]{Remark}
\numberwithin{equation}{section}
 
\hypersetup{
    colorlinks,
    linkcolor={red!60!black},
    citecolor={blue!80!black},
}

\usepackage{etoolbox}

\patchcmd{\subsection}{-.5em}{.5\baselineskip}{}{}
\patchcmd{\subsubsection}{-.5em}{.5\baselineskip}{}{}
\patchcmd{\section}
  {\normalfont\scshape\centering}
  {\normalfont\large\bfseries\centering}
  {}{}

\setlist[itemize]{leftmargin=3em}
\setlist[enumerate]{leftmargin=3em}
\mathtoolsset{showonlyrefs}

\title[Uniqueness and stability of Appleton's solitons]{Uniqueness and stability of Appleton's solitons}
\author[Michael B. Law]{Michael B. Law}
\address{MIT, Department of Mathematics, 77 Massachusetts Avenue, Cambridge, MA 02139, USA.}
\email{mikelaw@mit.edu}

\begin{document}

\begin{abstract}
    We prove that Appleton's noncollapsed steady gradient Ricci solitons on complex line bundles over positive K\"ahler--Einstein manifolds are the unique such solitons arising from his metric ansatz. We then show that the four-dimensional examples on $\O_{\CP^1}(-k)$, $k \geq 3$, are strictly linearly stable. By a result of Deruelle, this implies dynamical stability under the Ricci flow.
\end{abstract}

\maketitle

\tableofcontents

\section{Introduction}

Ricci solitons are self-similar solutions of Ricci flow and are key objects in the analysis of singularities. Finite-time type I singularities, characterized by the curvature bound $\sup_M |{\Rm}_{g(t)}|_{g(t)} \leq C(T-t)^{-1}$, have shrinking Ricci solitons as blowup models. Steady solitons often arise as models for Type II singularities and thus play a central role in their analysis. Steady and expanding solitons also provide natural candidates for gluing constructions aimed at continuing the flow past singularities.

The appearance of Ricci solitons as singularity models naturally leads to the question of their \emph{stability}: heuristically, whether flows starting near a soliton remain close or converge to it. Stability is a classical topic in the variational theory of Ricci solitons (and Einstein metrics), and plays an essential role in determining which solitons might occur in `generic' solutions to Ricci flow. In mean curvature flow, Colding and Minicozzi's classification of stable self-shrinkers \cite{cm-generic}, and numerous subsequent developments (see for instance \cites{ccms,cms,bamler-kleiner-m1}), have led to a deep understanding of generic mean curvature flows of embedded surfaces in $\R^3$. Developing a corresponding theory of generic Ricci flow remains an important goal, particularly in dimension four where the next frontier of Ricci flow stands.

We briefly survey some results on the stability of Ricci solitons. These involve several related notions of stability; the precise notion used in this paper will be given below. In the compact shrinking setting, Cao--Hamilton--Ilmanen \cite{cao-hamilton-ilmanen} and Cao--He \cite{cao-he} established linear (in)stability of various positive Einstein manifolds. Interestingly, $\CP^m$ is linearly stable but, as was discovered by Kr\"oncke \cites{kroncke-cpn,knopf-sesum-cpn}, dynamically unstable under Ricci flow when $m \geq 2$. In the noncompact steady and expanding settings, Deruelle \cite{deruelle-stability} established curvature conditions that imply strict linear stability, which in particular hold on the Bryant steady soliton and a sub-family of the Bryant expanding solitons. Deruelle--Kr\"oncke \cite{deruelle-kroncke} and Deruelle--Ozuch \cite{deruelle-ozuch} have studied the linear and dynamical (in)stability of Ricci-flat ALE metrics. More recently, Naff and Ozuch \cites{naff-ozuch-fik,naff-ozuch-kahler} proved linear stability of the four-dimensional FIK shrinker, linear instability of the four-dimensional BCCD shrinker, and linear stability of all four-dimensional steady and expanding K\"ahler--Ricci solitons. Lavoyer and Thompson \cite{lavoyer-thompson} recently extended the latter stability result to arbitrary dimension and proved strict linear stability of asymptotically conical K\"ahler--Ricci expanders. Very recently, Kim and Ozuch \cite{kim-ozuch} established dynamical stability results for the four-dimensional FIK shrinker.

We now turn to the complete, noncollapsed, non-K\"ahler steady gradient Ricci solitons constructed by Appleton \cite{appleton}, which are the subject of this paper. For any integers $k > m \geq 2$, his construction gives a $U(m)$-invariant soliton on the total space of the complex line bundle $\O_{\CP^{m-1}}(-k)$. The complement of the zero section is diffeomorphic to $(0,\infty) \times (\Sph^{2m-1}/\Z_k)$, where the $\Z_k$-action on $\Sph^{2m-1} \subset \C^m$ is generated by scalar multiplication by a $k$-th root of unity. On this set, the metric and soliton potential take the form
\begin{align} \label{eq:appleton-ansatz}
    g = ds^2 + a(s)^2 \sigma^2 + b(s)^2 \pi^*g_{\mathrm{FS}}, \qquad f=f(s).
\end{align}
Here $\pi: \Sph^{2m-1} \to \CP^{m-1}$ is the Hopf fibration which descends to $\Sph^{2m-1}/\Z_k$, $\sigma$ is the Hopf $1$-form on $\Sph^{2m-1}$ which also descends to $\Sph^{2m-1}/\Z_k$, and $g_{\mathrm{FS}}$ is the Fubini--Study metric with holomorphic sectional curvature $4$. Appleton's noncollapsed solitons satisfy
\begin{align} \label{eq:appleton-ab-asymptotics}
    a(s) \sim b(s) \sim C\sqrt{s} \qquad \text{as } s \to \infty
\end{align}
for some $C>0$.

More generally, Appleton constructed analogous steady gradient Ricci solitons on complex line bundles over positive K\"ahler--Einstein manifolds. We refer to \cite{appleton}*{\S 1.4} for the details of the construction. After normalizing the base metric, the soliton ODEs and boundary conditions take the same form as the $\O_{\CP^{m-1}}(-k)$ examples above. By analyzing these ODEs, we establish uniqueness of the noncollapsed solitons within Appleton's metric ansatz. This resolves the uniqueness part of Appleton's conjecture \cite{appleton}*{Conjecture 1}.

\begin{theorem} \label{thm:appleton-uniqueness}
    Fix a positive K\"ahler--Einstein base metric and a complex line bundle covered by \cite{appleton}*{Theorem 1.5}. Within Appleton's metric ansatz on this bundle, there is, up to scaling and isometry, exactly one complete noncollapsed steady gradient Ricci soliton.
\end{theorem}

\begin{remark}
    Shortly before releasing this paper, Appleton's full conjecture (which includes Theorem \ref{thm:appleton-uniqueness}) was independently resolved by Ribeiro and Tran \cite{ribeiro-tran}*{Theorem 6.1} in the four-dimensional case, and by Deng and Wang \cite{deng-wang}*{Theorem 1.3} in the general case. The ODE arguments presented in this paper appear quite similar to those in \cite{deng-wang}.
\end{remark}

Our other main result establishes \emph{strict linear stability} of the four-dimensional Appleton solitons, which occur on $\O_{\CP^1}(-k)$ for $k \geq 3$. To define this, let $(M,g,f)$ be a complete steady gradient Ricci soliton, with the convention
\begin{align}
    \Ric = \nabla^2 f.
\end{align}
Define the action of the curvature tensor on symmetric $2$-tensors $h$ by
\begin{align}
    \Rm(h)_{ij} = \Rm_{kijl} h^{kl},
\end{align}
where the curvature convention used is
\begin{align} \label{eq:Rm-convention}
    \Rm(X,Y)Z = \nabla_X \nabla_Y Z - \nabla_Y \nabla_X Z - \nabla_{[X,Y]}Z.
\end{align}
We say that $(M,g,f)$ is \emph{strictly linearly stable} if there exists $\mu > 0$ such that for every smooth compactly supported symmetric $2$-tensor $h$,
\begin{align} \label{eq:strict-linear-stab}
    \calE[h] := \int_M \left( |\nabla h|^2 - 2\inner{\Rm(h)}{h} \right) e^f \geq \mu \int_M |h|^2 e^f.
\end{align}
Here and throughout, all integrals of functions over $M$ are taken with respect to the Riemannian volume measure $\dvol_g$ unless otherwise stated.

\begin{theorem} \label{thm:appleton-stability}
    For every $k \geq 3$, the Appleton soliton on $\O_{\CP^1}(-k)$ is strictly linearly stable.
\end{theorem}

The main difficulties in proving stability are that Appleton's solitons have sectional curvatures of mixed signs and are non-K\"ahler. Thus, stability does not follow directly from existing results based on curvature conditions or K\"ahlerity. Instead, we exploit the cohomogeneity one symmetry and conformally K\"ahler structure of Appleton's solitons to reduce the stability inequality to ODE estimates and explicit matrix inequalities.

Theorem \ref{thm:appleton-uniqueness} is used in the author's classification of stable, asymptotically $SU(2)/\Gamma$-cylindrical steady gradient Ricci solitons as quotients of Bryant or Appleton solitons \cite{symmetry-law}. Theorem \ref{thm:appleton-stability} establishes the converse for the Appleton examples, whose stability properties were previously unknown.

Moreover, since Appleton's solitons have nonnegative Ricci curvature with bounded and decaying curvature at infinity (see \S\ref{subsec:local-frame}), Theorem \ref{thm:appleton-stability} combined with \cite{deruelle-stability}*{Theorems 1.2 and 1.3} yields dynamical stability under Ricci flow:

\begin{corollary} \label{cor:appleton-dynamically-stable}
    For every $k \geq 3$, the Appleton soliton $(M,g,f)$ on $\O_{\CP^1}(-k)$ is dynamically stable in the following sense.
    For each $I>0$, there exists $\epsilon = \epsilon(g,f,I) > 0$ such that whenever a smooth Riemannian metric $\tilde{g}$ satisfies
    \begin{align}
        \norm{\tilde{g}-g}_{L^2(e^f)} \leq I, \qquad \norm{\tilde{g}-g}_{L^\infty} < \epsilon,
    \end{align}
    there exists an immortal Ricci flow $(g(t))_{t \in [0,\infty)}$ with $g(0) = \tilde{g}$ and a smooth family of diffeomorphisms $(\phi_t)_{t \in [0,\infty)}$ of $M$, with $\phi_0 = \mathrm{Id}_M$, such that
    \begin{align}
        \phi_t^*g(t) \to g \quad \text{as } t \to \infty
    \end{align}
    in global $C^\ell$ norms for every integer $\ell \geq 0$. All norms are with respect to $g$.
\end{corollary}

The paper is organized as follows.

In \S\ref{sec:appleton-uniqueness}, we prove Theorem \ref{thm:appleton-uniqueness} using an ODE comparison argument.

In \S\ref{sec:further-preps-stability}, we establish further ODE estimates, record the connection and curvature in a convenient frame, and state additional conventions. We then give a sufficient criterion for strict linear stability: it is enough to show that \eqref{eq:strict-linear-stab} holds with an extra $b^{-2}$ factor in the right-hand integrand, and it suffices to show this for tensors supported away from the zero section.

In \S\ref{sec:appleton-separation}, we exploit the symmetry of the four-dimensional Appleton solitons to decompose certain tensor bundles into sectors labelled by pairs $(j,q)$, corresponding to irreducible representations of $U(2)$. For symmetric $2$-tensors, this decomposition allows stability to be addressed one sector at a time, provided the constant $\mu>0$ in \eqref{eq:strict-linear-stab} is uniform over all sectors. This generalizes the classical separation of variables method, and is essentially a reformulation of the method used in Naff and Ozuch's analysis of the FIK shrinker \cite{naff-ozuch-fik}.

In \S\ref{sec:charged-stability}, we establish the stability inequality in sectors with $q \neq 0$. The proof uses the Weitzenb\"ock formula of Naff--Ozuch \cite{naff-ozuch-kahler} and the conformally K\"ahler structure of Appleton's solitons to lower bound the energy quantity on the left-hand side of \eqref{eq:strict-linear-stab}. We conclude the argument using ODE estimates.

In \S\ref{sec:neutral-stability}, we establish the stability inequality in sectors with $q=0$. The proof uses tensor harmonics on the $2$-sphere to calculate $\calE[h]$ for $h$ belonging to such sectors, expressing it as a radial integral. The stability inequality reduces to a coercivity bound on this integral. After some algebraic manipulations, the desired estimates follow from a small number of matrix inequalities which are verified by exact symbolic computation.

In \S\ref{sec:conclusion}, we combine these results to establish strict linear stability, proving Theorem \ref{thm:appleton-stability}.

Accompanying the paper is a Wolfram Mathematica notebook containing symbolic computations which establish the matrix inequalities at the end of \S\ref{sec:neutral-stability}. The notebook also contains symbolic checks for many of the formulas appearing in the paper.

\begingroup
\renewcommand{\addcontentsline}[3]{}

\subsection*{Acknowledgements}

The author thanks William Minicozzi for continual support. The author also thanks Tristan Ozuch for bringing his work with Keaton Naff to the author's attention, and for suggesting conformal K\"ahlerity as an avenue to proving stability of Appleton's solitons. While this work was carried out, the author was supported by a Croucher Scholarship and a Simons Dissertation Fellowship.

\subsection*{Statement on AI usage}

Large language models (LLMs) suggested ODE arguments and explicit numerical estimates which appear in the following sections of the paper: \ref{subsec:new-variables}, \ref{subsec:uniqueness-proof}, \ref{subsec:4d-background-estimates}, \ref{subsec:scalar-coeff-bounds}, \ref{subsec:charged-proof}, \ref{subsec:reduction-pointwise}, \ref{subsec:neutral-proof}. In all cases, substantial human input was contributed to streamline the proofs and to render them human-readable. All other arguments were human-generated, with LLMs only used to help check correctness. LLMs also assisted in the writing and debugging of the accompanying Mathematica notebook, and in proofreading the manuscript. All AI-generated arguments occurring in the paper were manually vetted, and the author takes full responsibility for the correctness of all claims made.

\par
\endgroup

\section{Uniqueness of Appleton's solitons} \label{sec:appleton-uniqueness}

In this section, we prove Theorem \ref{thm:appleton-uniqueness} by comparing solutions of the soliton ODEs.

\subsection{The ansatz} \label{subsec:appleton-ODE-setup}

For integers $m \geq 2$ and $k>m$, let
\begin{align}
    M = \O_{\CP^{m-1}}(-k), \qquad B \cong \CP^{m-1}
\end{align}
denote the line bundle and its zero section. Write
\begin{align}
    M^\circ = M \setminus B \cong (0,\infty) \times L_{m,k}, \qquad L_{m,k} = \Sph^{2m-1}/\Z_k,
\end{align}
where $\Z_k$ acts by scalar multiplication by a $k$-th root of unity on $\Sph^{2m-1} \subset \C^m$.

The Hopf fibration
\begin{align}
    \pi: \Sph^{2m-1} \to \CP^{m-1}, \qquad (z_1,\ldots,z_m) \to [z_1: \ldots : z_m]
\end{align}
descends to a map $L_{m,k} \to \CP^{m-1}$ also denoted by $\pi$. Let $\sigma$ be the $1$-form on $\Sph^{2m-1}$ dual to the Hopf vector field $\xi(z_1,\ldots,z_m) = i(z_1,\ldots,z_m)$ with respect to the unit round metric, i.e. the round metric of sectional curvature $1$. Then $\sigma$ descends to a $1$-form on $L_{m,k}$ also denoted by $\sigma$. Thus, the unit round metric $g_{\Sph^{2m-1}} = \sigma^2 + \pi^*g_{\mathrm{FS}}$ induces the unit round quotient metric $g_{L_{m,k}} = \sigma^2 + \pi^*g_{\mathrm{FS}}$.

On $M^\circ$, we use Appleton's ansatz \eqref{eq:appleton-ansatz} for the metric and potential function. With our sign convention for $f$, the steady soliton equations are \cite{appleton}*{Appendix A}
\begin{align}
    f'' &= -\frac{a''}{a} - (2m-2) \frac{b''}{b}, \label{eq:appleton-soliton-f} \\
    a'' &= (2m-2)\left( \frac{a^3}{b^4} - \frac{a'b'}{b} \right) - a'f', \label{eq:appleton-soliton-a} \\
    b'' &= \frac{2m}{b} - 2\frac{a^2}{b^3} - \frac{a'b'}{a} - (2m-3)\frac{(b')^2}{b} - b'f'. \label{eq:appleton-soliton-b}
\end{align}
We require $a,b > 0$ on $(0,\infty)$ for a Riemannian metric. For $g$ and $f$ to extend smoothly across the zero section, we require $a$ to extend smoothly and oddly across $s=0$, while $b$ and $f$ extend smoothly and evenly. We also require $a'(0) = k$ and $b(0) > 0$. We fix the metric scale and the additive constant in $f$ by imposing $b(0) = 1$ and $f(0) = 0$. The remaining initial parameter is $f''(0)$. We are interested in solutions defining complete noncollapsed steady gradient Ricci solitons; by the analysis in \cite{appleton}, this requires that the solutions are defined for all $s>0$ and satisfy $\lim_{s \to \infty} a(s)/b(s) = 1$.

We also review the essential parts of Appleton's more general construction, which takes place on certain complex line bundles over a positive K\"ahler--Einstein base $(B,\bar{g})$ of complex dimension $m-1$. More details can be found in \cite{appleton}*{\S 1.4 and \S 2}. In this construction, the complement of the zero section is identified with $(0,\infty) \times P$, where $\pi: P \to B$ is the associated circle bundle, and Appleton's ansatz is
\begin{align} \label{eq:appleton-ansatz-general}
    g = ds^2 + a(s)^2 \sigma^2 + b(s)^2 \pi^*\bar{g}, \qquad f=f(s),
\end{align}
where $\sigma$ is a fixed connection one-form. After suitable normalization, the steady soliton equations are described by the same ODE system \eqref{eq:appleton-soliton-f}--\eqref{eq:appleton-soliton-b}. Smoothness at the zero section also demands the same boundary conditions as above, except that
\begin{align}
    a'(0) = \frac{mk}{p},
\end{align}
where $p$ and $k$ are integers associated with the choice of base and line bundle. Appleton's existence theorem \cite{appleton}*{Theorem 1.5} covers the cases where $k>p$, so $a'(0) = mk/p > m$.

We therefore analyze the ODEs \eqref{eq:appleton-soliton-f}--\eqref{eq:appleton-soliton-b} and impose the following boundary conditions throughout:
\begin{equation}
\label{eq:appleton-boundary-conditions}
    \begin{gathered}
        a(0)=0,\qquad a'(0)=\kappa > m,\qquad
        b(0)=1,\qquad b'(0)=0, \qquad f(0)=f'(0)=0, \\
        a(s),b(s)>0\quad(s>0),\qquad
        \lim_{s\to\infty}\frac{a(s)}{b(s)}=1, \\
        a \text{ extends smoothly and oddly across } s=0, \\
        b,f \text{ extend smoothly and evenly across } s=0.
    \end{gathered}
\end{equation}

We will establish the following result:
\begin{theorem} \label{thm:appleton-ODE-uniqueness}
    Fix an integer $m \geq 2$ and a real number $\kappa > m$. Let $(f_i,a_i,b_i)$, $i=1,2$, be smooth triples arising from complete steady gradient Ricci solitons in Appleton's ansatz. Assume that these triples satisfy \eqref{eq:appleton-soliton-f}--\eqref{eq:appleton-soliton-b} and \eqref{eq:appleton-boundary-conditions} with the same $m$ and $\kappa$. Then $f_1''(0) = f_2''(0)$.
\end{theorem}

Using this, we can prove Theorem \ref{thm:appleton-uniqueness} as follows.
\begin{proof}[Proof of Theorem \ref{thm:appleton-uniqueness}]
    Existence follows from \cite{appleton}*{Theorem 1.5}. We prove uniqueness.
    
    Fix a positive K\"ahler--Einstein base and a complex line bundle $M$ covered by Appleton's existence theorem. Let $(M,g_i,f_i)$, $i=1,2$ be two complete noncollapsed steady gradient Ricci solitons on the total space of this line bundle, both satisfying Appleton's ansatz \eqref{eq:appleton-ansatz-general}. Then, after normalization, the solitons are uniquely described by triples of smooth functions $(f_i,a_i,b_i)$ satisfying \eqref{eq:appleton-soliton-f}--\eqref{eq:appleton-soliton-b} and \eqref{eq:appleton-boundary-conditions} with the same $m$ and $\kappa = mk/p$. By Theorem \ref{thm:appleton-ODE-uniqueness}, we must have $f_1''(0) = f_2''(0)$.
    
    We claim that each $(f_i,a_i,b_i)$ extends analytically across $s=0$. Fix $i$ and temporarily suppress the index. By \cite{ivey-analytic}, $g$ and $f$ are real-analytic in suitable coordinates. Let $\gamma: (-\epsilon,\epsilon) \to M$ be a unit-speed geodesic through the zero section such that $\gamma'(s) = \del_s$ for $s > 0$. Then $\gamma$ and parallel vector fields along it also analytically depend on the signed parameter $s$. In particular, $f \circ \gamma$ is an analytic extension of $f$.

    At the point $\gamma(\epsilon/2)$, let $E_a = a\sigma^\sharp$, and let $E_b$ be a unit tangent vector orthogonal to $\del_s$ and $E_a$. Extend $E_a$ and $E_b$ to unit vector fields along $\gamma$ by parallel transport. Then $E_a(s) = a\sigma^\sharp$ for $0<s<\epsilon$, and $E_b$ remains orthogonal to $\del_s$ and $E_a$. The sectional curvatures
    \begin{align}
        K_a(s) = \operatorname{sec}(\gamma'(s), E_a(s)), \qquad K_b(s) = \operatorname{sec}(\gamma'(s), E_b(s))
    \end{align}
    are analytic functions on $(-\epsilon,\epsilon)$. A computation using the metric \eqref{eq:appleton-ansatz-general} shows that they satisfy, for $s>0$,
    \begin{align}
        a'' + K_a a = 0, \qquad b'' + K_b b = 0.
    \end{align}
    These are regular linear ODEs with analytic coefficients across $s=0$, so their solutions with the prescribed initial data \eqref{eq:appleton-boundary-conditions} are analytic and, by uniqueness, agree with $a,b$ for $0 \leq s < \epsilon$. This proves the claim.

    We may now apply the local uniqueness theorem of \cite{appleton}*{Theorem 2.1}, which gives that $(f_1,a_1,b_1)$ and $(f_2,a_2,b_2)$ agree on the interval $[0,\epsilon)$ for some $\epsilon>0$. Since the soliton equations \eqref{eq:appleton-soliton-f}--\eqref{eq:appleton-soliton-b} are regular whenever $a,b>0$, uniqueness of solutions to ODEs extends this agreement to all of $[0,\infty)$. Thus $(f_1,a_1,b_1) = (f_2,a_2,b_2)$, which means that $(M,g_1)$ and $(M,g_2)$ are isometric up to scaling.
\end{proof}

\subsection{Basic properties of solutions} \label{subsec:appleton-qualititative}

The rest of \S\ref{sec:appleton-uniqueness} is devoted to proving Theorem \ref{thm:appleton-ODE-uniqueness}. Throughout, let $(f,a,b)$ be a triple satisfying the assumptions of Theorem \ref{thm:appleton-ODE-uniqueness}. The parity conditions in \eqref{eq:appleton-boundary-conditions} and the soliton equations give the following expansions as $s \downarrow 0$:
\begin{align}
    f(s) &= \frac{f''(0)}{2} s^2 + \O(s^4), \label{eq:bolt-expansion-f} \\
    a(s) &= \kappa s - \frac{\kappa (f''(0) + 2m(m-1))}{6} s^3 + \O(s^5), \label{eq:bolt-expansion-a} \\
    b(s) &= 1 + \frac{m}{2} s^2 + \O(s^4). \label{eq:bolt-expansion-b}
\end{align}
Indeed, taking the limit in \eqref{eq:appleton-soliton-b} gives $b''(0) = 2m-b''(0)$, so $b''(0) = m$. Comparing the coefficients of $s$ in \eqref{eq:appleton-soliton-a} then gives $a'''(0) = -\kappa (f''(0) + 2m(m-1))$.

By \eqref{eq:appleton-boundary-conditions}, the quantity $Q = a/b$ satisfies $\lim_{s \to \infty} Q(s) = 1$. As computed in \cite{appleton}*{\S 3}, its evolution equation is
\begin{align} \label{eq:appleton-Q-ODE}
    Q'' + \left( f' + (2m-1) \frac{b'}{b} \right) Q' = \frac{2m}{b^2}Q(Q^2-1).
\end{align}

We collect several relevant properties proved in \cite{appleton}. Hereafter, $R$ denotes scalar curvature.

\begin{lemma} \label{lem:appleton-qualitative-properties}
    We have
    \begin{align}
        f''(0) > 0, \qquad R>0.
    \end{align}
    For every $s>0$, we also have
    \begin{align}
        0<Q<1, \qquad Q'>0, \qquad a'>0, \qquad b'>0, \qquad f'>0, \qquad f''>0.
    \end{align}
    Furthermore,
    \begin{align} \label{eq:appleton-z-bound}
        0 < bQ' \leq \kappa e^{-f} b^{-2m+2}.
    \end{align}
\end{lemma}
\begin{proof}
    By \eqref{eq:appleton-Q-ODE}, every critical point of $Q$ with $0<Q<1$ is a strict local maximum, while every critical point with $Q>1$ is a strict local minimum.
    
    Since $\lim_{s \to \infty} Q(s) = 1$, we must then have $Q \leq 1$. So if $Q=1$ anywhere, then $Q'=0$ at that point, but uniqueness of solutions to the ODE \eqref{eq:appleton-Q-ODE} would then imply $Q=1$ everywhere, contradicting the boundary condition $Q(0)=0$. Thus $0<Q<1$ whenever $s>0$.
    
    The boundary conditions also give $Q'(0) = \kappa > 0$. Since $0<Q<1$, a first contact with $Q' = 0$ must be a strict local maximum of $Q$. Moreover, by \eqref{eq:appleton-Q-ODE}, every zero of $Q'$ has $Q''<0$, so $Q'$ would remain strictly negative thereafter. This contradicts $Q \to 1$. Hence $Q'>0$.

    By Lemma 4.1 in \cite{appleton}, it follows that $a',b' > 0$ whenever $s>0$. By Corollary 3.5 and Lemma 7.1 in that paper, we also have $f''(0) > 0$, so by Lemma 4.2 in that paper, $f'>0$ and $f''>0$ for all $s > 0$. Also, Lemma 3.4 gives $R = 2f''(0)>0$ at the zero section. For a complete steady gradient Ricci soliton, it is well-known that this implies $R>0$ everywhere.
        
    Finally, \eqref{eq:appleton-Q-ODE} gives
    \begin{align}
        (e^f b^{2m-1} Q')' = -2m e^f b^{2m-3} Q(1-Q^2) < 0.
    \end{align}
    Since $e^f b^{2m-1} Q' \to \kappa$ as $s \to 0$, this proves $0 < bQ' \leq \kappa e^{-f} b^{-2m+2}$.
\end{proof}

\subsection{A change of variables} \label{subsec:new-variables}

We now make a change of variables that will be useful going forward. We take
\begin{align} \label{eq:x-uniqueness}
    x = Q^2
\end{align}
as a new independent variable. Since $Q'>0$, $0<Q<1$, and $Q \to 1$ as $s \to \infty$, this defines a smooth increasing coordinate from $s \in (0,\infty)$ to $x \in (0,1)$. Furthermore, let
\begin{align} \label{eq:appleton-comparison-variables}
    v = \frac{b'}{Q}, \qquad y = af' + 2a' + (2m-2)\frac{ab'}{b}, \qquad z=bQ'.
\end{align}
We will use dots (e.g. $\dot{v}, \dot{y}, \dot{z}$) to denote differentiation with respect to $x$. Primes continue to denote differentiation with respect to $s$.

The definitions give
\begin{align} \label{eq:comparison-basic-identities}
    ax' = 2xz, \qquad a' = z+xv, \qquad \frac{ab'}{b} = xv, \qquad af' = y-2z-2mxv.
\end{align}
We now derive the relations between $v,y,z$, i.e. \eqref{eq:general-v-equation}--\eqref{eq:general-z-equation} below. First, differentiating $v$ and using \eqref{eq:appleton-soliton-b} and \eqref{eq:comparison-basic-identities}, we obtain
\begin{align}
    av' &= bb'' - vz \\
    &= 2m-2x-v(a'+z+af') - (2m-3) xv^2 \\
    &= 2m-yv-2x(1-v^2).
\end{align}
Next, differentiating $y$ and using \eqref{eq:appleton-soliton-f} and then \eqref{eq:appleton-soliton-a} gives
\begin{align}
    y' &= a'f' + af'' + 2a'' + (2m-2)\left( \frac{a'b'}{b} + \frac{ab''}{b} - \frac{a(b')^2}{b^2} \right) \\
    &= a'f' + a'' + (2m-2)\left( \frac{a'b'}{b} - \frac{a(b')^2}{b^2} \right) \\
    &= \frac{(2m-2)(1-v^2)x^2}{a}.
\end{align}
Finally, \eqref{eq:appleton-Q-ODE} gives
\begin{align}
    z' &= b'Q' + bQ'' = -\left( f' + (2m-2)\frac{b'}{b} \right)z - \frac{2m}{b} Q(1-Q^2),
\end{align}
and hence, by \eqref{eq:comparison-basic-identities},
\begin{align}
    az' &= (2z+2xv-y)z - 2mx(1-x).
\end{align}
Dividing these identities by $ax' = 2xz$ gives
\begin{align}
    \dot{v} &= \frac{2m-yv-2x(1-v^2)}{2xz}, \label{eq:general-v-equation} \\
    \dot{y} &= \frac{(m-1)(1-v^2)x}{z}, \label{eq:general-y-equation} \\
    \dot{z} &= \frac{z}{x} + v - \frac{y}{2x} - \frac{m(1-x)}{z}. \label{eq:general-z-equation}
\end{align}
The parity conditions in \eqref{eq:appleton-boundary-conditions} imply that $v,y,z$ all extend smoothly as even functions of $s$ across $s=0$. Moreover, by \eqref{eq:bolt-expansion-a} and \eqref{eq:bolt-expansion-b}, we have $x = \kappa^2 s^2 + \O(s^4)$ as $s \to 0$. Thus, $v,y,z$ are smooth functions of $x$ down to $x=0$, and
\begin{align} \label{eq:comparison-bolt-values}
    v(0) = \frac{m}{\kappa}, \qquad y(0) = 2\kappa, \qquad z(0) = \kappa.
\end{align}

\begin{lemma} \label{lem:comparison-monotonicity}
    In the region $0<x<1$, we have
    \begin{align}
        0 < v < \frac{m}{\kappa} < 1, \qquad \dot{v} < 0, \qquad y > 2\kappa, \qquad \dot{y} > 0, \qquad z > 0, \qquad \lim_{x \uparrow 1} z = 0.
    \end{align}
\end{lemma}
\begin{proof}
    Lemma \ref{lem:appleton-qualitative-properties} implies $v>0$, $z>0$. Since $v(0) = m/\kappa < 1$, \eqref{eq:general-y-equation} shows that $y$ is strictly increasing whenever $v < 1$. At a first contact with $v=1$, we would have $y > y(0) = 2\kappa > 2m$, and so $\dot{v} = (2m-yv)/2xz < 0$, contradicting first contact from below. Therefore, $v<1$ throughout $x \in (0,1)$, and thus $\dot{y} > 0$ and $y > y(0) = 2\kappa$ on that region.

    We next prove $\dot{v} < 0$. By \eqref{eq:general-y-equation}, we have $\dot{y}(0) = 0$, and \eqref{eq:general-v-equation} gives
    \begin{align} \label{eq:2xz-dotv}
        2xz\dot{v} = 2m-yv-2x(1-v^2).
    \end{align}
    Differentiating this at $x=0$ and using \eqref{eq:comparison-bolt-values} gives
    \begin{align}
        2\kappa \dot{v}(0) = -2\kappa \dot{v}(0) - 2\left( 1 - \frac{m^2}{\kappa^2} \right).
    \end{align}
    Therefore,
    \begin{align} \label{eq:vdot-at-0}
        \dot{v}(0) = -\frac{1}{2\kappa} \left( 1 - \frac{m^2}{\kappa^2} \right) < 0.
    \end{align}
    Meanwhile, differentiating \eqref{eq:2xz-dotv} at a hypothetical first contact with $\dot{v} = 0$ would give
    \begin{align}
        2xz\ddot{v} = -\dot{y}v - 2(1-v^2) < 0,
    \end{align}
    hence $\ddot{v} < 0$. By \eqref{eq:vdot-at-0}, this contradicts first contact. Thus, in the interior we have $\dot{v} < 0$ and $v < v(0) = m/\kappa$.
    Finally, Lemma \ref{lem:appleton-qualitative-properties} gives, for $x \in (0,1)$,
    \begin{align}
        0 < z \leq \kappa e^{-f} b^{-2m+2}, \qquad b \geq b(0) = 1, \qquad f', f'' > 0.
    \end{align}
    The last inequalities imply $f(s) \to \infty$ as $s \to \infty$. Thus, $z \to 0$ as $s \to \infty$, equivalently as $x \uparrow 1$.
\end{proof}

\subsection{Proof of uniqueness} \label{subsec:uniqueness-proof}

\begin{proof}[Proof of Theorem \ref{thm:appleton-ODE-uniqueness}]
    Fix $m \geq 2$ and $\kappa > m$, and let $(f_i,a_i,b_i)$, $i=1,2$ be two triples satisfying the assumptions of the theorem. We will obtain a contradiction by assuming, without loss of generality, that
    \begin{align}
        F := f_2''(0) - f_1''(0) > 0.
    \end{align}
    Using these solutions, form the variables $x_i, v_i, y_i, z_i$ for $i=1,2$ as in \eqref{eq:x-uniqueness} and \eqref{eq:appleton-comparison-variables}. Using the change of variable $s \mapsto x_i$, we will regard $v_i,y_i,z_i$ as functions on the common interval $0<x<1$, and compare their profiles. Set
    \begin{align}
        V = v_1 - v_2, \qquad Y = y_2 - y_1, \qquad Z = z_1 - z_2, \qquad \alpha = 1 - \frac{m^2}{\kappa^2} > 0.
    \end{align}
    Note the reversed subtraction order in the definition of $Y$.

    \begin{lemma} \label{lem:short-time-VYZ-positivity}
        We have $V(x),Y(x),Z(x) > 0$ for all sufficiently small $x>0$. Moreover,
        \begin{align} \label{eq:Z-by-x-limit}
            \lim_{x \downarrow 0} \frac{Z(x)}{x} = \frac{F}{2\kappa}.
        \end{align}
    \end{lemma}
    \begin{proof}
        We expand $V,Y,Z$ near $x=0$. Firstly, inverting $x = \kappa^2 s^2 + \O(s^4)$ yields $s^2 = \frac{x}{\kappa^2} + \O(x^2)$. Thus, by \eqref{eq:bolt-expansion-a},
        \begin{align}
            a_i' = \kappa - \frac{f_i''(0) + 2m(m-1)}{2\kappa} x + \O(x^2).
        \end{align}
        Since $v_i = m/\kappa + \O(x)$, the identity $z_i = a_i' - xv_i$ gives
        \begin{align} \label{eq:zi-expansion}
            z_i = \kappa - \frac{f_i''(0) + 2m^2}{2\kappa} x + \O(x^2).
        \end{align}
        Hence the difference is given by
        \begin{align} \label{eq:Z-expansion}
            Z = \frac{F}{2\kappa} x + \O(x^2),
        \end{align}
        which is positive for sufficiently small $x>0$. This expansion also implies \eqref{eq:Z-by-x-limit}.

        To expand the other two differences, use \eqref{eq:comparison-bolt-values} and \eqref{eq:vdot-at-0} to get
        \begin{align} \label{eq:vi-expansion}
            v_i &= v_i(0) + \dot{v}_i(0)x + \O(x^2) = \frac{m}{\kappa} - \frac{\alpha}{2\kappa} x + \O(x^2).
        \end{align}
        Also, \eqref{eq:general-y-equation} gives $\dot{y}_i = \frac{(m-1)\alpha}{\kappa} x + \O(x^2)$. Integrating and using $y_i(0) = 2\kappa$ gives
        \begin{align} \label{eq:yi-expansion}
            y_i = 2\kappa + \frac{(m-1)\alpha}{2\kappa} x^2 + \O(x^3).
        \end{align}
        By \eqref{eq:vi-expansion} and \eqref{eq:yi-expansion}, we have
        \begin{align}
            V = h_2 x^2 + \O(x^3), \qquad Y = e_3 x^3 + \O(x^4)
        \end{align}
        for some constants $h_2$ and $e_3$. We now compute these constants. By \eqref{eq:2xz-dotv}, we have
        \begin{align} \label{eq:2xz-difference-eqn}
            2xz_1 \dot{V} + 2xZ \dot{v}_2 &= 2xz_1 \dot{v}_1 - 2xz_2 \dot{v}_2 = (-y_1 + 2x(v_1 + v_2)) V + Y v_2.
        \end{align}
        Plugging in the previously derived expansions and equating the coefficients of $x^2$, we get
        \begin{align}
            4\kappa h_2 - \frac{\alpha F}{2\kappa^2} = -2\kappa h_2, \qquad h_2 = \frac{\alpha F}{12\kappa^3}.
        \end{align}
        Similarly, taking differences of \eqref{eq:general-y-equation} yields
        \begin{align} \label{eq:Y-difference}
            \dot{Y} = (m-1)x \left( \frac{(v_1+v_2)V}{z_1} + \frac{(1-v_2^2)Z}{z_1 z_2} \right),
        \end{align}
        and equating coefficients of $x^2$ gives
        \begin{align}
            3e_3 = \frac{(m-1)\alpha F}{2\kappa^3}, \qquad e_3 = \frac{(m-1)\alpha F}{6\kappa^3}.
        \end{align}
        Thus
        \begin{align}
            V = \frac{\alpha F}{12\kappa^3} x^2 + \O(x^3), \qquad Y = \frac{(m-1)\alpha F}{6\kappa^3} x^3 + \O(x^4).
        \end{align}
        Both expressions are positive for sufficiently small $x>0$.
    \end{proof}

    We return to the proof of Theorem \ref{thm:appleton-ODE-uniqueness}. By Lemma \ref{lem:short-time-VYZ-positivity}, the differences $V,Y,Z$ are positive for sufficiently small $x > 0$. We claim that they remain positive throughout $(0,1)$.

    Suppose otherwise. By continuity, there is a first point $x_* \in (0,1)$ at which one of these differences vanishes. At that point, all three differences are nonnegative and at least one is zero. Any difference that is zero must have nonpositive derivative there.

    Taking differences in \eqref{eq:general-z-equation} gives
    \begin{align}
        \dot{Z} &= \frac{Z}{x} + V + \frac{Y}{2x} + \frac{m(1-x)Z}{z_1 z_2}.
    \end{align}
    Thus, on the interval $(0,x_*]$,
    \begin{align} \label{eq:Z-by-x-monotonicity}
        \frac{d}{dx} \left( \frac{Z}{x} \right) &= \frac{V}{x} + \frac{Y}{2x^2} + \frac{m(1-x)Z}{x z_1 z_2} \geq 0.
    \end{align}
    Using \eqref{eq:Z-by-x-limit}, we obtain
    \begin{align}
        Z(x_*) \geq \frac{F}{2\kappa} x_* > 0.
    \end{align}
    So $Z$ cannot be the vanishing difference at $x_*$.

    If $Y(x_*) = 0$, then \eqref{eq:Y-difference} gives $\dot{Y}(x_*) > 0$: the first term in brackets is nonnegative, the second term is strictly positive, and $(m-1)x > 0$. This is again a contradiction, so $Y$ cannot be the vanishing difference at $x_*$.

    Finally, if $V(x_*) = 0$, write $v_1(x_*) = v_2(x_*) = v_* > 0$. Evaluating \eqref{eq:2xz-difference-eqn} at $x_*$ yields
    \begin{align}
        2x_* z_1(x_*) \dot{V}(x_*) = Y(x_*) v_* - 2x_* Z(x_*) \dot{v}_2(x_*) > 0,
    \end{align}
    since $v_*, Y(x_*), Z(x_*) > 0$, and $\dot{v}_2(x_*) < 0$ by Lemma \ref{lem:comparison-monotonicity}. Thus $\dot{V}(x_*) > 0$, contradiction.

    We have therefore shown that $V,Y,Z > 0$ throughout $(0,1)$. Then \eqref{eq:Z-by-x-monotonicity} holds on the whole interval. Combining with \eqref{eq:Z-by-x-limit} gives
    \begin{align}
        z_1(x) - z_2(x) = Z(x) \geq \frac{F}{2\kappa} x > 0
    \end{align}
    for all $x \in (0,1)$. But this contradicts Lemma \ref{lem:comparison-monotonicity}, which asserts that $z_1, z_2 \to 0$ as $x \uparrow 1$. Hence $f_1''(0) = f_2''(0)$, proving Theorem \ref{thm:appleton-ODE-uniqueness}.
\end{proof}

\section{Further preparations for stability} \label{sec:further-preps-stability}

For the remainder of this paper, fix $k \geq 3$ and let $(M,g,f)$ be the four-dimensional Appleton soliton on $\O_{\CP^1}(-k)$. Thus, the triple $(f,a,b)$ describing the metric and potential via \eqref{eq:appleton-ansatz} satisfies the soliton equations \eqref{eq:appleton-soliton-f}--\eqref{eq:appleton-soliton-b} and the boundary conditions \eqref{eq:appleton-boundary-conditions}, with $m=2$ and $\kappa = k$.

\subsection{ODE estimates on four-dimensional Appleton solitons} \label{subsec:4d-background-estimates}

We reuse the variables defined in \S\ref{subsec:new-variables}, taking $m=2$ and $\kappa = k$. Thus
\begin{align}
    Q = \frac{a}{b}, \qquad x = Q^2, \qquad v = \frac{b'}{Q}, \qquad y = af' + 2a' + 2\frac{ab'}{b}, \qquad z = bQ'.
\end{align}
By Lemmas \ref{lem:appleton-qualitative-properties} and \ref{lem:comparison-monotonicity}, for any $s>0$ we have
\begin{align} \label{eq:4d-basic-bounds}
    0 < x < 1, \qquad 0 < v < \frac{2}{k} < 1, \qquad y > 2k, \qquad z > 0.
\end{align}
The equations derived in \S\ref{subsec:new-variables} give
\begin{align}
    ax' &= 2zx, \label{eq:4d-background-x} \\
    av' &= 4-yv-2x(1-v^2), \label{eq:4d-background-v} \\
    ay' &= 2x^2(1-v^2), \label{eq:4d-background-y} \\
    az' &= (2z+2xv-y)z - 4x(1-x). \label{eq:4d-background-z}
\end{align}

We first record the concavity of $a$ and some consequences for the growth of $a$ and $b$.
\begin{lemma} \label{lem:4d-concavity}
    For every $s>0$,
    \begin{align}
        a'' < 0, \qquad 0 < a' < k, \qquad a(s) \leq ks.
    \end{align}
    Moreover, there exists $c_* > 0$ such that as $s \to \infty$, we have
    \begin{align}
        a(s),b(s) \sim c_* \sqrt{s}, \qquad y(s) \to \infty.
    \end{align}
\end{lemma}
\begin{proof}
    The expansion \eqref{eq:bolt-expansion-a} gives $a'' < 0$ for sufficiently small $s>0$. Suppose $a'' = 0$ at some $s>0$. At this point, \eqref{eq:appleton-soliton-a} gives
    \begin{align} \label{eq:app0}
        a'f' + \frac{2a'b'}{b} = \frac{2a^3}{b^4}.
    \end{align}
    Lemma \ref{lem:appleton-qualitative-properties} gives $a',f' > 0$, so
    \begin{align}
        \frac{a'b'}{ab} < \frac{a^2}{b^4}.
    \end{align}
    Moreover, $Q' > 0$ implies $a'b > ab'$, so
    \begin{align}
        r := \frac{ab'}{a'b} \in (0,1).
    \end{align}
    Solving the equation \eqref{eq:app0} for $f'$ and substituting into \eqref{eq:appleton-soliton-b}, we obtain
    \begin{align}
        \frac{b''}{b} &= \frac{4}{b^2} - \frac{2a^2}{b^4} - \frac{a'b'}{ab} - \frac{(b')^2}{b^2} - \frac{b'}{b} f' \\
        &= \frac{4}{b^2} - \frac{2a^2}{b^4} (1+r) - \frac{a'b'}{ab} (1-r) \\
        &> \frac{4}{b^2} - (3+r) \frac{a^2}{b^4} = \frac{4-(3+r)x}{b^2} > 0.
    \end{align}
    However, \eqref{eq:appleton-soliton-f} and Lemma \ref{lem:appleton-qualitative-properties} give $b''/b = -f''/2 < 0$ at the same point. This contradiction shows that $a''$ cannot vanish, and hence $a''<0$ everywhere. Since $a'(0) = k$ and $a' > 0$, it follows that $0<a'<k$. Integrating gives $a(s) \leq ks$.

    Finally, \cite{appleton}*{Corollary 6.7} proves $a,b \sim c_*\sqrt{s}$ as $s \to \infty$ for some $c_* > 0$. Also, $f''>0$ implies $f'(s) \geq f'(1) > 0$ for $s \geq 1$. Since $y = af' + 2a' + 2ab'/b \geq af'$, we also have $y \to \infty$.
\end{proof}

We now introduce two further quantities:
\begin{align} \label{eq:Z-epsilon-defn}
    Z = \frac{bf'}{Q} + 2v, \qquad \xi = 4-x+xv^2-vy.
\end{align}
From now on, $Z$ denotes the quantity defined here, instead of the one introduced in \S\ref{subsec:uniqueness-proof}. 
These definitions give
\begin{align} \label{eq:background-algebraic-identities}
    a' = z+xv, \qquad bf' = \sqrt{x}(Z-2v), \qquad y = 2z+2xv+xZ.
\end{align}
In particular, $Z-2v > 0$. The following estimates will be used in the proof of stability.
\begin{lemma} \label{lem:sharp-background-bounds}
    For every $s>0$,
    \begin{align} \label{eq:sharp-background-bounds}
        z > k(1-x), \qquad vZ > 2, \qquad 0<\xi< x(1-v^2)^2 < 1.
    \end{align}
\end{lemma}
\begin{proof}
    We first prove $z > k(1-x)$. Set $u = z/(1-x) > 0$. By \eqref{eq:4d-background-x} and \eqref{eq:4d-background-z},
    \begin{align}
        au' = u(2u+2xv-y) - 4x.
    \end{align}
    At any point $s>0$ where $u=k$, the right-hand side equals $k(2k-y) + 2x(kv-2)$, which is negative by \eqref{eq:4d-basic-bounds}.
    Thus, if $u \leq k$ anywhere, it remains below $k$ thereafter. Since $y \to \infty$ and $x \to 1$ as $s \to \infty$, this would imply, for all sufficiently large $s$,
    \begin{align}
        au' \leq u\left(2k + \frac{4}{k} - y\right) - 4x \leq -2.
    \end{align}
    Using $a(s) \leq ks$ gives $u' \leq -2/(ks)$. Integrating this out to infinity shows that $u < 0$ eventually, a contradiction. This proves that $u > k$ and hence $z > k(1-x)$.

    Next, differentiating the definition of $\xi$ using \eqref{eq:4d-background-x}--\eqref{eq:4d-background-y}, and applying the last identity in \eqref{eq:background-algebraic-identities}, we get after simplification
    \begin{align} \label{eq:epsilon-evolution}
        a\xi' &= (2xv-y)\xi + x^2(1-v^2)(Z-2v).
    \end{align}
    Since $y - 2xv = af' + 2a'$ and $Z - 2v = bf'/Q$, this is equivalent to
    \begin{align}
        (a^2 e^f \xi)' = a^2 e^f x(1-v^2) f'.
    \end{align}
    Recall the initial values $v(0) = 2/k$ and $y(0) = 2k$ from \eqref{eq:comparison-bolt-values}, as well as $x(0) = 0$. Thus $\xi(0) = 0$, so for any $s > 0$,
    \begin{align} \label{eq:epsilon-integral}
        \xi(s) = \frac{1}{a(s)^2 e^{f(s)}} \int_0^s a(t)^2 e^{f(t)} x(t) (1-v(t)^2) f'(t) \, dt > 0.
    \end{align}
    This proves the claimed lower bound for $\xi$. To prove the upper bound, put $d = x(1-v^2)^2 - \xi$. The bolt expansions \eqref{eq:vi-expansion} and \eqref{eq:yi-expansion} imply $\xi = \O(x^2)$. Since $v = 2/k + \O(x)$, we have
    \begin{align}
        d = x \left( 1 - \frac{4}{k^2} \right)^2 + \O(x^2) > 0
    \end{align}
    for sufficiently small $x>0$. At a hypothetical first contact with $d=0$ we have, using \eqref{eq:4d-background-v} and the definition of $\xi$,
    \begin{align}
        \xi &= x(1-v^2)^2, \label{eq:epsilon-value} \\
        av' &= \xi - x(1-v^2) = -xv^2 (1-v^2), \label{eq:avp-value} \\
        yv &= 4-x+xv^2-\xi = 4-2x+3xv^2 - xv^4. \label{eq:yv-value}
    \end{align}
    Differentiating $d$ then using \eqref{eq:epsilon-evolution}, \eqref{eq:epsilon-value}, and \eqref{eq:avp-value}, we obtain at the same point
    \begin{align}
        ad' &= 2xz(1-v^2)^2 - 4xv(1-v^2) av' - a\xi' \\
        &= 2xz(1-v^2)^2 + 4x^2 v^3(1-v^2)^2 - x(2xv-y)(1-v^2)^2 - x^2(1-v^2)(Z-2v).
    \end{align}
    Substituting the identity $xZ = y-2z-2xv$ as well as \eqref{eq:yv-value}, we obtain after simplification,
    \begin{align}
        ad' &= x(1-v^2) \left[ (4-2v^2)(z+xv) - 4v + xv(5v^2 - 3v^4) \right].
    \end{align}
    Using $z > k(1-x)$, $k \geq 3$, and $0<v<1$, the term in square brackets is strictly greater than
    \begin{align}
        &(4-2v^2)(k(1-x)+xv)-4v + xv(5v^2-3v^4) \\
        &\qquad = (1-x) \left[ 2k(2-v^2) - 4v \right] + 3xv^3 (1-v^2) > 0.
    \end{align}
    Thus, $d' > 0$ at its first zero, a contradiction. This proves that $d>0$ whenever $s>0$, so $0 < \xi < x(1-v^2)^2 < 1$.

    It remains to prove $vZ > 2$. Differentiating $y=2z + 2xv + xZ$ and using \eqref{eq:4d-background-x}--\eqref{eq:4d-background-z} yields
    \begin{align}
        xaZ' &= ay' - 2az' - 2xav' - (2v+Z)ax' \\
        &= 2x^2(1-v^2) + 8x(1-x) - 2x(4-yv-2x(1-v^2)) - 4vxz.
    \end{align}
    Dividing by $x$, substituting $y=2z + 2xv + xZ$ again, then simplifying, we get
    \begin{align}
        aZ' = 2x(vZ-1-v^2).
    \end{align}
    Hence, setting $F = vZ - 2$, we obtain
    \begin{align}
        aF' &= Zav' + vaZ' = Z (\xi - x + xv^2) + 2vx(vZ - 1 - v^2) \\
        &= Z\xi + xZ(3v^2 - 1) - 2vx(1+v^2).
    \end{align}
    Substituting $3v^2 - 1 = v^2(1+v^2) - (1-v^2)^2$, we get
    \begin{align} \label{eq:Fprime}
        aF' = vx(1+v^2)F - Zd.
    \end{align}
    Since $Z,d > 0$, if $F$ were nonpositive at a point, then \eqref{eq:Fprime} gives $F'<0$ there. Thus, $F$ would be strictly negative at some later point $s_1$. Then for all $s \geq s_1$,
    \begin{align}
        F(s) \leq F(s_1) \exp\left\{ \int_{s_1}^s \frac{vx}{a}(1+v^2) \, dt \right\}.
    \end{align}
    But $vx/a = b'/b$, so the exponential factor is at least $b(s)/b(s_1)$, which goes to $\infty$ as $s \to \infty$. Thus $F(s) \to -\infty$, contradicting $F = vZ - 2 > -2$. Therefore $F > 0$, so $vZ > 2$.
\end{proof}

\subsection{Local frames, connection, and curvature} \label{subsec:local-frame}

We introduce a local orthonormal frame on $M^\circ = M \setminus B$ and record the connection and curvature in this frame. Identify $SU(2)$ with $\Sph^3 \subset \C^2$ by
\begin{align} \label{eq:SU2-S3-identification}
    SU(2) \ni \theta \mapsto \theta \cdot (0,1)^T \in \Sph^3.
\end{align}
Define
\begin{align} \label{eq:T-matrices}
    T_1 = \begin{pmatrix}
        -i & 0 \\ 0 & i
    \end{pmatrix},
    \qquad
    T_2 = \begin{pmatrix}
        0 & 1 \\ -1 & 0
    \end{pmatrix},
    \qquad
    T_3 = \begin{pmatrix}
        0 & -i \\ -i & 0
    \end{pmatrix}.
\end{align}
These form a basis of $\mathfrak{su}(2)$ with
\begin{align} \label{eq:T-brackets}
    [T_1, T_2] = 2T_3, \qquad [T_2, T_3] = 2T_1, \qquad [T_3, T_1] = 2T_2.
\end{align}
Let $X_1, X_2, X_3$ be the corresponding left-invariant vector fields on $SU(2)$, and let $\sigma_1, \sigma_2, \sigma_3$ be the dual coframe. Under the identification \eqref{eq:SU2-S3-identification}, the curve $\theta \exp(t T_1)$ corresponds to $e^{it} \theta \cdot (0,1)^T$. Thus, $X_1$ is the unit Hopf vector field, and
\begin{align}
    \sigma_1 = \sigma, \qquad \sigma_2^2 + \sigma_3^2 = \pi^*g_{\mathrm{FS}}.
\end{align}
On the cover $(0,\infty) \times SU(2)$ of $M^\circ$, the lifted soliton metric is
\begin{align}
    \tilde{g} = ds^2 + a(s)^2 \sigma_1^2 + b(s)^2 (\sigma_2^2 + \sigma_3^2).
\end{align}
A global orthonormal frame on the cover is therefore given by
\begin{align} \label{eq:appleton-E-frame}
    E_0 = \del_s, \qquad E_1 = \frac{X_1}{a}, \qquad E_2 = \frac{X_2}{b}, \qquad E_3 = \frac{X_3}{b},
\end{align}
with dual coframe
\begin{align} \label{eq:appleton-epsilon-frame}
    \epsilon^0 = ds, \qquad \epsilon^1 = a \sigma_1, \qquad \epsilon^2 = b\sigma_2, \qquad \epsilon^3 = b\sigma_3.
\end{align}
The scalar $\Z_k$-action on $\Sph^3 \cong SU(2)$ leaves $E_0$ and $E_1$ invariant, so $E_0,E_1,\epsilon^0,\epsilon^1$ descend globally to $M^\circ$. Meanwhile, the $\Z_k$-action rotates the $(E_2,E_3)$-plane. Thus, $E_2,E_3,\epsilon^2,\epsilon^3$ only descend locally, and we understand the frame \eqref{eq:appleton-E-frame} and \eqref{eq:appleton-epsilon-frame} to be defined locally on $M^\circ$. On the other hand, $\Span\{E_2,E_3\}$ and $\Span\{\epsilon^2,\epsilon^3\}$ are globally defined tangent and cotangent subbundles on $M^\circ$, and the tensor $\sigma_2^2 + \sigma_3^2 = \pi^*g_{\mathrm{FS}}$ is also globally defined on $M^\circ$.

Define the quantities
\begin{align} \label{eq:connection-coeffs-defn}
    A = \frac{a'}{a}, \qquad B = \frac{b'}{b}, \qquad C = \frac{2}{a} - \frac{a}{b^2}, \qquad D = \frac{a}{b^2}.
\end{align}
Left-invariance implies that $\{X_\alpha\}_{\alpha=1}^3$ also satisfy the relations \eqref{eq:T-brackets}. We then calculate
\begin{gather}
    [E_0,E_1] = -AE_1, \qquad [E_0,E_\mu] = -B E_\mu \quad (\mu=2,3), \\
    [E_1,E_2] = \frac{2}{a} E_3, \qquad [E_2,E_3] = 2D E_1, \qquad [E_3,E_1] = \frac{2}{a} E_2.
\end{gather}
The Koszul formula then gives the following connection coefficients:
\begin{align}\label{eq:frame-connection}
    \begin{array}{c|cccc}
        &E_0&E_1&E_2&E_3\\ \hline
        \nabla_{E_0}
        &0&0&0&0\\[2mm]
        \nabla_{E_1}
        &AE_1&-AE_0&CE_3&-CE_2\\[2mm]
        \nabla_{E_2}
        &BE_2&-DE_3&-BE_0&DE_1\\[2mm]
        \nabla_{E_3}
        &BE_3&DE_2&-DE_1&-BE_0
    \end{array}.
\end{align}
In particular, this frame is parallel along $\del_s$. By metric compatibility, the same table describes $\nabla_{E_i} \epsilon^j$ after replacing each $E_i$ on the right by $\epsilon^i$.

Using our curvature convention \eqref{eq:Rm-convention} and $\Rm_{ijkl} := \inner{\Rm(E_i,E_j)E_k}{E_l}$, the nonzero coefficients of $\Rm$ up to symmetries are given by
\begin{equation} \label{eq:frame-curvature}
    \begin{aligned}
        \Rm_{0110} &= -\frac{a''}{a} =: K_1, \\
        \Rm_{0220} = \Rm_{0330} &= -\frac{b''}{b} =: K_2, \\
        \Rm_{1221} = \Rm_{1331} &= D^2 - AB =: H_{12}, \\
        \Rm_{2332} &= \frac{4}{b^2} - 3D^2 - B^2 =: H_{23}, \\
        \Rm_{0231} = \Rm_{0312} &= D(A-B) =: M_2, \\
        \Rm_{0123} &= -2D(A-B) = -2M_2.
    \end{aligned}
\end{equation}
Thus, the Ricci and scalar curvatures are given by
\begin{align}
    \Ric &= (K_1 + 2K_2)(\epsilon^0)^2 + (K_1 + 2H_{12})(\epsilon^1)^2 \\
    &\quad + (K_2 + H_{12} + H_{23}) ((\epsilon^2)^2 + (\epsilon^3)^2), \label{eq:appleton-frame-Ricci} \\
    R &= 2(K_1 + 2K_2 + 2H_{12} + H_{23}). \label{eq:appleton-frame-scalar}
\end{align}
Here and below, for a $1$-form $\alpha$ we write $\alpha^2 = \alpha \otimes \alpha$. We also compute
\begin{align} \label{eq:hess-f-calculation}
    \nabla^2 f = f''(\epsilon^0)^2 + Af'(\epsilon^1)^2 + Bf'((\epsilon^2)^2 + (\epsilon^3)^2),
\end{align}
so that Lemma \ref{lem:appleton-qualitative-properties} and the soliton equation give $\Ric = \nabla^2 f > 0$ on $M^\circ$. By smoothness,
\begin{align}
    \Ric \geq 0 \quad \text{on } M.
\end{align}
Tracing $\Ric = \nabla^2 f$ also gives
\begin{align} \label{eq:R-intermsof-f}
    R = f'' + (A+2B)f'.
\end{align}
The soliton equations \eqref{eq:appleton-soliton-f}--\eqref{eq:appleton-soliton-b} and \eqref{eq:background-algebraic-identities} express the curvature coefficients in terms of the variables $x,v,z,Z$ from \S\ref{subsec:new-variables} and \eqref{eq:Z-epsilon-defn}:
\begin{equation}
\label{eq:scaled-curvature-coefficients}
\begin{aligned}
    b^2 K_1 &= (z+xv)Z - 2x, \\
    b^2 K_2 &= -4+2x+vz+vxZ, \\
    b^2 H_{12} &= x-vz-xv^2, \\
    b^2 H_{23} &= 4-3x-xv^2, \\
    b^2 M_2 &= z.
\end{aligned}
\end{equation}

We also record the linear decay of curvature:
\begin{lemma} \label{lem:curvature-linear-decay}
    There exists $C_{\Rm} > 0$ such that $b^2|{\Rm}| \leq C_{\Rm}$ on $M$.
\end{lemma}
\begin{proof}
    By Lemma \ref{lem:4d-concavity}, $a'' < 0$ and $a(0)=0$, which imply $sa'(s) \leq a(s)$. Since $Q'>0$, we therefore have $0 < B < A = a'/a \leq 1/s$. Also, $a,b \sim c_*\sqrt{s}$ give $b^{-2}=\O(s^{-1})$ and $D=\O(s^{-1/2})$. It follows from \eqref{eq:frame-curvature} that
    \begin{align}
        H_{12},H_{23},M_2=\O(s^{-1}).
    \end{align}
    Next, equating \eqref{eq:appleton-frame-Ricci} and \eqref{eq:hess-f-calculation}, we get
    \begin{align}
        K_1 &= Af'-2H_{12}, \\
        K_2 &= Bf'-H_{12}-H_{23}.
    \end{align}
    Lemma \ref{lem:appleton-qualitative-properties} gives $R>0$, so the standard identity $R+(f')^2=\text{const.}$ for complete steady gradient Ricci solitons implies that $f'$ is bounded. Thus, $K_1,K_2=\O(s^{-1})$ as well. By \eqref{eq:frame-curvature}, every curvature component is therefore $\O(s^{-1})$, proving that $b^2|{\Rm}|$ is bounded.
\end{proof}

\subsection{Additional conventions; reductions of the stability inequality} \label{subsec:stability-reductions}

Let us specify some additional notation and conventions used throughout the rest of this paper. From here onwards, all tensors on $M$ may be complex-valued unless otherwise stated. We use the induced Hermitian inner products $\inner{\cdot}{\cdot}$ on complex tensor bundles, so that they are linear in their first arguments and conjugate-linear in the second. Our convention for the Laplacian and weighted Laplacian on tensors is
\begin{align}
    \Delta = \tr_g \nabla^2, \qquad \Delta_f = \Delta + \nabla_{\nabla f}.
\end{align}
For $1$-forms $\alpha,\beta$, we write
\begin{gather}
    \alpha \odot \beta = \frac{1}{2}(\alpha \otimes \beta + \beta \otimes \alpha), \label{eq:symm-prod} \\
    \alpha \wedge \beta = \alpha \otimes \beta - \beta \otimes \alpha. \label{eq:wedge-prod}
\end{gather}
We abbreviate $\epsilon^{ij} = \epsilon^i \wedge \epsilon^j$. Interior multiplication on differential forms is given by
\begin{align}
    \iota_V \alpha = \alpha(V, \ldots)
\end{align}
for any tangent vector $V$. On differential $p$-forms, we use the Hermitian inner product
\begin{align} \label{eq:form-inner-product}
    \inner{\alpha}{\beta} = \frac{1}{p!} \sum_{i_1,\ldots,i_p} \alpha_{i_1 \cdots i_p} \overline{\beta_{i_1 \cdots i_p}},
\end{align}
where the sum is over an orthonormal frame. This differs from the inner product on general covariant $p$-tensors by a factor of $1/p!$.

The quantity $\calE[h]$ from \eqref{eq:strict-linear-stab} makes sense for compactly supported (complex) symmetric $2$-tensors $h$. Then for compactly supported real symmetric $2$-tensors $h_1$ and $h_2$, we have $\calE[h_1 + ih_2] = \calE[h_1] + \calE[h_2]$ and $|h_1 + ih_2|^2 = |h_1|^2 + |h_2|^2$. Consequently, strict linear stability of $(M,g,f)$ is equivalent to the existence of $\mu > 0$ such that for any smooth compactly supported symmetric $2$-tensor $h$ on $M$,
\begin{align}
    \calE[h] \geq \mu \int_M |h|^2 e^f.
\end{align}

Let $d\omega$ be the volume measure of the unit round quotient metric $\sigma^2 + \pi^*g_{\mathrm{FS}}$ on $L_{2,k}$. The metric ansatz \eqref{eq:appleton-ansatz} gives
\begin{align}
    \dvol_g = ab^2 \, ds \, d\omega.
\end{align}
Set $w = e^f ab^2$, so that $w$ is a function of $s$ and
\begin{align} \label{eq:dvolg-decomp}
    e^f \dvol_g = w \, ds \, d\omega.
\end{align}
Recall that $R+(f')^2 = \text{const.} =: c > 0$. We now prove a radial Poincar\'e inequality:

\begin{lemma} \label{lem:poincare-E0}
    For every smooth compactly supported symmetric $2$-tensor $h$ on $M$,
    \begin{align}
        \int_M |\nabla_{\del_s} h|^2 e^f \geq \frac{c}{4} \int_M |h|^2 e^f.
    \end{align}
\end{lemma}
\begin{proof}
    We have
    \begin{align} \label{eq:sq-completion-poincare}
        \int_M \left| \nabla_{\del_s} h + \frac{f'}{2} h \right|^2 e^f &= \int_M \left( |\nabla_{\del_s} h|^2 + \frac{(f')^2}{4} |h|^2 + f' \Re\inner{\nabla_{\del_s}h}{h} \right) e^f.
    \end{align}
    Using $f'\del_s = \nabla f$, the identity $\Delta_f f = R + (f')^2 = c$, and integrating by parts, we have
    \begin{align}
        \int_M f' \Re\inner{\nabla_{\del_s} h}{h} e^f &= \frac{1}{2} \int_M \inner{\nabla f}{\nabla |h|^2} e^f = -\frac{1}{2} \int_M (\Delta_f f) |h|^2 e^f = -\frac{c}{2} \int_M |h|^2 e^f.
    \end{align}
    Inserting this into \eqref{eq:sq-completion-poincare} and using again that $R+(f')^2 = c$, we obtain
    \begin{align}
        0 \leq \int_M \left| \nabla_{\del_s} h + \frac{f'}{2} h \right|^2 e^f = \int_M \left( |\nabla_{\del_s} h|^2 - \frac{R+c}{4} |h|^2 \right) e^f.
    \end{align}
    Since $R > 0$, the conclusion follows.
\end{proof}

We now show that strict linear stability follows from an ostensibly weaker estimate for symmetric 2-tensors supported in $M^\circ$, i.e. supported away from the zero section of $\O_{\CP^1}(-k)$.
\begin{proposition} \label{prop:stability-reduction}
    Suppose that there exists $\mu > 0$ such that for every smooth symmetric $2$-tensor compactly supported in $M^\circ$,
    \begin{align} \label{eq:b-weighted-stability}
        \calE[h] \geq \mu \int_M \frac{|h|^2}{b^2} e^f.
    \end{align}
    Then $(M,g,f)$ is strictly linearly stable.
\end{proposition}
\begin{proof}
    Let $h$ be a smooth compactly supported symmetric $2$-tensor on $M$. For $0 < \epsilon < 1/2$, choose a smooth radial cutoff function $\chi_\epsilon = \chi_\epsilon(s)$ such that
    \begin{align}
        0 \leq \chi_\epsilon \leq 1, \qquad \chi_\epsilon = 0 \text{ on } \{s \leq \epsilon^2\}, \qquad \chi_\epsilon = 1 \text{ on } \{s \geq \epsilon\}, \qquad |\chi_\epsilon'| \leq \frac{C}{s|{\log\epsilon}|},
    \end{align}
    where $C>0$ is independent of $\epsilon$. Applying \eqref{eq:b-weighted-stability} to $\chi_\epsilon h$, which is compactly supported in $M^\circ$, gives
    \begin{align} \label{eq:b-weighted-stability-cutoff}
        \calE[\chi_\epsilon h] \geq \mu \int_M \frac{|\chi_\epsilon h|^2}{b^2} e^f.
    \end{align}
    We have
    \begin{align}
        \int_M |(\chi_\epsilon - 1)h|^2 e^f \leq \int_{\{s < \epsilon\}} |h|^2 e^f \to 0 \quad \text{as } \epsilon \to 0.
    \end{align}
    Moreover, the expansions \eqref{eq:bolt-expansion-f}--\eqref{eq:bolt-expansion-b} give $w(s) = ks(1+\O(s^2))$, so
    \begin{align}
        \int_M |\nabla((\chi_\epsilon - 1) h)|^2 e^f &\lesssim \int_M (\chi_\epsilon - 1)^2 |\nabla h|^2 e^f + \int_M |\nabla\chi_\epsilon|^2 |h|^2 e^f \\
        &= \int_M (\chi_\epsilon-1)^2 |\nabla h|^2 e^f + \int_{L_{2,k}} \int_0^\infty w|\nabla\chi_\epsilon|^2 |h|^2 \, ds \, d\omega \\
        &\lesssim \int_{\{s < \epsilon\}} |\nabla h|^2 e^f + \frac{1}{|{\log\epsilon}|^2} \int_{L_{2,k}} \int_{\epsilon^2}^{\epsilon} \frac{ds}{s}\, d\omega \to 0 \quad \text{as } \epsilon \to 0.
    \end{align}
    Here we have used \eqref{eq:dvolg-decomp} to break up the integral into radial and cross-sectional parts. The above shows that $\chi_\epsilon h \to h$ in the weighted Sobolev space $H^1(e^f)$ as $\epsilon \to 0$. Since the curvature tensor is bounded, and $b$ is bounded from below (Lemma \ref{lem:appleton-qualitative-properties}), we can therefore take limits in \eqref{eq:b-weighted-stability-cutoff} to get
    \begin{align} \label{eq:b-weighted-stability-proof}
        \calE[h] \geq \mu \int_M \frac{|h|^2}{b^2} e^f.
    \end{align}
    Next, we apply Lemma \ref{lem:poincare-E0}, Lemma \ref{lem:curvature-linear-decay}, and \eqref{eq:b-weighted-stability-proof} to get
    \begin{align}
        \frac{c}{4} \int_M |h|^2 e^f &\leq \int_M |\nabla_{\del_s} h|^2 e^f \leq \int_M |\nabla h|^2 e^f \\
        &= \calE[h] + 2\int_M \inner{\Rm(h)}{h} e^f \\
        &\leq \calE[h] + 2C_{\Rm} \int_M \frac{|h|^2}{b^2} e^f \\
        &\leq (1+2C_{\Rm}\mu^{-1}) \calE[h].
    \end{align}
    Thus, we have shown that for all compactly supported symmetric 2-tensors $h$ on $M$,
    \begin{align}
        \calE[h] \geq \frac{c}{4(1+2C_{\Rm}\mu^{-1})} \int_M |h|^2 e^f.
    \end{align}
    Since the constant factor is positive and independent of $h$, it follows that $(M,g,f)$ is strictly linearly stable.
\end{proof}

\section{Separation of variables} \label{sec:appleton-separation}

In this section, we decompose tensor fields on the four-dimensional Appleton soliton $(M,g,f)$ into sectors corresponding to the irreducible representations of $U(2)$, and show that it suffices to establish the stability inequality in each sector.

\subsection{The group action and Peter--Weyl decomposition} \label{subsec:group-action-peter-weyl}

The standard $G=U(2)$-action on $L_{2,k} = \Sph^3/\Z_k$ is transitive, with isotropy subgroup at $o=[(0,1)]$ given by
\begin{align} \label{eq:GH-defns}
    H = \{ \diag(u,\zeta) : u \in U(1), \, \zeta^k = 1 \}.
\end{align}
Thus, there is a $G$-equivariant identification $L_{2,k} \cong G/H$. We write
\begin{align}
    \Sigma_s = \{s\} \times L_{2,k} \subset M^\circ, \qquad \Sigma_s \cong G/H.
\end{align}
The natural action of $G$ on $M^\circ$ is isometric for the Appleton metric \eqref{eq:appleton-ansatz}, preserves $f$, and has cohomogeneity one with orbits $\Sigma_s$. Moreover, the measure $d\omega$ from \eqref{eq:dvolg-decomp} is $G$-invariant.

For $\theta \in G$, let $\phi_\theta$ denote the corresponding isometry of $M^\circ$. The induced action on covariant tensors (i.e. sections of tensor powers of $T^*_{\C}M^\circ$) is given by inverse pullback:
\begin{align} \label{eq:inverse-pullback-action}
    \theta \cdot u = \phi_{\theta^{-1}}^* u.
\end{align}

Let $\scrE \to M^\circ$ be a complex subbundle of a covariant tensor bundle preserved by the action \eqref{eq:inverse-pullback-action}, equipped with the Hermitian metric induced by $g$. For each $s>0$, regard $\scrE_s := \scrE|_{\Sigma_s}$ as a bundle over $L_{2,k}$, and let $F_s$ be its fiber above $o$. Since $H$ fixes $o$, restricting the action \eqref{eq:inverse-pullback-action} to $H$ and $F_s$ gives a unitary representation $\rho_s: H \to U(F_s)$ called the \emph{isotropy representation}. Then the action of $G$ on $\scrE_s$ identifies it as a homogeneous vector bundle
\begin{align} \label{eq:E-homog-vb}
    G \times_{\rho_s} F_s \cong \scrE_s, \qquad [\theta,\eta] \mapsto \theta \cdot \eta.
\end{align}
Under this identification, a section $u$ of $\scrE_s$ corresponds to a function $U: G \to F_s$ satisfying
\begin{align} \label{eq:U-identification}
    u(\theta H) = [\theta, U(\theta)], \qquad U(\theta h) = \rho_s(h)^{-1} U(\theta).
\end{align}
Then the inverse pullback action of $G$ is given by $(\theta \cdot U)(\theta') = U(\theta^{-1} \theta')$. We note that this action is unitary on $L^2(\scrE_s, d\omega)$.

To state the decomposition of $L^2(\scrE_s, d\omega)$, we recall the representation theory of $G = U(2)$. The irreducible unitary representations of $G$ are indexed by pairs $(j,q)$ satisfying
\begin{align} \label{eq:jq-representation-labels}
    j \in \tfrac{1}{2}\Z_{\geq 0}, \qquad q \in \Z, \qquad q \equiv 2j \pmod 2.
\end{align}
Denote the corresponding representation by $(V_{j,q}, \pi_{j,q})$. Its restriction to $SU(2)$ is the spin-$j$ representation of dimension $2j+1$, while the central elements $zI \in G$ ($z \in U(1)$) act by multiplication by $z^q$. For each such pair $(j,q)$, define the space
\begin{align}
    \calM_{j,q}(F_s) &:= \Hom_H(V_{j,q}, F_s) \\
    &= \left\{ \sfA: V_{j,q} \to F_s \text{ linear} : \sfA \pi_{j,q}(h) = \rho_s(h) \sfA \text{ for all } h \in H \right\},
\end{align}
and equip $V_{j,q} \otimes \calM_{j,q}(F_s)$ with the $G$-action $\pi_{j,q} \otimes I$. The Peter--Weyl theorem for homogeneous vector bundles (e.g. \cite{wallach}*{Chapter 5}) gives the following decomposition.

\begin{proposition} \label{prop:sector-peter-weyl}
    For every $s>0$, there is an orthogonal decomposition
    \begin{align} \label{eq:bundle-peter-weyl}
        L^2(\scrE_s, d\omega) = \widehat{\bigoplus_{j,q}} \, \calH_{j,q,s},
    \end{align}
    where the sum runs over pairs satisfying \eqref{eq:jq-representation-labels}. Each summand is $G$-equivariantly isomorphic to $V_{j,q} \otimes \calM_{j,q}(F_s)$ via the map
    \begin{align} \label{eq:normalized-section-map}
        \scrT_{j,q,s}: V_{j,q} \otimes \calM_{j,q}(F_s) \to \calH_{j,q,s}, \qquad
        \scrT_{j,q,s}(v \otimes \sfA)(\theta H) = [\theta, \sfA(\pi_{j,q}(\theta^{-1})v)].
    \end{align}
    The orthogonal projection $\Pi_{j,q,s}: L^2(\scrE_s, d\omega) \to \calH_{j,q,s}$ is given by
    \begin{align} \label{eq:jq-orthog-projection}
        \Pi_{j,q,s} u = (2j+1) \int_G \overline{\chi_{j,q}(\theta)} (\theta \cdot u) \, d\theta,
    \end{align}
    where $\chi_{j,q}: G \to \C$ is the character of $\pi_{j,q}$, and $d\theta$ is normalized Haar measure on $G$.
\end{proposition}

The subspace $\calH_{j,q,s}$ in \eqref{eq:bundle-peter-weyl} is called the \emph{$(j,q)$-sector} of $\scrE_s$. We will also say that a smooth section $u \in C^\infty(\scrE)$ belongs to the \emph{$(j,q)$-sector} of $\scrE$ if $u|_{\Sigma_s}$ belongs to the $(j,q)$-sector of $\scrE_s$ for every $s>0$.

\begin{corollary} \label{cor:equivariant-operator-splitting}
    Let $\scrE,\scrE'$ be bundles as above, and let $P: C^\infty(\scrE) \to C^\infty(\scrE')$ be a $G$-equivariant linear differential operator. If $u$ belongs to the $(j,q)$-sector of $\scrE$, then $Pu$ belongs to the $(j,q)$-sector of $\scrE'$.
\end{corollary}
\begin{proof}
    A smooth section $u \in C^\infty(\scrE)$ belongs to the $(j,q)$-sector if and only if
    \begin{align} \label{eq:u-jq-sector}
        u = (2j+1) \int_G \overline{\chi_{j,q}(\theta)} (\theta \cdot u) \, d\theta.
    \end{align}
    Indeed, since the action of $G$ on $M^\circ$ fixes $s$, restricting \eqref{eq:u-jq-sector} to $\Sigma_s$ reads $u|_{\Sigma_s} = \Pi_{j,q,s} u|_{\Sigma_s}$ by \eqref{eq:jq-orthog-projection}. This holding for each $s>0$ is precisely the condition that $u$ belongs to the $(j,q)$-sector.
    
    Applying $P$ to both sides of \eqref{eq:u-jq-sector}, commuting $P$ and the integral, then applying $G$-equivariance, we get
    \begin{align}
        Pu = (2j+1) \int_G \overline{\chi_{j,q}(\theta)} (\theta \cdot Pu) \, d\theta,
    \end{align}
    which is precisely the condition that $Pu$ belongs to the $(j,q)$ sector of $\scrE'$.
\end{proof}

Finally, we relate tensors in $(j,0)$-sectors to tensors on the base of the Hopf fibration. For each $j \in \Z_{\geq 0}$, the center of $G$ acts trivially on $V_{j,0}$ and hence on $V_{j,0} \otimes \calM_{j,0}(F_s)$. The center of $G$ also induces the Hopf circle action on $L^2(\scrE_s, d\omega)$. Thus, by $G$-equivariance in Proposition \ref{prop:sector-peter-weyl}, any section in a $(j,0)$-sector of $\scrE_s$ is invariant under the Hopf circle action.

Concretely, let $\calH = \ker\sigma \subset T\Sigma_s$ be the horizontal distribution, and let $P_{\calH}: T\Sigma_s \to \calH$ be orthogonal projection with respect to the induced metric. If $u$ is a covariant $r$-tensor, define its horizontal part by
\begin{align}
    u_{\calH}(v_1, \ldots, v_r) = u(P_{\calH} v_1, \ldots, P_{\calH} v_r), \qquad v_1,\ldots,v_r \in T\Sigma_s.
\end{align}
Invariance under the Hopf circle action means the following:
\begin{lemma} \label{lem:hopf-invariance}
    Let $u$ be a section in a $(j,0)$-sector of $\scrE_s$, where $j \in \Z_{\geq 0}$. Then there is a unique covariant $r$-tensor $\bar{u}$ on $\CP^1$ such that $u_{\calH} = \pi^*\bar{u}$.
\end{lemma}

\subsection{Isotropy weights and examples} \label{subsec:isotropy-weights-examples}

Let $\scrE \to M^\circ$ be a bundle as in \S\ref{subsec:group-action-peter-weyl}, and define $F_s$ and $\calM_{j,q}(F_s)$ as before.
We now explain how to determine $\calM_{j,q}(F_s)$ from the isotropy representation, then apply this to the bundles needed below. Throughout this subsection, we fix $s>0$, write $F=F_s$ for the fiber of $\scrE_s \to \Sigma_s \cong L_{2,k}$ above $o$, and $\rho=\rho_s$ for the isotropy representation of $H$ on $F$.

For each $(j,q)$ satisfying \eqref{eq:jq-representation-labels}, the standard representation theory of $G$ allows us to find an orthonormal weight basis $\{v_\ell\}_{\ell=-j}^j$ of $V_{j,q}$ such that
\begin{align} \label{eq:representation-torus-weights}
    \pi_{j,q}(\diag(z_1,z_2)) v_\ell = z_1^{q/2+\ell} z_2^{q/2-\ell} v_\ell.
\end{align}
Define the character
\begin{align}
    \chi: H \to U(1), \qquad \chi(\diag(u,\zeta)) = u\zeta^{-1}.
\end{align}
Restricting the representation to $H$ gives
\begin{align} \label{eq:H-representation-on-V}
    \pi_{j,q}(h) v_\ell = \zeta^q \chi(h)^{q/2+\ell} v_\ell, \qquad h = \diag(u,\zeta).
\end{align}
Suppose the isotropy representation $\rho: H \to U(F)$ admits an orthogonal decomposition
\begin{align} \label{eq:fiber-weight-decomp}
    F = \bigoplus_{r \in \Z} F_r, \qquad \rho(h)|_{F_r} = \chi(h)^r I_{F_r},
\end{align}
with only finitely nonzero summands. This will hold for all the bundles considered below. Then the space $\calM_{j,q}(F)$ is given as follows.

\begin{lemma} \label{lem:weight-matching}
    Suppose $F$ has the decomposition \eqref{eq:fiber-weight-decomp}, and let $(j,q)$ satisfy \eqref{eq:jq-representation-labels}. If $q \not\equiv 0 \pmod k$, then $\calM_{j,q}(F) = 0$. Otherwise, evaluation on the weight basis defines a linear isomorphism
    \begin{align} \label{eq:weight-multiplicity-identification}
        \calM_{j,q}(F) \to \bigoplus_{\substack{r \in \Z \\ |r-q/2| \leq j}} F_r, \qquad \sfA  \mapsto \left( \sfA(v_{r-q/2}) \right)_r.
    \end{align}
\end{lemma}
\begin{proof}
    Let $\sfA \in \calM_{j,q}(F)$, and take a weight vector $v_\ell$, $-j \leq \ell \leq j$. For $r \in \Z$, the intertwining property of $\sfA$ implies that it can only map $v_\ell$ nontrivially into $F_r$ if their $H$-characters agree. By \eqref{eq:H-representation-on-V} and \eqref{eq:fiber-weight-decomp}, this condition reads
    \begin{align}
        u^{\frac{q}{2}+\ell} \zeta^{\frac{q}{2}-\ell} = u^r \zeta^{-r}
    \end{align}
    for every $\diag(u,\zeta) \in H$. Since $u$ can be arbitrary in $U(1)$, this condition implies $r = \frac{q}{2} + \ell$. This in turn implies $\zeta^q = 1$ for every $k$-th root of unity $\zeta$, equivalently $q \equiv 0 \pmod k$.

    Thus, if $q \not\equiv 0 \pmod k$, then $\sfA$ must vanish and hence $\calM_{j,q}(F) = 0$. On the other hand, if $q \equiv 0 \pmod k$, then the image of $v_\ell$ must lie in $F_{q/2 + \ell}$ but may otherwise be chosen arbitrarily. This applies to each $\ell \in \{-j,-j+1,\ldots,j\}$, so evaluation on the weight basis $\{v_\ell\}_{\ell=-j}^j$ defines a linear isomorphism $\calM_{j,q}(F) \to \bigoplus_{\ell=-j}^j F_{q/2+\ell}$. The lemma follows.
\end{proof}

Next, we explain the isotropy representation when $\scrE$ is the cotangent bundle of $M^\circ$. We then apply the above result to symmetric $2$-tensors and anti-self-dual $2$-forms. Finally, we discuss the case of functions.

\subsubsection{The cotangent isotropy representation} \label{subsubsec:cotangent-isotropy-rep}

Consider the complexified cotangent bundle $\scrE = T^*_{\C} M^\circ$, so that $F = T^*_{\C} M^\circ|_{(s,o)}$. We compute its isotropy representation $\rho: H \to U(F)$ using the coframe $\epsilon^0,\ldots,\epsilon^3$ from \eqref{eq:appleton-epsilon-frame}. Recall that $\rho$ is given by inverse pullback \eqref{eq:inverse-pullback-action}.

Since $G$ fixes $s$ and preserves $\sigma$, the isotropy representation fixes $\epsilon^0 = ds$ and $\epsilon^1 = a\sigma$. To determine its action on the other covectors, define
\begin{align}
    \epsilon^\pm = \frac{\epsilon^2 \pm i\epsilon^3}{\sqrt{2}}.
\end{align}
Consider the coordinate $w([z_1:z_2])=z_1/z_2$ on $\CP^1$ near $\pi(o)=[0:1]$. Recall the matrices $T_\alpha \in \mathfrak{su}(2)$ from \eqref{eq:T-matrices}, the left-invariant vector fields $X_\alpha$, and the Hopf fibration $\pi: \Sph^3 \subset \C^2 \to \CP^1$ which descends to $\pi: L_{2,k} \to \CP^1$. The curves in $\CP^1$ defined by
\begin{align}
    \gamma_1(t) &:= \pi(\exp(t T_1) (0,1)^T) = [0:e^{it}], \\
    \gamma_2(t) &:= \pi(\exp(t T_2) (0,1)^T) = [\sin t: \cos t], \\
    \gamma_3(t) &:= \pi(\exp(t T_3) (0,1)^T) = [-i\sin t: \cos t]
\end{align}
satisfy $\gamma_\alpha(0) = \pi(o)$ and $\gamma_\alpha'(0) = \pi_*(X_\alpha)$. Using these curves, we calculate at $o$,
\begin{align}
    \pi^* dw(X_1) = 0, \qquad \pi^*dw(X_2) = 1, \qquad \pi^*dw(X_3) = -i.
\end{align}
It is also clear that $\pi^*dw(\del_s) = 0$, so at the point $(s,o)$,
\begin{align} \label{eq:pullback-dw}
    \pi^*dw = \sigma_2 - i\sigma_3 = \frac{\sqrt{2}}{b} \epsilon^-.
\end{align}
Similarly,
\begin{align} \label{eq:pullback-dwbar}
    \pi^*d\bar{w} = \sigma_2 + i\sigma_3 = \frac{\sqrt{2}}{b} \epsilon^+.
\end{align}
Meanwhile, the inverse pullback of $h = \diag(u,\zeta) \in H$ on the function $\pi^*w$ is
\begin{align}
    (\phi_{h^{-1}}^* \pi^*w)(z_1,z_2) = (w \circ \pi)(u^{-1} z_1, \zeta^{-1} z_2) = u^{-1} \zeta \frac{z_1}{z_2} = \chi(h)^{-1} (\pi^*w)(z_1,z_2).
\end{align}
Thus, taking the exterior derivative gives $\phi_{h^{-1}}^* d\pi^*w = \chi(h)^{-1} d\pi^*w$. Likewise, $\phi_{h^{-1}}^* d\pi^*\bar{w} = \chi(h) d\pi^*\bar{w}$. Combining with \eqref{eq:pullback-dw} and \eqref{eq:pullback-dwbar} gives
\begin{align}
     \rho(h)\epsilon^\pm = \phi_{h^{-1}}^*\epsilon^\pm = \chi(h)^{\pm 1} \epsilon^\pm.
\end{align}
Together with $\rho(h)\epsilon^0 = \epsilon^0$ and $\rho(h)\epsilon^1 = \epsilon^1$, this fully specifies the isotropy representation. In particular, $F$ decomposes as \eqref{eq:fiber-weight-decomp} with nonzero weight spaces
\begin{align}
    F_0 = \Span_{\C}\{\epsilon^0, \epsilon^1\}, \qquad F_{\pm 1} = \Span_{\C}\{\epsilon^\pm\}.
\end{align}

\subsubsection{Symmetric 2-tensors} \label{subsubsec:symm-2-tensors}

Now let $\scrE = S^2 T^*_{\C} M^\circ$ be the bundle of symmetric $2$-tensors, so that $F = S^2 T^*_{\C} M^\circ|_{(s,o)}$. The isotropy representation of $H$ on $F$ is obtained by tensoring and symmetrizing the cotangent isotropy representation determined above. Thus, in the notation of \eqref{eq:fiber-weight-decomp}, $F$ decomposes into weight spaces
\begin{align}
    F = F_{-2} \oplus F_{-1} \oplus F_0 \oplus F_1 \oplus F_2
\end{align}
where, using the symmetric product notation \eqref{eq:symm-prod},
\begin{align}
    F_0 &= \Span_{\C}\{(\epsilon^0)^2, \, (\epsilon^1)^2, \, \epsilon^0 \odot \epsilon^1, \, \epsilon^+ \odot \epsilon^- \}, \\
    F_{\pm 1} &= \Span_{\C} \{ \epsilon^0 \odot \epsilon^\pm, \, \epsilon^1 \odot \epsilon^\pm \}, \\
    F_{\pm 2} &= \Span_{\C} \{ (\epsilon^\pm)^2 \}.
\end{align}
Thus
\begin{align} \label{eq:symmetric-tensor-dims}
    \dim_{\C} F_0 = 4, \qquad \dim_{\C} F_{\pm 1} = 2, \qquad \dim_{\C} F_{\pm 2} = 1.
\end{align}
By Lemma \ref{lem:weight-matching}, $\calM_{j,q}(F) = 0$ unless  $q \equiv 0 \pmod k$, in which case
\begin{align} \label{eq:symmetric-tensor-multiplicity}
    \calM_{j,q}(F) \cong \bigoplus_{\substack{r \in \{-2,-1,0,1,2\} \\ |r-q/2| \leq j}} F_r.
\end{align}
In particular, among pairs $(j,q)$ satisfying \eqref{eq:jq-representation-labels}, the $(j,q)$-sector of $\scrE_s$ is nonzero if and only if
\begin{align} \label{eq:jq-admissibility-conditions}
    q \equiv 0 \pmod k, \qquad |q| \leq 2j+4.
\end{align}

\subsubsection{Anti-self-dual 2-forms} \label{subsubsec:asd-2-forms}

Equip $M^\circ$ with the orientation for which $-\epsilon^0 \wedge \epsilon^1 \wedge \epsilon^2 \wedge \epsilon^3$ is positive. Let $\scrE = \Lambda^-_{\C} \to M^\circ$ denote the bundle of anti-self-dual (ASD) $2$-forms. Then $F = \Lambda^-_{\C}|_{(s,o)}$. Since $G$ acts by orientation-preserving isometries, its action on $2$-forms commutes with the Hodge star and hence preserves $\scrE$. Therefore, the results of this section apply to $\scrE$.

Recall the notation $\epsilon^{ij} = \epsilon^i \wedge \epsilon^j$, where the wedge product is unnormalized \eqref{eq:wedge-prod}. With respect to the form inner product \eqref{eq:form-inner-product}, an orthonormal basis of $F$ is given by
\begin{align} \label{eq:asd-fiber-weight-basis}
    \sigma_0 = -\frac{\epsilon^{01} + \epsilon^{23}}{\sqrt{2}}, \qquad \sigma_{\pm} = \frac{1}{2}(\epsilon^{02} - \epsilon^{13} \pm i(\epsilon^{03} + \epsilon^{12})) = \frac{1}{\sqrt{2}} (\epsilon^0 \pm i\epsilon^1) \wedge \epsilon^{\pm}.
\end{align}
The isotropy representation of $H$ on $F$ is the exterior product of the cotangent isotropy representation above. Thus, in the notation of \eqref{eq:fiber-weight-decomp}, $F$ decomposes into weight spaces
\begin{align}
    F = F_{-1} \oplus F_0 \oplus F_1,
\end{align}
where $F_0 = \Span_{\C}\{\sigma_0\}$ and $F_{\pm 1} = \Span_{\C}\{\sigma_\pm\}$. By Lemma \ref{lem:weight-matching}, $\calM_{j,q}(F) = 0$ unless $q \equiv 0 \pmod k$, in which case
\begin{align}
    \calM_{j,q}(F) \cong \bigoplus_{\substack{r \in \{-1,0,1\} \\ |r-q/2| \leq j}} F_r. 
\end{align}

\subsubsection{Complex-valued functions} \label{subsubsec:complex-valued-functions}

Finally, consider the trivial complex line bundle $\scrE \to M^\circ$, whose sections are complex-valued functions. The isotropy representation of $H$ on $F$ is trivial, so $F$ decomposes as \eqref{eq:fiber-weight-decomp} with $F_0 = \C$ and $F_r = 0$ for $r \neq 0$. Lemma \ref{lem:weight-matching} therefore gives, for every $(j,q)$ satisfying \eqref{eq:jq-representation-labels},
\begin{align} \label{eq:function-multiplicity}
    \calM_{j,q}(\C) \cong \begin{cases}
        \C &\text{if } q \equiv 0 \pmod k \text{ and } |q| \leq 2j, \\
        0 &\text{otherwise}.
    \end{cases}
\end{align}
Thus, Proposition \ref{prop:sector-peter-weyl} gives
\begin{align} \label{eq:L2-functions-identification}
    L^2(L_{2,k}, d\omega) = \widehat{\bigoplus}_{(j,q) \text{ satisfying } \eqref{eq:jq-representation-labels}, \, q \equiv 0 \text{ mod } k, \, |q| \leq 2j} \calH_{j,q}, \qquad \calH_{j,q} \cong_G V_{j,q}.
\end{align}
We call the $(j,q)$-sector of $\scrE_s$, i.e. $\calH_{j,q}$, the \emph{scalar $(j,q)$-sector}.

For each $s>0$, the scalar Laplacian $\Delta_{\Sigma_s}$ with respect to the induced metric is $G$-equivariant, so \eqref{eq:L2-functions-identification} and Schur's lemma imply that $\Delta_{\Sigma_s}$ acts by scalar multiplication on each scalar $(j,q)$-sector. The corresponding eigenvalues are shown in Corollary \ref{cor:scalar-eigenvalues}.

For $j \in \Z_{\geq 0}$, Lemma \ref{lem:hopf-invariance} shows that functions in $\calH_{j,0}$ are Hopf-invariant. Thus, if $u \in \calH_{j,0}$, then $u$ descends to a function $\bar{u}$ on $\CP^1$. The map $u \mapsto \bar{u}$ is then injective on $\calH_{j,0}$. Since $u = \pi^*\bar{u}$ and each round Hopf fiber has length $2\pi/k$,
\begin{align}
    \norm{u}_{L^2(L_{2,k}, d\omega)}^2 = \frac{2\pi}{k} \lVert \bar{u} \rVert_{L^2(\CP^1,g_{\mathrm{FS}})}^2.
\end{align}
We now relate $\bar{u}$ to spherical harmonics on $\CP^1$. Consider the natural action of $G$ on $\CP^1$, given by $\theta \cdot [z] = [\theta z]$. This action preserves $g_{\mathrm{FS}}$ and its orientation. Moreover, the Hopf projection $\pi: L_{2,k} \to \CP^1$ is $G$-equivariant. The isotropy subgroup at $[0:1] \in \CP^1$ is $K = U(1) \times U(1)$, giving a $G$-equivariant identification $\CP^1 \cong G/K$. Under this identification and $L_{2,k} \cong G/H$, the Hopf projection is the canonical map $G/H \to G/K$.

We can apply the Peter--Weyl decomposition (Proposition \ref{prop:sector-peter-weyl}) to the trivial line bundle over $G/K$, and redo the analysis at the start of \S\ref{subsec:isotropy-weights-examples} with $K$ in place of $H$. Combining the results as we have done to obtain \eqref{eq:L2-functions-identification}, we get
\begin{align} \label{eq:spherical-harmonics-decomp}
    L^2(\CP^1, g_{\mathrm{FS}}) = \widehat{\bigoplus}_{j \in \Z_{\geq 0}} \, \bar{\calH}_j, \qquad \bar{\calH}_j \cong_G V_{j,0}.
\end{align}
This is exactly the standard decomposition into spherical harmonics. It is well-known that
\begin{align}
    \bar{\calH}_j = \{ Y \in C^\infty(\CP^1;\C) : -\Delta_{g_{\mathrm{FS}}} Y = 4j(j+1) Y \},
\end{align}
as $g_{\mathrm{FS}}$ has constant sectional curvature $4$.

\begin{lemma} \label{lem:hopf-projection-harmonics}
    For each $j \in \Z_{\geq 0}$, the map $\calH_{j,0} \to L^2(\CP^1, g_{\mathrm{FS}})$ defined by $u \mapsto \sqrt{2\pi/k} \cdot \bar{u}$ has image $\bar{\calH}_j$. It defines a unitary $G$-equivariant isomorphism from $\calH_{j,0}$ to $\bar{\calH}_j$, with inverse $u \mapsto \sqrt{k/2\pi} \, \pi^*u$.
\end{lemma}
\begin{proof}
    As discussed above, the Hopf projection is the canonical $G$-equivariant map $G/H \to G/K$. This implies that the map $\calH_{j,0} \to L^2(\CP^1, g_{\mathrm{FS}})$, $u \mapsto \sqrt{2\pi/k} \, \bar{u}$ is isometric, $G$-equivariant, and injective. The claim now follows from \eqref{eq:spherical-harmonics-decomp} and Schur's lemma.
\end{proof}

\subsection{Differential identities}

As in \S\ref{subsec:group-action-peter-weyl} and \S\ref{subsec:isotropy-weights-examples}, fix a bundle $\scrE \to M^\circ$ and $s>0$, and abbreviate $F=F_s$. We also assume that the isotropy representation splits $F$ as in \eqref{eq:fiber-weight-decomp}.

For each weight space $F_r$, choose an orthonormal basis $\{\tau_{r,n}\}_{n=1}^{\dim F_r}$. The identification \eqref{eq:E-homog-vb} allows us to uniquely expand any smooth section $u \in C^\infty(\scrE_s)$ as
\begin{align} \label{eq:u-identification}
    u(\theta H) = \sum_{r,n} u_{r,n}(\theta) [\theta, \tau_{r,n}],
\end{align}
where $u_{r,n}: G \to \C$ are functions satisfying
\begin{align} \label{eq:coeff-equivariance}
    u_{r,n}(\theta h) = \chi(h)^{-r} u_{r,n}(\theta) \quad \text{for all } h \in H.
\end{align}
The next lemma records differential identities for the coefficients $u_{r,n}(\theta)$ when $u$ belongs to a fixed sector. Below, we regard the generators $T_\alpha \in \mathfrak{su}(2)$ from \S\ref{subsec:local-frame} as elements of $\mathfrak{u}(2)$, and use $X_\alpha$ to denote the corresponding left-invariant vector fields on $G=U(2)$.

\begin{lemma} \label{lem:differential-identities}
    In the setting above, if $u$ belongs to the $(j,q)$-sector, then 
    \begin{align}
        X_1 u_{r,n} &= -i(q-2r) u_{r,n}, \label{eq:X1-differential-identity} \\
        (X_1^2 + X_2^2 + X_3^2) u_{r,n} &= -4j(j+1) u_{r,n}. \label{eq:X2-differential-identity}
    \end{align}
    Moreover, for each $r,n$, the functions $|u_{r,n}|^2$, $|X_1 u_{r,n}|^2$, and $|X_2 u_{r,n}|^2 + |X_3 u_{r,n}|^2$ descend globally to $L_{2,k}$, and
    \begin{align} \label{eq:X-integral-identity}
        \int_{L_{2,k}} \left( |X_2 u_{r,n}|^2 + |X_3 u_{r,n}|^2 \right) \, d\omega = (4j(j+1) - (q-2r)^2) \int_{L_{2,k}} |u_{r,n}|^2 \, d\omega.
    \end{align}
\end{lemma}
\begin{proof}
    To prove \eqref{eq:X1-differential-identity} and \eqref{eq:X2-differential-identity}, we may assume by linearity that $u = \scrT_{j,q,s}(v \otimes \sfA)$, where $\scrT_{j,q,s}$ is defined in Proposition \ref{prop:sector-peter-weyl}, and $v \in V_{j,q}$, $\sfA \in \calM_{j,q}(F)$. Then using \eqref{eq:u-identification} and \eqref{eq:U-identification}, the function $U: G \to F$ representing $u$ is given by the matrix coefficient
    \begin{align}
        U(\theta) = \sum_{r,n} u_{r,n}(\theta) \tau_{r,n} = \sfA(\pi(\theta^{-1})v),
    \end{align}
    where $\pi := \pi_{j,q}$. Standard formulas for differentiating matrix coefficients give, for $\alpha=1,2,3$,
    \begin{align}
        X_\alpha U(\theta) &= -\sfA(d\pi(T_\alpha) \pi(\theta^{-1}) v), \label{eq:diff1-identity} \\
        X_\alpha^2 U(\theta) &= \sfA(d\pi(T_\alpha)^2 \pi(\theta^{-1}) v). \label{eq:diff2-identity}
    \end{align}
    Let $\Pi_r: F \to F_r$ be the orthogonal projection. By Lemma \ref{lem:weight-matching}, $\Pi_r\sfA$ vanishes on every weight vector $v_{\ell} \in V_{j,q}$ except possibly $v_{r-q/2}$. By \eqref{eq:representation-torus-weights}, we have $d\pi(T_1) v_\ell = -2i\ell v_\ell$ and so
    \begin{align}
        \Pi_r \sfA \circ d\pi(T_1) = -2i(r-q/2) \Pi_r \sfA.
    \end{align}
    Applying $\Pi_r$ to \eqref{eq:diff1-identity} therefore yields
    \begin{align}
        \sum_n (X_1 u_{r,n})(\theta) \tau_{r,n} &= \Pi_r(X_1 U(\theta)) \\
        &= -i(q-2r) \Pi_r(U(\theta)) = -i(q-2r) \sum_n u_{r,n}(\theta) \tau_{r,n}.
    \end{align}
    Comparing coefficients gives \eqref{eq:X1-differential-identity}. To prove \eqref{eq:X2-differential-identity}, simply substitute into \eqref{eq:diff2-identity} the well-known Casimir identity (which may also be deduced using \eqref{eq:representation-torus-weights})
    \begin{align}
        d\pi(T_1)^2 + d\pi(T_2)^2 + d\pi(T_3)^2 = -4j(j+1) I_{V_{j,q}}.
    \end{align}

    Next, fix a single coefficient $u_{r,n}$. By \eqref{eq:coeff-equivariance}, the function $|u_{r,n}|^2$ descends to $\Sigma_s$. By \eqref{eq:X1-differential-identity}, so does $|X_1 u_{r,n}|^2$. For $\mu=2,3$, differentiating along the curve $t \mapsto \theta h \exp(t T_\mu) = \theta \exp(t \operatorname{Ad}_h(T_\mu)) h$ gives
    \begin{align} \label{eq:Xmu-on-coeff}
        (X_\mu u_{r,n})(\theta h) = \chi(h)^{-r} (X_{\operatorname{Ad}_h(T_\mu)} u_{r,n})(\theta),
    \end{align}
    where $X_T$ denotes the left-invariant vector field generated by $T \in \mathfrak{u}(2)$. It is routine to check that the adjoint action of $H$ rotates the plane spanned by $T_2,T_3$. Thus, by \eqref{eq:Xmu-on-coeff}, the function $|X_2 u_{r,n}|^2 + |X_3 u_{r,n}|^2$ descends to $\Sigma_s$.

    It remains to prove \eqref{eq:X-integral-identity}. Let $d\theta$ be normalized Haar measure on $G$. Integration by parts and \eqref{eq:X1-differential-identity}, \eqref{eq:X2-differential-identity} give
    \begin{align} \label{eq:X-integral-identity-G}
        \int_G \left( |X_2 u_{r,n}|^2 + |X_3 u_{r,n}|^2 \right) \, d\theta = (4j(j+1) - (q-2r)^2) \int_G |u_{r,n}|^2 \, d\theta.
    \end{align}
    The pushforward of $d\theta$ under $G \to G/H \cong L_{2,k}$ is the $G$-invariant probability measure, which is $d\omega/\vol_\omega(L_{2,k})$. The integrands descend to $L_{2,k}$, as already shown. Thus, \eqref{eq:X-integral-identity-G} holds with $G$ replaced by $L_{2,k}$, proving the asserted identity.
\end{proof}

\begin{corollary} \label{cor:scalar-eigenvalues}
    Suppose $u$ belongs to the scalar $(j,q)$-sector, and let $\Delta_{\Sigma_s}$ be the scalar Laplacian on $\Sigma_s \cong L_{2,k}$ with respect to the induced metric. Then
    \begin{align}
        -\Delta_{\Sigma_s} u = \left( \frac{q^2}{a^2} + \frac{4j(j+1)-q^2}{b^2} \right) u.
    \end{align}
\end{corollary}
\begin{proof}
    By the discussion in \S\ref{subsubsec:complex-valued-functions}, if the scalar $(j,q)$-sector is nontrivial, then $F$ is one-dimensional with $F = F_0$. Thus, the only coefficient function $u_{r,n}: G \to \C$ descends to $G/H \cong L_{2,k}$ and is given by $u$ itself. The connection table \eqref{eq:frame-connection} shows that $\nabla_{E_\alpha}^{\Sigma_s} E_\alpha = 0$ for $\alpha=1,2,3$. Thus,
    \begin{align}
        -\Delta_{\Sigma_s} u &= -\sum_{\alpha=1}^3 E_\alpha^2 u = -\frac{1}{b^2}(X_1^2 + X_2^2 + X_3^2) u + \left( \frac{1}{b^2} - \frac{1}{a^2} \right) X_1^2 u.
    \end{align}
    Substituting \eqref{eq:X1-differential-identity} and \eqref{eq:X2-differential-identity} then establishes the asserted formula.
\end{proof}

\subsection{The separated stability problem}

Call a pair $(j,q)$ \emph{admissible} if it satisfies \eqref{eq:jq-representation-labels} and \eqref{eq:jq-admissibility-conditions}. We have seen in \S\ref{subsubsec:symm-2-tensors} that $(j,q)$ is admissible if and only if its corresponding sector in $S^2 T^*_{\C}M^\circ$ is nonzero. We call the corresponding sector an \emph{admissible sector}.

We now apply the preceding results to reduce the strict linear stability problem to one that can be studied separately in each admissible sector.

\begin{proposition} \label{prop:sector-stability-reduction}
    Suppose there exists $\mu > 0$ such that for every smooth symmetric $2$-tensor compactly supported in $M^\circ$ and belonging to a single admissible sector,
    \begin{align} \label{eq:sector-coercivity-hypothesis}
        \calE[h] \geq \mu \int_M \frac{|h|^2}{b^2} e^f.
    \end{align}
    Then $(M,g,f)$ is strictly linearly stable.
\end{proposition}
\begin{proof}
    Let $h$ be a smooth symmetric $2$-tensor compactly supported in $M^\circ$. For each admissible pair $(j,q)$, let $\Pi_{j,q}$ denote the orthogonal projection onto the $(j,q)$-sector of $\scrE_s$ for each $s>0$. Set $h_{j,q} = \Pi_{j,q} h$. By the formula \eqref{eq:jq-orthog-projection}, and smoothness of the $G$-action, $h_{j,q}$ is smooth and compactly supported in $M^\circ$.
    
    Since $G$ acts by isometries, the operators
    \begin{align}
        \nabla: C^\infty(S^2 T^*_{\C}M^\circ) &\to C^\infty(T^*_{\C}M^\circ \otimes S^2 T^*_{\C}M^\circ), \\
        \Rm: C^\infty(S^2 T^*_{\C}M^\circ) &\to C^\infty(S^2 T^*_{\C}M^\circ)
    \end{align}
    are linear and $G$-equivariant. By Corollary \ref{cor:equivariant-operator-splitting}, $\nabla h_{j,q}$ and $\Rm(h_{j,q})$ also belong to the $(j,q)$-sector. The Peter--Weyl decomposition (Proposition \ref{prop:sector-peter-weyl}) also gives pairwise orthogonality between sectors along each level set of $s$, so
    \begin{align}
        \calE[h] &= \sum_{j,q,j',q'} \int_M \left( \inner{\nabla h_{j,q}}{\nabla h_{j',q'}} - 2\inner{\Rm(h_{j,q})}{h_{j',q'}} \right) e^f \\
        &= \sum_{j,q,j',q'} \int_0^\infty w \int_{L_{2,k}} \left( \inner{\nabla h_{j,q}}{\nabla h_{j',q'}} - 2\inner{\Rm(h_{j,q})}{h_{j',q'}} \right) \, d\omega \, ds \\
        &= \sum_{j,q} \int_0^\infty w \int_{L_{2,k}} \left( \inner{\nabla h_{j,q}}{\nabla h_{j,q}} - 2\inner{\Rm(h_{j,q})}{h_{j,q}} \right) \, d\omega \, ds = \sum_{j,q} \calE[h_{j,q}].
    \end{align}
    Here the sum runs over all admissible pairs. Now using \eqref{eq:sector-coercivity-hypothesis}, and the fact that $b$ and $f$ are functions of $s$, we get
    \begin{align}
        \calE[h] &\geq \mu \sum_{j,q} \int_M \frac{|h_{j,q}|^2}{b^2} e^f \\
        &= \mu \sum_{j,q} \int_0^\infty \frac{w}{b^2} \int_{L_{2,k}}  |h_{j,q}|^2 \, d\omega \, ds \\
        &= \mu \int_0^\infty \frac{w}{b^2} \int_{L_{2,k}} |h|^2 \, d\omega \, ds = \mu \int_M \frac{|h|^2}{b^2} e^f.
    \end{align}
    Strict linear stability now follows from Proposition \ref{prop:stability-reduction}.
\end{proof}

\section{The stability inequality in $q \neq 0$ sectors} \label{sec:charged-stability}

In this section, we establish the stability estimate \eqref{eq:sector-coercivity-hypothesis} in admissible sectors with $q \neq 0$:
\begin{theorem} \label{thm:charged-coercivity}
    There exists $\mu > 0$ such that every smooth symmetric $2$-tensor $h$ compactly supported in $M^\circ$ and belonging to an admissible $(j,q)$-sector with $q \neq 0$ satisfies
    \begin{align}
        \calE[h] \geq \mu \int_M \frac{|h|^2}{b^2} e^f.
    \end{align}
\end{theorem}

\subsection{Conformally K\"ahler structure and Hermitian decomposition} \label{subsec:conformally-kahler}

Recall that $E_0$ and $E_1$ from \S\ref{subsec:local-frame} are globally defined on $M^\circ$, while $E_2$ and $E_3$ are locally defined using lifts to $(0,\infty) \times SU(2)$. Also recall the orientation fixed in \S\ref{subsubsec:asd-2-forms}, which makes $-\epsilon^0 \wedge \epsilon^1 \wedge \epsilon^2 \wedge \epsilon^3$ positive. There is a well-defined orthogonal almost complex structure $J$ on $M^\circ$ given by
\begin{align}
    JE_0 = -E_1, \qquad JE_1 = E_0, \qquad JE_2 = E_3, \qquad JE_3 = -E_2.
\end{align}
Indeed, we can alternatively define $J$ by setting $JE_0 = -E_1$ and $JE_1 = E_0$, then defining $J$ to be the standard complex structure on $\Span_{\R}\{E_0,E_1\}^\perp \cong \C$ where an orientation-preserving isometry is chosen. This description is global. Notice that $J$ is $G$-invariant, since the action of $G$ fixes $E_0,E_1$ and rotates the $(E_2,E_3)$-plane in an orientation-preserving manner.

The fundamental form associated to $g$ and $J$ is the self-dual $2$-form
\begin{align} \label{eq:fund-2form}
    \Omega = g(J\cdot,\cdot) = -\epsilon^{01} + \epsilon^{23}.
\end{align}

\begin{lemma} \label{lem:conformally-kahler}
    The almost complex structure $J$ is integrable and satisfies $\nabla_{E_0} J = 0$. If $u=u(s)$ satisfies $u' = D-B$, then $\check{g} = e^{2u}g$ is K\"ahler with respect to $J$ on $M^\circ$.
\end{lemma}
\begin{proof}
    The complex covectors
    \begin{align}
        \zeta^0 = \frac{\epsilon^0 - i\epsilon^1}{\sqrt{2}}, \qquad \zeta^1 = \frac{\epsilon^2 + i\epsilon^3}{\sqrt{2}}
    \end{align}
    satisfy $\zeta^0 \circ J = i\zeta^0$ and $\zeta^1 \circ J = i\zeta^1$, so $(\zeta^0,\zeta^1)$ is a unitary coframe of $(1,0)$-forms for $J$. The structure equations obtained from \eqref{eq:frame-connection} give
    \begin{align}
        d\zeta^0 &= -\frac{1}{\sqrt{2}}(A \zeta^0 \wedge \bar{\zeta}^0 + 2D\zeta^1 \wedge \bar{\zeta}^1), \qquad d\zeta^1 = \left( B\epsilon^0 - \frac{2i}{a}\epsilon^1 \right) \wedge \zeta^1.
    \end{align}
    Since $\epsilon^0 = \frac{1}{\sqrt{2}} (\zeta^0 + \bar{\zeta}^0)$ and $\epsilon^1 = \frac{i}{\sqrt{2}} (\zeta^0 - \bar{\zeta}^0)$, neither $d\zeta^0$ nor $d\zeta^1$ has a $(0,2)$-component. This proves integrability. Also, the frame $\{E_i\}_{i=0}^3$ is parallel along $E_0$, and $J$ has constant coefficients in this frame, so $\nabla_{E_0}J = 0$.

    Finally, using the structure equations again, we compute $d\Omega = 2(B-D) \, ds \wedge \Omega$. Hence if $u' = D-B$, then $d(e^{2u} \Omega) = 2e^{2u} (u' + B-D) \, ds \wedge \Omega = 0$. Thus, $\check{g} = e^{2u}g$ is K\"ahler.
\end{proof}

For a symmetric $2$-tensor $h$, define $(J^*h)(X,Y) = h(JX,JY)$. Since $J$ is orthogonal, $J^*$ is a self-adjoint involution on $S^2 T^*_{\C}M^\circ$. Its $+1$ and $-1$ eigenspaces are the Hermitian and skew-Hermitian subbundles, respectively. The corresponding orthogonal projections are
\begin{align}
    h_H = \frac{1}{2}(h + J^*h), \qquad h_A = \frac{1}{2}(h - J^*h).
\end{align}
Since $J^*g = g$, every skew-Hermitian tensor is orthogonal to $g$, thus is trace-free. We therefore have the pointwise orthogonal decomposition
\begin{align} \label{eq:charged-tensor-decomposition}
    h = \phi\sfE + h_{0,H} + h_A, \qquad |h|^2 = |\phi|^2 + |h_{0,H}|^2 + |h_A|^2,
\end{align}
where $\sfE = g/2$ is a unit tensor, $\phi = \inner{h}{\sfE} = \frac{1}{2} \tr_g h$, and $h_{0,H} = h_H - \phi\sfE$. Write $S_{0,H}^2$ and $S_A^2$ for the bundles of trace-free Hermitian and skew-Hermitian tensors. Then the bundle of trace-free symmetric 2-tensors $S_0^2 T^*_{\C}M^\circ$ splits as
\begin{align} \label{eq:tracefree-H-A-splitting}
    S_0^2 T^*_{\C}M^\circ = S_{0,H}^2 \oplus S_A^2.
\end{align}
Because $g$ and $J$ are $G$-invariant, projections onto $S_{0,H}^2$ and $S_A^2$ are $G$-equivariant, and so is taking the trace. Thus, by Corollary \ref{cor:equivariant-operator-splitting}, if $h$ belongs to a $(j,q)$-sector, then $h_{0,H}$ and $h_A$ also belong to the same $(j,q)$-sector, while $\phi$ belongs to the scalar $(j,q)$-sector.

\subsection{The Naff--Ozuch Weitzenb\"ock formula}

We recall the Naff--Ozuch Weitzenb\"ock formula on a four-dimensional gradient Ricci soliton. First, we introduce the relevant objects.

Let $\Lambda^2_{\C}$ denote the bundle of $2$-forms, and denote by $\Lambda_{\C}^\pm$ the self-dual and anti-self-dual (ASD) subbundles. We equip $\Lambda^2_{\C} \otimes \Lambda^2_{\C}$ and subbundles thereof with the inner product induced by the form inner products \eqref{eq:form-inner-product} on each factor. Define the trace map
\begin{align}
    \tr: \Lambda^2_{\C} \otimes \Lambda^2_{\C} \to S^2 T^*_{\C} M^\circ, \qquad \tr(\scrS)_{ij} = \frac{1}{2} \sum_k (\scrS_{ikjk} + \scrS_{jkik}),
\end{align}
where the sum is taken over an orthonormal frame. (This operation is not to be confused with the trace of a symmetric $2$-tensor.)

The trace map is $G$-equivariant. Its restriction to $\Lambda^-_{\C} \otimes \Lambda^+_{\C}$ defines an isometric bundle isomorphism
\begin{align} \label{eq:trace-isomorphism}
    \Lambda^-_{\C} \otimes \Lambda^+_{\C} \cong S_0^2 T^*_{\C}M^\circ.
\end{align}
Meanwhile, its restriction to $\Lambda^+_{\C} \otimes \Lambda^+_{\C}$ takes values in $\C g$, which is pointwise orthogonal to trace-free symmetric tensors. Consequently, for $\scrT \in \Lambda^2_{\C} \otimes \Lambda^+_{\C}$ and $\scrS \in \Lambda^-_{\C} \otimes \Lambda^+_{\C}$,
\begin{align} \label{eq:trace-isometry-2}
    \inner{\tr(\scrT)}{\tr(\scrS)} = \inner{\scrT}{\scrS}.
\end{align}
These facts can be verified using the orthonormal bases $(\sigma_0,\sigma_+,\sigma_-)$ for $\Lambda^-_{\C}$ from \eqref{eq:asd-fiber-weight-basis}, and
\begin{align} \label{eq:lambda-plus-onb}
    \left\{ \hat{\Omega}, \, \frac{\epsilon^{02} + \epsilon^{13}}{\sqrt{2}}, \, \frac{\epsilon^{03} - \epsilon^{12}}{\sqrt{2}} \right\}
\end{align}
for $\Lambda^+_{\C}$. Here $\hat{\Omega} = \Omega/\sqrt{2}$ is the normalized fundamental form, cf. \eqref{eq:fund-2form}.

We now describe the Hermitian splitting \eqref{eq:tracefree-H-A-splitting} under the isomorphism \eqref{eq:trace-isomorphism}. One checks that $J^*$ acts as the identity on $\Lambda_{\C}^-$ and on $\C \hat{\Omega}$, while it acts as minus the identity on $(\C\hat{\Omega})^\perp \subset \Lambda^+_{\C}$. Since $J$ is orthogonal, we have $J^*\tr(\alpha \otimes \beta) = \tr(J^*\alpha \otimes J^*\beta)$ for any $2$-forms $\alpha,\beta$. Thus, the map \eqref{eq:trace-isomorphism} restricts to isomorphisms
\begin{align}
    \Lambda^-_{\C} \otimes \C\hat{\Omega} \cong S_{0,H}^2, \qquad \Lambda^-_{\C} \otimes (\C\hat{\Omega})^\perp \cong S_A^2.
\end{align}
In particular, writing a trace-free symmetric $2$-tensor as $h_0 = h_{0,H} + h_A = \tr(\scrS)$, we have a corresponding decomposition
\begin{align} \label{eq:scrS-decomp}
    \scrS = \scrS_H + \scrS_A, \qquad \scrS_H = \alpha \otimes \hat{\Omega}
\end{align}
for a unique complex ASD $2$-form $\alpha$ satisfying $h_{0,H} = \tr(\alpha \otimes \hat{\Omega})$ and $|h_{0,H}| = |\alpha|$. Because the trace map is $G$-equivariant and $\hat{\Omega}$ is $G$-invariant, the identification $h_{0,H} \leftrightarrow \alpha$ is $G$-equivariant. Hence, if $h_{0,H}$ belongs to a $(j,q)$-sector, then so does $\alpha$.

Regard sections of $\Lambda^-_{\C} \otimes \Lambda^+_{\C}$ as $\Lambda^+_{\C}$-valued $2$-forms. Let $d^L$ be the exterior covariant derivative on $\Lambda^+_{\C}$-valued forms, and let $\delta^L$ be its (unweighted) formal adjoint. Explicitly, for a $p$-form $\alpha$ and a section $\beta$ of $\Lambda^+_{\C}$,
\begin{alignat}{2}
    d^L(\alpha \otimes \beta)
    &= (d\alpha) \otimes \beta
       + (-1)^p \sum_{i=0}^3
       (\alpha \wedge \epsilon^i) \otimes \nabla_{E_i}\beta,
    &\qquad d\alpha
    &= \epsilon^i \wedge \nabla_{E_i}\alpha, \\
    \delta^L(\alpha \otimes \beta)
    &= (\delta\alpha) \otimes \beta
       - \sum_{i=0}^3
       \iota_{E_i}\alpha \otimes \nabla_{E_i}\beta,
    &\qquad \delta\alpha
    &= -\iota_{E_i}(\nabla_{E_i}\alpha).
\end{alignat}
Moreover, let
\begin{align}
    \delta^L_f = e^{-f} \delta^L e^f = \delta^L - \iota_{\nabla f}^L
\end{align}
be the formal adjoint of $d^L$ with respect to the weighted measure $e^f\, \dvol_g$. Here $\iota_{\nabla f}^L(\alpha \otimes \beta) = (\iota_{\nabla f}\alpha) \otimes \beta$. Using these operators, we define the $L^2(e^f)$-self-adjoint weighted Hodge Laplacian
\begin{align}
    \Delta^L_{H,f} = -(d^L \delta^L_f + \delta^L_f d^L).
\end{align}
Our convention for $d^L$ differs from that in \cite{naff-ozuch-kahler} by a constant factor depending on the form degree, while our $\delta^L$ and $\iota$ differ by the reciprocal factor, so our $\Delta^L_{H,f}$ agrees with theirs.

Let $W^+: \Lambda_{\C}^+ \to \Lambda_{\C}^+$ be the self-dual Weyl curvature operator, normalized as in \cite{naff-ozuch-kahler}.
Finally, for $\scrS \in \Lambda^2_{\C} \otimes \Lambda^2_{\C}$ and an endomorphism $U$ of $\Lambda^2_{\C}$, write $\scrS \circ U$ for its action on the right form factor, i.e. $(\alpha \otimes \beta) \circ U = \alpha \otimes U\beta$.

We state the Naff--Ozuch Weitzenb\"ock formula, specialized to steady gradient Ricci solitons. 

\begin{proposition}[\cite{naff-ozuch-kahler}*{Proposition 8}] \label{prop:naff-ozuch-weitzenbock}
    Let $(M,g,f)$ be a oriented four-dimensional steady gradient Ricci soliton. Let $h$ be a smooth symmetric $2$-tensor on $M$ and decompose it as
    \begin{align}
        h = \phi\sfE + \tr(\scrS), \qquad \sfE = \frac{g}{2}, \qquad \scrS \in C^\infty(\Lambda_{\C}^- \otimes \Lambda_{\C}^+).
    \end{align}
    Then
    \begin{align} \label{eq:naff-ozuch-weitzenbock}
        \Delta_f h + 2\Rm(h) = (\Delta_f \phi) \sfE + \phi\Ric + \tr\left( \Delta^L_{H,f} \scrS + \scrS \circ \left( W^+ - \frac{R}{3} I_{\Lambda_{\C}^+} \right) \right).
    \end{align}
\end{proposition}

Returning to our four-dimensional Appleton soliton on $\O_{\CP^1}(-k)$, the above Weitzenb\"ock formula and the conformally K\"ahler structure yield a useful decomposition of the energy of a trace-free tensor. To state this, for a smooth compactly supported section $\scrS$ of $\Lambda_{\C}^- \otimes \Lambda_{\C}^+$, set
\begin{align} \label{eq:Qf-defn}
    \calQ_f[\scrS] &:= -\int_M \inner{\tr(\Delta^L_{H,f} \scrS)}{\tr(\scrS)} e^f.
\end{align}
Using \eqref{eq:trace-isometry-2} and integrating by parts,
\begin{align}
    \calQ_f[\scrS] = -\int_M \inner{\Delta^L_{H,f} \scrS}{\scrS} e^f = \int_M \left( |d^L \scrS|^2 + |\delta^L_f \scrS|^2 \right) e^f. \label{eq:Qf-identity-1}
\end{align}

\begin{proposition} \label{prop:tracefree-energy}
    Let $h_0 = \tr(\scrS)$ be a smooth trace-free symmetric $2$-tensor compactly supported in $M^\circ$. Decompose $h_0 = h_{0,H} + h_A$, where $h_{0,H} \in C^\infty(S_{0,H}^2)$ and $h_A \in C^\infty(S_A^2)$. Then
    \begin{align} \label{eq:tracefree-energy}
        \calE[h_0] = \calQ_f[\scrS] + \int_M \left( V_H|h_{0,H}|^2 + V_A|h_A|^2 \right) e^f,
    \end{align}
    where, in terms of the curvature coefficients \eqref{eq:frame-curvature},
    \begin{align} \label{eq:V-potentials}
        V_H &= 2(K_2 + H_{12} + 2M_2), \\
        V_A &= K_1 + K_2 + H_{12} + H_{23} - 2M_2.
    \end{align}
\end{proposition}
\begin{proof}
    Let $\check{g} = e^{2u}g$ be the K\"ahler metric from Lemma \ref{lem:conformally-kahler}, let $\check{R}$ be its scalar curvature, and write $\kappa = e^{2u} \check{R}$. The standard conformal change formula for scalar curvature gives
    \begin{align} \label{eq:scalar-conformal-change}
        \kappa = R - 6\Delta u - 6|\nabla u|^2.
    \end{align}
    The self-dual and ASD bundles are invariant under conformal change, so the self-dual Weyl operator $\check{W}^+$ of $\check{g}$ also maps $\Lambda^+_{\C} \to \Lambda^+_{\C}$. With our normalization, which agrees with \cite{naff-ozuch-kahler}, it is a well-known fact of K\"ahler surfaces that $\check{W}^+$ has eigenvalue $\check{R}/3$ on the line $\C\hat{\Omega}$, and $-\check{R}/6$ on its $\check{g}$- (hence $g$-)orthogonal complement in $\Lambda^+_{\C}$. Moreover, under conformal change, $\check{W}^+ = e^{-2u} W^+$. Thus,
    \begin{align}
        W^+|_{\C\hat{\Omega}} = \frac{\kappa}{3} I, \qquad W^+|_{(\C\hat{\Omega})^\perp} = -\frac{\kappa}{6} I.
    \end{align}
    Writing $\scrS = \scrS_H + \scrS_A$ as in \eqref{eq:scrS-decomp}, we therefore obtain
    \begin{align}
        \scrS \circ \left( W^+ - \frac{R}{3} I_{\Lambda_{\C}^+} \right) = \frac{\kappa-R}{3} \scrS_H - \frac{\kappa+2R}{6} \scrS_A.
    \end{align}
    Applying Proposition \ref{prop:naff-ozuch-weitzenbock}, we therefore get
    \begin{align}
        \Delta_f h_0 + 2\Rm(h_0) = \tr(\Delta^L_{H,f}\scrS) + \frac{\kappa-R}{3} h_{0,H} - \frac{\kappa+2R}{6} h_A,
    \end{align}
    so integration by parts gives
    \begin{align}
        \calE[h_0] &= -\int_M \inner{\Delta_f h_0 + 2\Rm(h_0)}{h_0} e^f \\
        &= \calQ_f[\scrS] + \int_M \left( \frac{R-\kappa}{3} |h_{0,H}|^2 + \frac{\kappa+2R}{6} |h_A|^2 \right) e^f.
    \end{align}
    It remains to express the coefficients of $|h_{0,H}|^2$ and $|h_A|^2$ in terms of curvature coefficients. Since $u' = D-B$, $D' = D(A-2B)$, and $B' = -K_2 - B^2$, we compute
    \begin{align}
        \Delta u + |\nabla u|^2 &= u'' + (A + 2B)u' + (u')^2 \\
        &= D(A-2B)+K_2+B^2 + (A+2B)(D-B) + (D-B)^2 \\
        &= K_2 + D^2 - AB + 2D(A-B) \\
        &= K_2 + H_{12} + 2M_2.
    \end{align}
    This and \eqref{eq:scalar-conformal-change} give
    \begin{align}
        \frac{R-\kappa}{3} = 2(K_2 + H_{12} + 2M_2).
    \end{align}
    Finally, using \eqref{eq:appleton-frame-scalar}, we obtain
    \begin{align}
        \frac{\kappa+2R}{6} = \frac{R}{2} - (\Delta u + |\nabla u|^2) = K_1 + K_2 + H_{12} + H_{23} - 2M_2.
    \end{align}
\end{proof}

For the estimates below, we also record $V_H$ and $V_A$ in terms of the variables $v,x,z,Z$ introduced in \S\ref{subsec:new-variables} and \eqref{eq:Z-epsilon-defn}. By \eqref{eq:scaled-curvature-coefficients}, we have
\begin{align}
    b^2 V_H &= 2(-4+3x-v^2 x + vxZ + 2z), \label{eq:scaled-H-A-potentials-1} \\
    b^2 V_A &= a'Z + vxZ - 2x - 2v^2 x - 2z. \label{eq:scaled-H-A-potentials-2}
\end{align}

\subsection{A lower bound for $\calQ_f$}

We first note down an elementary identity:
\begin{lemma} \label{lem:ASD-norm-identity}
    For every ASD $2$-form $\eta$ and real tangent vector $V$, we have $|\iota_V \eta|^2 = \frac{1}{2} |V|^2 |\eta|^2$.
\end{lemma}
\begin{proof}
    If $V=0$ the identity is trivial, so assume $V \neq 0$. Choose an oriented orthonormal frame $\{F_0,F_1,F_2,F_3\}$ such that $F_0 = V/|V|$. As $\eta$ is ASD, its components in this frame satisfy
    \begin{align}
        |\eta_{00}| = 0, \qquad |\eta_{01}| = |\eta_{23}|, \qquad |\eta_{02}| = |\eta_{13}|, \qquad |\eta_{03}| = |\eta_{12}|.
    \end{align}
    Thus, $|\eta|^2 = \sum_{i<j} |\eta_{ij}|^2 = 2\sum_{j=0}^3 |\eta_{0j}|^2$, and hence
    \begin{align}
        |\iota_V \eta|^2 = |V|^2 |\iota_{F_0} \eta|^2 = |V|^2 \sum_{j=0}^3 |\eta_{0j}|^2 = \frac{1}{2} |V|^2 |\eta|^2.
    \end{align}
\end{proof}

The next lemma provides an alternative expression and a lower bound for the $\calQ_f[\scrS]$ term in \eqref{eq:tracefree-energy}.

\begin{lemma} \label{lem:Qf-lower-bound}
    Let $h_0 = \tr(\scrS)$ be a smooth trace-free symmetric $2$-tensor compactly supported in $M^\circ$. Then
    \begin{align} \label{eq:Qf-square-completion}
        \calQ_f[\scrS] &= 2\int_M \left| \delta^L\scrS - \frac{1}{2}\iota_{\nabla f}^L \scrS \right|^2 e^f + \frac{1}{4} \int_M (f')^2 |h_0|^2 e^f.
    \end{align}
    Moreover, decomposing $h_0 = h_{0,H} + h_A$, and writing $h_{0,H} = \tr(\scrS_H)$ where $\scrS_H = \alpha \otimes \hat{\Omega}$ as in \eqref{eq:scrS-decomp}, we have
    \begin{align} \label{eq:Qf-sharper-bound}
        \calQ_f[\scrS] \geq \frac{1}{4} \int_M (f')^2 |h_0|^2 e^f - 4\int_M (B-D)^2|h_A|^2 e^f + \frac{4}{3} \int_M |\delta(e^{f/2}\alpha)|^2.
    \end{align}
\end{lemma}
\begin{proof}
    Let $\star^L$ denote the Hodge star on the left form factor. Since $\star^L \scrS = -\scrS$, we have $\delta^L \scrS = -\star^L d^L \star^L \scrS = \star^L d^L \scrS$, so $|d^L \scrS| = |\delta^L \scrS|$. Consequently,
    \begin{equation}
        |d^L\scrS|^2 + |\delta^L_f \scrS|^2 = |\delta^L \scrS|^2 + |\delta^L \scrS - \iota_{\nabla f}^L \scrS|^2 = 2\left| \delta^L\scrS - \frac{1}{2} \iota_{\nabla f}^L \scrS \right|^2 + \frac{1}{2} |\iota_{\nabla f}^L \scrS|^2. \label{eq:Qf-bound-1}
    \end{equation}
    Applying Lemma \ref{lem:ASD-norm-identity} to each left form component of $\scrS$ gives $|\iota_{\nabla f}^L \scrS|^2 = \frac{1}{2} (f')^2 |\scrS|^2 = \frac{1}{2}(f')^2|h_0|^2$. Combining this with \eqref{eq:Qf-bound-1} and then \eqref{eq:Qf-identity-1} proves \eqref{eq:Qf-square-completion}.

    Next, we will bound the first integral in \eqref{eq:Qf-square-completion} from below by its Hermitian projection. Define $\beta_1 = \frac{1}{\sqrt{2}} (\epsilon^{02} + \epsilon^{13})$ and $\beta_2 = \frac{1}{\sqrt{2}} (\epsilon^{03} - \epsilon^{12})$. Then $\{\hat{\Omega}, \beta_1, \beta_2\}$ is a local orthonormal basis for $\Lambda_{\C}^+$ (cf. \eqref{eq:lambda-plus-onb}). Thus, we can write $h_A = \tr(\scrS_A)$ where
    \begin{align}
        \scrS_A = \alpha_1 \otimes \beta_1 + \alpha_2 \otimes \beta_2, \qquad |h_A|^2 = |\alpha_1|^2 + |\alpha_2|^2,
    \end{align}
    and $\alpha_1,\alpha_2$ are ASD. The connection table \eqref{eq:frame-connection} gives
    \begin{align}
        \nabla_{E_0} \hat{\Omega} = \nabla_{E_1} \hat{\Omega} = 0, \qquad \nabla_{E_2} \hat{\Omega} = -(B-D)\beta_2, \qquad \nabla_{E_3} \hat{\Omega} = (B-D) \beta_1.
    \end{align}
    Denote by $\Pi_H^R$ the right-acting orthogonal projection $\Lambda^1_{\C} \otimes \Lambda^+_\C \to \Lambda^1_{\C} \otimes \C\hat{\Omega}$, where $\Lambda^1_{\C}$ is the bundle of $1$-forms. Using the definition of $\delta^L$ and the fact that $\langle \nabla_{E_i}\hat{\Omega}, \hat{\Omega} \rangle = 0$ for any $E_i$, we compute
    \begin{equation}
        \Pi_H^R(\delta^L \scrS_H) = \Pi_H^R \left( \delta^L(\alpha \otimes \hat{\Omega}) \right) = \Pi_H^R \left( \delta \alpha \otimes \hat{\Omega} - \sum_{i=0}^3 (\iota_{E_i}\alpha) \otimes \nabla_{E_i}\hat{\Omega} \right) = \delta\alpha \otimes \hat{\Omega}. \label{eq:Pi-H-R-Hermitian}
    \end{equation}
    On the other hand, using that $\langle \nabla_{E_i} \beta_n, \hat{\Omega} \rangle = -\langle \beta_n, \nabla_{E_i} \hat{\Omega} \rangle$ for each $n=1,2$,
    \begin{align}
        \Pi_H^R(\delta^L \scrS_A) &= \sum_{n=1}^2 \Pi_H^R \left( \delta \alpha_n \otimes \beta_n - \sum_{i=0}^3 ((\iota_{E_i}\alpha_n) \otimes \nabla_{E_i} \beta_n) \right) \\
        &= \sum_{n=1}^2 \sum_{i=0}^3 (\iota_{E_i} \alpha_n) \otimes \langle \beta_n, \nabla_{E_i}\hat{\Omega} \rangle \hat{\Omega} \\
        &= (B-D)(\iota_{E_3} \alpha_1 - \iota_{E_2} \alpha_2) \otimes \hat{\Omega}. \label{eq:Pi-H-R-skewHermitian}
    \end{align}
    Since $\iota_{\nabla f}^L$ only contracts in the left form factor, and it commutes with $\Pi_H^R$, we have $\Pi_H^R(\iota_{\nabla f}^L \scrS) = \iota_{\nabla f}\alpha \otimes \hat{\Omega}$. Combining this with \eqref{eq:Pi-H-R-Hermitian} and \eqref{eq:Pi-H-R-skewHermitian} gives
    \begin{align}
        \Pi_H^R \left( \delta^L \scrS - \frac{1}{2} \iota_{\nabla f}^L \scrS \right) &= \left( \delta\alpha + (B-D)(\iota_{E_3}\alpha_1 - \iota_{E_2}\alpha_2) - \frac{1}{2} \iota_{\nabla f}\alpha \right) \otimes \hat{\Omega}. \label{eq:Q-projection-formula}
    \end{align}
    For $n=1,2$, anti-self-duality of $\alpha_n$ and Lemma \ref{lem:ASD-norm-identity} imply
    \begin{align}
        |\iota_{E_3}\alpha_1 - \iota_{E_2}\alpha_2|^2 \leq 2(|\iota_{E_3}\alpha_1|^2 + |\iota_{E_2}\alpha_2|^2) = |\alpha_1|^2 + |\alpha_2|^2 = |h_A|^2. \label{eq:alphas-bound}
    \end{align}
    Using \eqref{eq:Q-projection-formula} and \eqref{eq:alphas-bound} in \eqref{eq:Qf-square-completion}, together with the inequality $|x+y|^2 \geq \frac{2}{3}|x|^2 - 2|y|^2$, we get
    \begin{align}
        \calQ_f[\scrS] - \frac{1}{4} \int_M (f')^2 |h_0|^2 e^f &\geq 2\int_M \left| \Pi_H^R\left( \delta^L \scrS - \frac{1}{2} \iota_{\nabla f}^L \scrS \right) \right|^2 e^f \\
        &\geq \frac{4}{3} \int_M \left| \delta\alpha - \frac{1}{2} \iota_{\nabla f}\alpha \right|^2 e^f - 4\int_M (B-D)^2 |h_A|^2 e^f.
    \end{align}
    The bound \eqref{eq:Qf-sharper-bound} follows from this and the calculation $\delta(e^{f/2}\alpha) = e^{f/2}(\delta\alpha - \frac{1}{2} \iota_{\nabla f}\alpha)$.
\end{proof}

We now estimate the divergence term in \eqref{eq:Qf-sharper-bound} in a fixed sector. For this, we use the orthonormal frame $(\sigma_0, \sigma_+, \sigma_-)$ for $\Lambda_{\C}^-$ from \eqref{eq:asd-fiber-weight-basis}. Also recall the variables $v,x,z,Z$ introduced in \S\ref{subsec:new-variables} and \eqref{eq:Z-epsilon-defn}.

\begin{lemma} \label{lem:divergence-weighted-sectorial-asd}
    Let $\eta$ be a smooth ASD $2$-form compactly supported in $M^\circ$ and belonging to a single $(j,q)$-sector, where $(j,q)$ is admissible. Expand
    \begin{align} \label{eq:eta-expansion}
        \eta = u_0 \sigma_0 + u_+ \sigma_+ + u_- \sigma_-.
    \end{align}
    Then
    \begin{align}
        2\int_M |\delta(e^{f/2}\eta)|^2 \geq \int_M b^{-2} \left( d_{j,q}^0 |u_0|^2 + d_{j,q}^+ |u_+|^2 + d_{j,q}^- |u_-|^2 \right) e^f,
    \end{align}
    where
    \begin{align}
        d_{j,q}^0 &= \frac{q^2}{x} + \lambda_{j,q} + 2vxZ - 4z + 4x(1-v-v^2) - 8, \label{eq:djq0-defn} \\
        d_{j,q}^\pm &= \frac{(q \pm a')^2}{x} + \lambda_{j,q} + zZ + 2vxZ - x(3v^2 + 4v + 2) - 2, \label{eq:djqpm-defn} \\
        \lambda_{j,q} &= 2j(j+1)-\frac{q^2}{2}.
    \end{align}
\end{lemma}
\begin{proof}
    Let $\tilde{\eta} = e^{f/2}\eta$. Since $\tilde{\eta}$ is ASD, we have $\delta\tilde{\eta} = -\star d \star \tilde{\eta} = \star d\tilde{\eta}$, so $|d\tilde{\eta}|^2 = |\delta\tilde{\eta}|^2$. Write $\Delta_H = -(d\delta + \delta d)$ for the (unweighted) Hodge Laplacian. By the Bochner formula,
    \begin{align} \label{eq:divergence-equality}
        2 \int_M |\delta\tilde{\eta}|^2 &= \int_M |d\tilde{\eta}|^2 + |\delta\tilde{\eta}|^2 = -\int_M \inner{\Delta_H \tilde{\eta}}{\tilde{\eta}} = \int_M |\nabla\tilde{\eta}|^2 + \inner{\calK\tilde{\eta}}{\tilde{\eta}},
    \end{align}
    where $\calK$ acts on $2$-forms by
    \begin{align}
        (\calK \eta)_{ij} = \sum_k \left( \Ric_{ik} \eta_{kj} + \Ric_{jk} \eta_{ik} \right) - \sum_{kl} \Rm_{ijlk} \eta_{kl}
    \end{align}
    in an orthonormal frame. Dropping the $|\nabla_{\del_s} \tilde{\eta}|^2$ term in \eqref{eq:divergence-equality}, then breaking up the integral into radial and cross-sectional parts using \eqref{eq:dvolg-decomp}, we get
    \begin{align}
        2 \int_M |\delta(e^{f/2}\eta)|^2 &\geq \int_M \sum_{\alpha=1}^3 |\nabla_{E_\alpha} \tilde{\eta}|^2 + \inner{\calK \tilde{\eta}}{\tilde{\eta}} = \int_M \left( \sum_{\alpha=1}^3 |\nabla_{E_\alpha} \eta|^2 + \inner{\calK \eta}{\eta} \right) e^f \\
        &= \int_0^\infty w \int_{L_{2,k}} \sum_{\alpha=1}^3 |\nabla_{E_\alpha} \eta|^2 \, d\omega \, ds + \int_M \inner{\calK \eta}{\eta} e^f. \label{eq:divergence-lb-proof}
    \end{align}

    We first calculate the curvature integral. Substituting \eqref{eq:frame-curvature} and then \eqref{eq:scaled-curvature-coefficients}, we compute
    \begin{align}
        b^2 \calK \sigma_0 &= 2b^2 (K_2 + H_{12} - 2M_2) \sigma_0 \\
        &= 2(-4+3x-v^2 x + vxZ - 2z) \sigma_0, \label{eq:ASD-Bochner-curvature-1} \\
        b^2 \calK \sigma_\pm &= b^2(K_1 + K_2 + H_{12} + H_{23} + 2M_2) \sigma_\pm \\
        &= (a'Z + vxZ - 2x - 2v^2 x + 2z) \sigma_{\pm}. \label{eq:ASD-Bochner-curvature-2}
    \end{align}
    Since $(\sigma_0,\sigma_+,\sigma_-)$ is orthonormal, it follows that
    \begin{align}
        \int_M \inner{\calK\eta}{\eta} e^f = \int_M b^{-2} \Big[& 2(-4+3x-v^2 x + vxZ - 2z) |u_0|^2 \\
        & + (a'Z + vxZ - 2x - 2v^2 x + 2z) (|u_+|^2 + |u_-|^2) \Big] e^f. \label{eq:K-ASD-integral}
    \end{align}

    Next, we estimate the inner integral in the first term of \eqref{eq:divergence-lb-proof}. First, Lemma \ref{lem:differential-identities} gives $X_1 u_0 = -iq u_0$ and $X_1 u_\pm = -i(q \mp 2) u_\pm$, while \eqref{eq:frame-connection} implies $\nabla_{X_1} \sigma_0 = 0$ and $\nabla_{X_1} \sigma_\pm = \mp ia(A+C) \sigma_\pm$. Combining these identities with $a(A+C) = a'+2-x$ gives
    \begin{align}
        \nabla_{X_1} \eta = -iq u_0 \sigma_0 - i(q+a'-x) u_+ \sigma_+ - i(q-a'+x) u_- \sigma_-.
    \end{align}
    Since $E_1 = X_1/a$ and $x = a^2/b^2$, it follows that along each level set of $s$,
    \begin{equation}
        \int_{L_{2,k}} |\nabla_{E_1} \eta|^2 \, d\omega = \frac{1}{b^2} \int_{L_{2,k}} \left[ \frac{q^2}{x} |u_0|^2 + \frac{(q+a'-x)^2}{x} |u_+|^2 + \frac{(q-a'+x)^2}{x} |u_-|^2 \right] \, d\omega. \label{eq:ASD-Hopf-derivative}
    \end{equation}
    Next, using the elementary inequality $|x+y|^2 \geq \frac{1}{2}|x|^2 - |y|^2$,
    \begin{align}
        \sum_{\mu=2}^3 |\nabla_{E_\mu} \eta|^2
        &\geq \frac{1}{2} \sum_{\mu=2}^3 (|E_\mu u_0|^2 + |E_\mu u_+|^2 + |E_\mu u_-|^2) \\
        &\quad - \sum_{\mu=2}^3 |u_0 \nabla_{E_\mu} \sigma_0 + u_+ \nabla_{E_\mu} \sigma_+ + u_- \nabla_{E_\mu} \sigma_-|^2. \label{eq:horizontal-bd}
    \end{align}
    Let $\lambda_{j,q} = 2j(j+1) - q^2/2$. Using Lemma \ref{lem:differential-identities}, the first sum on the right integrates over any level set of $s$ to
    \begin{align}
        &\frac{1}{2} \int_{L_{2,k}} \sum_{\mu=2}^3 (|E_\mu u_0|^2 + |E_\mu u_+|^2 + |E_\mu u_-|^2) \, d\omega \\
        &= \frac{1}{2b^2} \int_{L_{2,k}} \sum_{\mu=2}^3 (|X_\mu u_0|^2 + |X_\mu u_+|^2 + |X_\mu u_-|^2) \, d\omega \\
        &= \frac{1}{b^2} \int_{L_{2,k}} \left[\lambda_{j,q}|u_0|^2 + (\lambda_{j,q} + 2q-2)|u_+|^2 + (\lambda_{j,q} - 2q - 2) |u_-|^2 \right] \, d\omega. \label{eq:term1-div}
    \end{align}
    Meanwhile, the connection table \eqref{eq:frame-connection} allows us to compute
    \begin{align}
        u_0 \nabla_{E_2} \sigma_0 + u_+ \nabla_{E_2} \sigma_+ + u_- \nabla_{E_2} \sigma_- &= -\frac{i(B+D)}{\sqrt{2}} [(u_+ - u_-) \sigma_0 + u_0(\sigma_+ - \sigma_-)], \\
        u_0 \nabla_{E_3} \sigma_0 + u_+ \nabla_{E_3} \sigma_+ + u_- \nabla_{E_3} \sigma_- &= \frac{B+D}{\sqrt{2}} [(u_+ + u_-) \sigma_0 - u_0(\sigma_+ + \sigma_-)].
    \end{align}
    Since $B = b'/b$, $D = Q/b$, and $b' = vQ$, we have $b^2(B+D)^2 = (b'+Q)^2 = x(1+v)^2$. Thus,
    \begin{equation} \label{eq:term2-div}
        \sum_{\mu=2}^3 |u_0 \nabla_{E_\mu} \sigma_0 + u_+ \nabla_{E_\mu} \sigma_+ + u_- \nabla_{E_\mu} \sigma_-|^2 = b^{-2} x(1+v)^2 \left( 2|u_0|^2 + |u_+|^2 + |u_-|^2 \right).
    \end{equation}
    Using \eqref{eq:term1-div} and \eqref{eq:term2-div} in \eqref{eq:horizontal-bd}, we conclude that over any level set of $s$,
    \begin{align}
        \int_{L_{2,k}} \sum_{\mu=2}^3 |\nabla_{E_\mu} \eta|^2 \, d\omega \geq \frac{1}{b^2} \int_{L_{2,k}} \Big[& (\lambda_{j,q} - 2x(1+v)^2) |u_0|^2 \\
        & + (\lambda_{j,q} + 2q - 2 - x(1+v)^2) |u_+|^2 \\
        & + (\lambda_{j,q} - 2q - 2 - x(1+v)^2) |u_-|^2 \Big] \, d\omega. \label{eq:ASD-integrated-horizontal-bound}
    \end{align}
    By \eqref{eq:ASD-Hopf-derivative} and \eqref{eq:ASD-integrated-horizontal-bound}, we have
    \begin{align}
        \int_0^\infty w \int_{L_{2,k}} \sum_{\alpha=1}^3 |\nabla_{E_\alpha} \eta|^2 \, d\omega \, ds \\
        \geq \int_M b^{-2} \Bigg[& \left( \frac{q^2}{x} + \lambda_{j,q} - 2x(1+v)^2 \right) |u_0|^2 \\
        & + \left( \frac{(q+a'-x)^2}{x} + \lambda_{j,q} + 2q - 2 - x(1+v)^2 \right) |u_+|^2 \\
        & + \left( \frac{(q-a'+x)^2}{x} + \lambda_{j,q} - 2q - 2 - x(1+v)^2 \right) |u_-|^2 \Bigg] e^f.
    \end{align}
    Inserting this and \eqref{eq:K-ASD-integral} back into \eqref{eq:divergence-lb-proof}, then simplifying the coefficients using $a' = z + vx$, the asserted inequality follows.
\end{proof}

\subsection{Scalar coefficient bounds} \label{subsec:scalar-coeff-bounds}

In this subsection, we will bound the coefficients that appear in the stability proof below. Recall the quantities $V_H$ and $V_A$ defined in Proposition \ref{prop:tracefree-energy}, as well as the quantities $d_{j,q}^0$ and $d_{j,q}^\pm$ defined in Lemma \ref{lem:divergence-weighted-sectorial-asd}.

\begin{lemma} \label{lem:modified-potential-bounds}
    For every $s>0$, we have
    \begin{align}
        \tilde{V}_H &:= b^2 \left( V_H + \frac{(f')^2}{4} \right) \geq 4 + \frac{10}{9}x, \label{eq:tildeVH-bound} \\
        \tilde{V}_A &:= b^2 \left( V_A + \frac{(f')^2}{4} - 4(B-D)^2 \right) \geq \frac{49}{36}. \label{eq:tildeVA-bound}
    \end{align}
    Moreover, if $(j,q)$ is an admissible pair such that $q \neq 0$ and $j \geq |q|/2$, we have at every $s>0$
    \begin{align}
        \tilde{V}_H + \frac{2}{3} d_{j,q}^0 &\geq \frac{2q^2}{3x} + \frac{xZ^2}{4} - \frac{22}{3} + \frac{160x}{27}, \label{eq:tildeVH-0-bound} \\
        \tilde{V}_H + \frac{2}{3} d_{j,q}^\pm &\geq \frac{40}{9}. \label{eq:tildeVH-pm-bound}
    \end{align}
\end{lemma}
\begin{proof}
    Using \eqref{eq:scaled-H-A-potentials-1}, \eqref{eq:scaled-H-A-potentials-2}, and the identities $a' = z + vx$ and $bf' = \sqrt{x}(Z-2v)$, we have
    \begin{align}
        \tilde{V}_H &= 2(-4 + 3x - v^2 x + vxZ + 2z) + \frac{x}{4}(Z-2v)^2, \label{eq:VH-expanded} \\
        b^2 \left( V_A + \frac{(f')^2}{4} \right) &= z(Z-2) + 2x(vZ-1-v^2) + \frac{x}{4} (Z-2v)^2. \label{eq:VA-expanded}
    \end{align}
    Lemma \ref{lem:sharp-background-bounds} gives $z>k(1-x) \geq 3(1-x)$ and $vZ > 2$. Also recall from Lemma \ref{lem:comparison-monotonicity} that $v < 2/k \leq 2/3$. Thus, $Z > 2/v \geq 3$. Using these bounds in \eqref{eq:VA-expanded}, we get
    \begin{align}
        b^2 \left( V_A + \frac{(f')^2}{4} \right) &\geq 3(1-x) + 2x(1-v^2) + x(v^{-1} - v)^2 \\
        &= 3(1-x) + x(v^{-2} - v^2). \label{eq:Va-preliminary}
    \end{align}
    Since $B = b'/b$, $D = Q/b$, and $b' = vQ$, we have $b^2(B-D)^2 = (b'-Q)^2 = x(1-v)^2$. Subtracting $4b^2(B-D)^2 = 4x(1-v)^2$ from \eqref{eq:Va-preliminary} gives
    \begin{align} \label{eq:VA-first-lb}
        \tilde{V}_A \geq 3(1-x) + x \underbrace{(v^{-2} - 5v^2 + 8v - 4)}_{=: F(v)}.
    \end{align}
    We have $F'(v) = -2v^{-3} - 10v + 8$, and $F''(v) = 6v^{-4} - 10 > 0$ on $0 < v \leq 2/3$. Since $F'(2/3) = -65/12 < 0$, it follows that $F$ is decreasing on this interval. Thus $F(v) \geq F(2/3) = 49/36$, and \eqref{eq:VA-first-lb} gives
    \begin{align}
        \tilde{V}_A \geq 3(1-x) + \frac{49}{36}x \geq \frac{49}{36}.
    \end{align}
    Also, using $y=2z+2vx+xZ$, we can rewrite \eqref{eq:VH-expanded} as
    \begin{align} \label{eq:VH-4plusx}
        \tilde{V}_H = 2y - 8 + x\left[ \left( \frac{Z}{2} + v - 2 \right)^2 + 2(1-v^2) \right] \geq 4+\frac{10}{9}x,
    \end{align}
    where the inequality uses $y \geq 2k \geq 6$ and $v \leq 2/3$. This proves \eqref{eq:tildeVH-bound} and \eqref{eq:tildeVA-bound}.

    Write $\lambda_{j,q} = 2j(j+1) - q^2/2$. Since $q \neq 0$, the admissibility condition \eqref{eq:jq-admissibility-conditions} gives $|q| \geq k \geq 3$. Since $j \geq |q|/2$,
    \begin{align} \label{eq:charged-lambda-lower-bound}
        \lambda_{j,q} \geq 2\frac{|q|}{2} \left( \frac{|q|}{2} + 1 \right) - \frac{q^2}{2} = |q| \geq 3.
    \end{align}
    Using \eqref{eq:VH-expanded} and the definition \eqref{eq:djq0-defn} of $d_{j,q}^0$, we obtain
    \begin{align}
        \tilde{V}_H + \frac{2}{3} d_{j,q}^0 &= \frac{2q^2}{3x} + \frac{xZ^2}{4} + \frac{2}{3} \lambda_{j,q} + \frac{4}{3}z - \frac{40}{3} + x\left( \frac{26}{3} - \frac{11}{3}v^2 + \frac{7}{3}vZ - \frac{8}{3}v \right).
    \end{align}
    Substituting \eqref{eq:charged-lambda-lower-bound}, $z \geq 3(1-x)$, and $vZ \geq 2$, and simplifying, this gives
    \begin{align}
        \tilde{V}_H + \frac{2}{3} d_{j,q}^0 &\geq \frac{2q^2}{3x} + \frac{xZ^2}{4} - \frac{22}{3} + x\left( \frac{28}{3} - \frac{11}{3}v^2 - \frac{8}{3}v \right).
    \end{align}
    Since $0 < v \leq 2/3$, we have $\frac{28}{3} - \frac{11}{3}v^2 - \frac{8}{3}v \geq \frac{160}{27}$. Using this above proves \eqref{eq:tildeVH-0-bound}.

    Discarding the nonnegative square in the definition \eqref{eq:djqpm-defn} of $d_{j,q}^\pm$, then using $zZ \geq 9(1-x)$ and $vZ \geq 2$, we get
    \begin{align}
        d_{j,q}^{\pm} \geq \lambda_{j,q} + 7 - x(7+3v^2+4v).
    \end{align}
    Since $v \leq 2/3$, we have $7+3v^2+4v \leq 11$. Together with \eqref{eq:charged-lambda-lower-bound}, this gives $d_{j,q}^\pm \geq 10-11x$. Finally, applying \eqref{eq:VH-4plusx} and using $x<1$ gives
    \begin{align}
        \tilde{V}_H + \frac{2}{3} d_{j,q}^\pm \geq 4 + \frac{10}{9}x + \frac{2}{3}(10-11x) = \frac{32}{3} - \frac{56}{9}x \geq \frac{40}{9},
    \end{align}
    proving \eqref{eq:tildeVH-pm-bound}.
\end{proof}

\begin{lemma} \label{lem:rho-lemma}
    Define $\rho = b^2 f'\frac{Q'}{Q} = z(Z-2v)$. Then $-1 < \rho-\xi < kZ$ for all $s>0$.
\end{lemma}
\begin{proof}
    Since $f',Q',Q>0$, we have $\rho > 0$. Lemma \ref{lem:sharp-background-bounds} gives $0<\xi<1$, so $\rho-\xi > -\xi > -1$. Using $\xi = 4-x+xv^2-vy$, $y=2z+2vx+xZ$, and $a'=z+xv$, we compute
    \begin{align}
        \rho-\xi &= z(Z-2v) - (4-x+xv^2-v(2z+2vx+xZ)) \\
        &= a'Z - 4 + x(1+v^2). \label{eq:rho-epsilon-identity}
    \end{align}
    Since $0<x<1$ and $v \leq 2/3$, we have $x(1+v^2)<4$, and hence $\rho-\xi < a'Z$. Finally, $Z>0$ and Lemma \ref{lem:4d-concavity} gives $a' < k$, so $\rho-\xi < kZ$.
\end{proof}

\subsection{Proof of Theorem \ref{thm:charged-coercivity}} \label{subsec:charged-proof}

We are now ready to prove the main result of this section.

\begin{proof}[Proof of Theorem \ref{thm:charged-coercivity}]
    Let $h$ be a smooth symmetric $2$-tensor compactly supported in $M^\circ$ and belonging to an admissible $(j,q)$-sector with $q \neq 0$. As in \eqref{eq:charged-tensor-decomposition}, decompose
    \begin{align}
        h = \phi\sfE + h_0, \qquad \sfE = \frac{g}{2}, \qquad h_0 = h_{0,H} + h_A.
    \end{align}
    As explained at the end of \S\ref{subsec:conformally-kahler}, $\phi, h_{0,H}, h_A$ all belong to the $(j,q)$-sector of the respective bundles.

    First suppose $j < |q|/2$. From \eqref{eq:L2-functions-identification}, the scalar $(j,q)$-sector vanishes. Thus $\phi = 0$ and $h = h_0$. Proposition \ref{prop:tracefree-energy}, \eqref{eq:Qf-square-completion}, and Lemma \ref{lem:modified-potential-bounds} give
    \begin{align}
        \calE[h] &\geq \int_M \left[ \left( V_H + \frac{(f')^2}{4} \right) |h_{0,H}|^2 + \left( V_A + \frac{(f')^2}{4} \right) |h_A|^2 \right] e^f \\
        &\geq \int_M b^{-2} \left( \tilde{V}_H |h_{0,H}|^2 + \tilde{V}_A |h_A|^2 \right) e^f \\
        &\geq \int_M b^{-2} \left( 4|h_{0,H}|^2 + \frac{49}{36} |h_A|^2 \right) e^f \\
        &\geq \frac{49}{36} \int_M \frac{|h|^2}{b^2} e^f,
    \end{align}
    proving the theorem in this case.

    Henceforth suppose $j \geq |q|/2$. We first decompose $\calE[h]$. Since $\sfE$ is parallel and $h_0$ is trace-free, $\inner{\nabla(\phi\sfE)}{\nabla h_0} = 0$. Also, $\Rm(\sfE) = {\Ric}/2$. Expanding $\calE[h]$ and using self-adjointness of $\Rm$ therefore gives
    \begin{align} \label{eq:charged-energy-decomposition}
        \calE[h] &= \calE[\phi\sfE] + \calE[h_0] - 2\Re\int_M \phi\inner{\Ric_0}{h_0} e^f,
    \end{align}
    where $\Ric_0 = \Ric - \frac{R}{4}g$ is the trace-free Ricci tensor. We estimate the two diagonal terms and then control the cross-term.

    \vspace*{1em}
    \noindent\textbf{Estimating $\calE[\phi\sfE]$:}
    \vspace*{1em}

    Since $|\sfE| = 1$, $\nabla\sfE = 0$, and $\Rm(\sfE) = {\Ric}/2$, we have
    \begin{align}
        \calE[\phi\sfE] = \int_M \left( |\nabla\phi|^2 - \frac{R}{2}|\phi|^2 \right) e^f.
    \end{align}
    Using $\Delta f = R$ and integrating by parts, we compute
    \begin{align}
        \int_M |\nabla\phi + \phi\nabla f|^2 e^f &= \int_M \left( |\nabla\phi|^2 + \inner{\nabla|\phi|^2}{\nabla f} + |\phi|^2 |\nabla f|^2 \right) e^f = \int_M \left( |\nabla\phi|^2 - R|\phi|^2 \right) e^f.
    \end{align}
    Consequently,
    \begin{align}
        \calE[\phi\sfE] &= \int_M \left( |\nabla\phi + \phi\nabla f|^2 + \frac{R}{2} |\phi|^2 \right) e^f \\
        &= \int_M \left( |\del_s \phi + f'\phi|^2 + \sum_{\alpha=1}^3 |E_\alpha \phi|^2 + \frac{R}{2} |\phi|^2 \right) e^f \\
        &\geq \int_M \left( \sum_{\alpha=1}^3 |E_\alpha \phi|^2 + \frac{R}{2}|\phi|^2 \right) e^f, \label{eq:traceful-energy}
    \end{align}
    Define
    \begin{align}
        \Lambda_{j,q} := 4j(j+1) + (x^{-1}-1)q^2.
    \end{align}
    Since $\phi$ belongs to the scalar $(j,q)$-sector, integration by parts on each level set of $s$ and Corollary \ref{cor:scalar-eigenvalues} give
    \begin{align} \label{eq:IBP-scalar}
        \int_M \sum_{\alpha=1}^3 |E_\alpha \phi|^2 e^f = \int_M \frac{\Lambda_{j,q}}{b^2} |\phi|^2 e^f.
    \end{align}
    Meanwhile, \eqref{eq:scaled-curvature-coefficients} and \eqref{eq:rho-epsilon-identity} give
    \begin{align}
        \frac{b^2 R}{2} &= b^2(K_1 + 2K_2 + 2H_{12} + H_{23}) \\
        &= a'Z + 2vxZ - 4 + x - 3v^2 x \\
        &= \rho-\xi + 2vx(Z-2v) \\
        &= \rho-\xi + 2bb'f' \geq \rho-\xi. \label{eq:b-squared-R-estimate}
    \end{align}
    Combining \eqref{eq:IBP-scalar} and \eqref{eq:b-squared-R-estimate} into \eqref{eq:traceful-energy} yields
    \begin{align} \label{eq:charged-trace-energy-bound}
        \calE[\phi\sfE] \geq \int_M b^{-2} \left( \Lambda_{j,q} + \rho-\xi \right) |\phi|^2 e^f.
    \end{align}

    \vspace*{1em}
    \noindent\textbf{Estimating $\calE[h_0]$:}
    \vspace*{1em}

    Recall the local orthonormal basis \eqref{eq:asd-fiber-weight-basis} of ASD $2$-forms $(\sigma_0,\sigma_+,\sigma_-)$. Define the unit trace-free Hermitian tensor
    \begin{align}
        \sfT := \tr(\sigma_0 \otimes \hat{\Omega}) = \frac{1}{2} \left( (\epsilon^0)^2 + (\epsilon^1)^2 - (\epsilon^2)^2 - (\epsilon^3)^2 \right)
    \end{align}
    and the unit skew-Hermitian tensor $\tau = ((\epsilon^0)^2 - (\epsilon^1)^2)/\sqrt{2}$. Indeed, $\tau$ is skew-Hermitian since $J^*$ interchanges $(\epsilon^0)^2$ and $(\epsilon^1)^2$, so $J^*\tau = -\tau$. Decompose pointwise
    \begin{alignat}{2}
        h_{0,H} &= \psi\sfT + w_H, \qquad & w_H &\perp \sfT, \\
        h_A &= \chi\tau + w_A, \qquad & w_A &\perp \tau.
    \end{alignat}
    These are orthogonal decompositions, so
    \begin{align} \label{eq:charged-component-norm}
        |h|^2 = |\phi|^2 + |\psi|^2 + |w_H|^2 + |\chi|^2 + |w_A|^2.
    \end{align}

    Write $h_{0,H} = \tr(\alpha \otimes \hat{\Omega})$. As explained shortly after \eqref{eq:scrS-decomp}, $\alpha$ belongs to the $(j,q)$-sector of ASD $2$-forms. Since the trace map is an isometry, the basis expansion
    \begin{align}
        \alpha = \psi \sigma_0 + u_+ \sigma_+ + u_- \sigma_-
    \end{align}
    has the same coefficient $\psi$, and $|w_H|^2 = |u_+|^2 + |u_-|^2$. Proposition \ref{prop:tracefree-energy} and \eqref{eq:Qf-sharper-bound} give
    \begin{align}
        \calE[h_0] \geq \int_M b^{-2} \left( \tilde{V}_H |\alpha|^2 + \tilde{V}_A |h_A|^2 \right) e^f + \frac{4}{3} \int_M |\delta(e^{f/2}\alpha)|^2.
    \end{align}
    Applying Lemma \ref{lem:divergence-weighted-sectorial-asd} therefore yields
    \begin{align}
        \calE[h_0] \geq \int_M b^{-2} \Bigg[
            &\left( \tilde{V}_H + \frac{2}{3} d_{j,q}^0 \right) |\psi|^2
            + \left( \tilde{V}_H + \frac{2}{3} d_{j,q}^+ \right) |u_+|^2 \\
            &+ \left( \tilde{V}_H + \frac{2}{3} d_{j,q}^- \right) |u_-|^2
            + \tilde{V}_A |h_A|^2
            \Bigg] e^f.
    \end{align}
    Since $j \geq |q|/2$, we may apply all the bounds in Lemma \ref{lem:modified-potential-bounds} to conclude
    \begin{equation}
        \calE[h_0] \geq \int_M b^{-2} \Bigg[
            \left( \frac{2q^2}{3x} + \frac{xZ^2}{4} - \frac{22}{3} + \frac{160x}{27} \right) |\psi|^2 + \frac{40}{9} |w_H|^2 + \frac{49}{36} \left( |\chi|^2 + |w_A|^2 \right)
        \Bigg] \, e^f. \label{eq:charged-tracefree-energy-bound}
    \end{equation}

    \vspace*{1em}
    \noindent\textbf{Expressing the cross-term:}
    \vspace*{1em}

    Finally, we express the cross-term in \eqref{eq:charged-energy-decomposition}. By \eqref{eq:appleton-frame-Ricci},
    \begin{align}
        \Ric_0 = \Ric - \frac{R}{4}g = (K_1 - H_{23}) \sfT + \sqrt{2}(K_2 - H_{12}) \tau.
    \end{align}
    Using \eqref{eq:scaled-curvature-coefficients} and \eqref{eq:rho-epsilon-identity}, we have $b^2(K_1-H_{23}) = \rho-\xi$ and $b^2(K_2 - H_{12}) = -\xi$. Thus,
    \begin{align} \label{eq:charged-Ricci-decomposition}
        b^2 \Ric_0 = (\rho-\xi)\sfT - \sqrt{2}\xi \tau,
    \end{align}
    and hence
    \begin{align} \label{eq:charged-cross-term}
        -2\Re\int_M \phi \inner{\Ric_0}{h_0}e^f = \int_M b^{-2} \left[ -2(\rho-\xi) \Re(\phi\bar{\psi}) + 2\sqrt{2}\xi \Re(\phi\bar{\chi}) \right] e^f.
    \end{align}

    \vspace*{1em}
    \noindent\textbf{Combining the estimates and concluding:}
    \vspace*{1em}

    Combining the estimates \eqref{eq:charged-trace-energy-bound}, \eqref{eq:charged-tracefree-energy-bound}, and \eqref{eq:charged-cross-term} into \eqref{eq:charged-energy-decomposition}, we obtain
    \begin{align} \label{eq:charged-matrix-energy}
        \calE[h] \geq \int_M b^{-2} \left[ \bfu^* \bfM \bfu + \frac{40}{9} |w_H|^2 + \frac{49}{36} |w_A|^2 \right] e^f,
    \end{align}
    where
    \begin{gather}
        \bfu = \begin{pmatrix} \phi \\ \psi \\ \chi \end{pmatrix}, \qquad \bfM = \begin{pmatrix}
            \Lambda_{j,q}+\rho-\xi & -(\rho-\xi) & \sqrt{2}\xi \\
            -(\rho-\xi) & \calB + \frac{1}{4} & 0 \\
            \sqrt{2}\xi & 0 & \frac{49}{36}
        \end{pmatrix}, \\
        \calB := \frac{2q^2}{3x} + \frac{xZ^2}{4} - \frac{91}{12} + \frac{160x}{27}.
    \end{gather}
    We will show that $\bfM - \frac{1}{4}I$ is nonnegative definite.

    Since the admissibility condition \eqref{eq:jq-admissibility-conditions} gives $|q| \geq k \geq 3$, the AM-GM inequality gives
    \begin{align}
        \calB - \left( \frac{3q^2}{8x} + \frac{xZ^2}{4} \right) &= \frac{7q^2}{24x} + \frac{160x}{27} - \frac{91}{12} \\
        &\geq \frac{21}{8x} + \frac{160x}{27} - \frac{91}{12} \\
        &\geq 2\sqrt{\frac{21}{8} \frac{160}{27}} - \frac{91}{12} > 0. \label{eq:B-positive}
    \end{align}
    In particular, $\calB > 0$.
    The lower-right $2 \times 2$ block of $\bfM - \frac{1}{4} I$ is therefore the positive definite matrix $\mathrm{diag}(\calB, 10/9)$. By the Schur complement criterion, $\bfM - \frac{1}{4} I \succeq 0$ if and only if
    \begin{align} \label{eq:charged-Schur-complement}
        S := \Lambda_{j,q} + \rho - \xi - \frac{1}{4} - \frac{(\rho-\xi)^2}{\calB} - \frac{9}{5}\xi^2 \geq 0.
    \end{align}

    Recall that $v \leq 2/3$ and $vZ > 2$, so $Z \geq 3$ and $kZ > 1$. Lemma \ref{lem:rho-lemma} then gives $|\rho-\xi| < kZ \leq |q|Z$. Combining this with \eqref{eq:B-positive} gives
    \begin{align}
        \calB \geq \frac{3q^2}{8x} + \frac{x(\rho-\xi)^2}{4q^2}.
    \end{align}
    We use the elementary calculus inequality
    \begin{align} \label{eq:charged-elementary-inequality}
        \frac{t^2}{3/8 + t^2/4} \leq 1+t
    \end{align}
    for all $t \geq 0$. Taking $t=x|\rho-\xi|/q^2$ in \eqref{eq:charged-elementary-inequality} gives
    \begin{align}
        \frac{(\rho-\xi)^2}{\calB} \leq \frac{q^2}{x} \frac{t^2}{3/8 + t^2/4} \leq \frac{q^2}{x} + |\rho-\xi|.
    \end{align}
    Moreover, $j \geq |q|/2$ implies $\Lambda_{j,q} = \frac{q^2}{x} + 4j(j+1) - q^2 \geq \frac{q^2}{x} + 2|q|$. Substituting these bounds and $0<\xi<1$ into \eqref{eq:charged-Schur-complement}, we obtain
    \begin{align}
        S &\geq 2|q| + (\rho-\xi - |\rho-\xi|) - \frac{1}{4} - \frac{9}{5} \geq 2|q| - 2 - \frac{41}{20} \geq \frac{39}{20} > 0,
    \end{align}
    where we used $\rho-\xi > -1$ and $|q| \geq 3$. Thus the Schur complement criterion is met, and $\bfM \succeq \frac{1}{4}I$. Applying this in \eqref{eq:charged-matrix-energy} and using \eqref{eq:charged-component-norm} now gives
    \begin{align}
        \calE[h] \geq \frac{1}{4} \int_M \frac{|h|^2}{b^2} e^f,
    \end{align}
    completing the proof of Theorem \ref{thm:charged-coercivity}.
\end{proof}

\section{The stability inequality in $q=0$ sectors} \label{sec:neutral-stability}

In this section, we establish the stability estimate \eqref{eq:sector-coercivity-hypothesis} in admissible sectors with $q=0$.
\begin{theorem} \label{thm:neutral-coercivity}
    There exists $\mu > 0$ such that every smooth symmetric $2$-tensor $h$ compactly supported in $M^\circ$ and belonging to an admissible $(j,0)$-sector satisfies
    \begin{align}
        \calE[h] \geq \mu \int_M \frac{|h|^2}{b^2} e^f.
    \end{align}
\end{theorem}

The estimates of the previous section do not suffice to establish the stability inequality in $q=0$ sectors. We therefore compute the energy $\calE[h]$ directly. To motivate the method, consider the decomposition of Euclidean Dirichlet energy into spherical harmonics. If $\phi(r,\theta) = u(r)Y(\theta)$ is a smooth function compactly supported in $\R^n \setminus \{0\}$, where
\begin{align}
    -\Delta_{\Sph^{n-1}} Y = \lambda Y, \qquad \norm{Y}_{L^2(\Sph^{n-1})} = 1,
\end{align}
then
\begin{align} \label{eq:dirichlet-energy-decomp}
    \int_{\R^n} |\nabla\phi|^2 \, dx = \int_0^\infty r^{n-1} \left( |u'|^2 + \frac{\lambda}{r^2} |u|^2 \right) \, dr.
\end{align}
Using spherical tensor harmonics on $\CP^1 \cong \Sph^2$, we obtain an analogous expression for $\calE[h]$, where $h$ belongs to an admissible $(j,0)$-sector. The expression is a radial quadratic form whose zeroth-order part is governed by an explicit Hermitian matrix $\sfC_j(s)$. The remaining task is then to establish coercivity of these quadratic forms uniformly in $j$; see Corollary \ref{cor:Cj-positivity-implies-stability} and its proof for the precise estimate needed.

\subsection{Hopf-invariant tensors} \label{subsec:hopf-invariant-tensors}

We first describe symmetric $2$-tensors in a $(j,0)$-sector using their radial, Hopf, and horizontal components.

Write $\bar{g} = g_{\mathrm{FS}}$. Extending the Hopf projection $\pi: L_{2,k} \to \CP^1$ to $M^\circ$ via $\pi(s,x) = \pi(x)$, define $g_H = b^2 \pi^*\bar{g}$. With $\epsilon^0 = ds$ and $\epsilon^1 = a\sigma$ as in \S\ref{subsec:local-frame}, the soliton metric is
\begin{align}
    g = (\epsilon^0)^2 + (\epsilon^1)^2 + g_H.
\end{align}
For a $1$-form $\alpha$ and a symmetric $2$-tensor $T$ on $\CP^1$, define their normalized horizontal lifts to $M^\circ$ by
\begin{align}
    \alpha^H = b \, \pi^*\alpha, \qquad T^H = b^2 \, \pi^*T.
\end{align}
These lifts are Hopf-invariant. Moreover, the normalizations give
\begin{align}
    |\alpha^H|_g^2 = \pi^*(|\alpha|_{\bar{g}}^2), \qquad |T^H|_g^2 = \pi^*(|T|_{\bar{g}}^2).
\end{align}
If $\alpha$ and $T$ are independent of $s$, then $\nabla_{\del_s} \alpha^H = 0$ and $\nabla_{\del_s} T^H = 0$. Indeed, write $\pi^*\alpha = c_2 \sigma_2 + c_3 \sigma_3$ where $c_2,c_3$ are independent of $s$. Using that $\epsilon^\mu = b \sigma_\mu$ for $\mu=2,3$, we have $\alpha^H = c_2 \epsilon^2 + c_3 \epsilon^3$, so $\nabla_{\del_s} \alpha^H = 0$ by \eqref{eq:frame-connection}. The proof for $T^H$ is similar.

By Lemma \ref{lem:hopf-invariance}, every symmetric $2$-tensor on $M^\circ$ belonging to a $(j,0)$-sector admits a unique decomposition
\begin{align} 
    h ={}& u_1 (\epsilon^0)^2 + \sqrt{2} \, u_2 \, \epsilon^0 \odot \epsilon^1 + u_3 (\epsilon^1)^2 + \frac{u_4}{\sqrt{2}} g_H \\
    &+ \sqrt{2} \, \epsilon^0 \odot \alpha^H + \sqrt{2} \, \epsilon^1 \odot \beta^H + T^H, \label{eq:neutral-base-decomposition}
\end{align}
where $u_1,\ldots,u_4$ are $s$-dependent functions on $\CP^1$, $\alpha,\beta$ are $s$-dependent $1$-forms on $\CP^1$, and $T$ is an $s$-dependent trace-free symmetric $2$-tensor on $\CP^1$. The summands in \eqref{eq:neutral-base-decomposition} are pointwise orthogonal with
\begin{align}
    |h|_g^2 = \sum_{n=1}^4 |u_n|^2 + |\alpha|_{\bar{g}}^2 + |\beta|_{\bar{g}}^2 + |T|_{\bar{g}}^2.
\end{align}

\subsection{A basis of the $(j,0)$-sector via spherical tensor harmonics}

In this subsection, we construct a convenient basis for any given $(j,0)$-sector of symmetric $2$-tensors, then reduce Theorem \ref{thm:neutral-coercivity} to the task of proving a radial coercivity estimate.

Orient $(\CP^1,\bar{g})$ so that $\pi^*\dvol_{\bar{g}} = \sigma_2 \wedge \sigma_3$. The bracket relations \eqref{eq:T-brackets} then imply
\begin{align} \label{eq:dsigma-base}
    d\sigma = -2\pi^* \dvol_{\bar{g}}.
\end{align}
Let $\bar{\nabla}$ and $\bar{\star}$ denote the Levi-Civita connection and Hodge star of $(\CP^1, \bar{g})$, and write $\bar{\Delta} = \tr_{\bar{g}} \bar{\nabla}^2$ for the rough Laplacian on tensors.

Fix an integer $j \geq 0$. Recall from \eqref{eq:spherical-harmonics-decomp} that $\bar{\calH}_j$ is the $(2j+1)$-dimensional space of eigenfunctions of $-\bar{\Delta}$ with eigenvalue $4j(j+1)$. Set
\begin{align}
    \lambda = \lambda_j = 4j(j+1), \qquad c_k = \sqrt{\frac{k}{2\pi}}.
\end{align}
Then $-\bar{\Delta} Y = \lambda Y$ for $Y \in \bar{\calH}_j$, and by Lemma \ref{lem:hopf-projection-harmonics}, the map $Y \mapsto c_k \pi^*Y$ identifies $\bar{\calH}_j$ unitarily and $G$-equivariantly with the scalar $(j,0)$-sector on $L_{2,k}$.

The construction of $1$-form and symmetric tensor harmonics from $Y \in \bar{\calH}_j$ is well-known; see for instance \cite{boucetta-sphere}. If $j \geq 1$, define the $1$-forms
\begin{align} \label{eq:1form-harmonics}
    \alpha_E[Y] = \frac{dY}{\sqrt{\lambda}}, \qquad \alpha_B[Y] = \frac{\bar{\star} dY}{\sqrt{\lambda}}.
\end{align}
If $j \geq 2$, also define $\kappa = \sqrt{(\lambda-8)/2}$ and the trace-free symmetric $2$-tensors
\begin{align} \label{eq:sym2-harmonics}
    T_E[Y] = \frac{\bar{\nabla}^2 Y + \frac{\lambda}{2} Y \bar{g}}{\kappa \sqrt{\lambda}}, \qquad T_B[Y] = \frac{\operatorname{sym}\bar{\nabla}(\bar{\star} dY)}{\kappa\sqrt{\lambda}}.
\end{align}
Here $(\operatorname{sym} \bar{\nabla} \alpha)_{\mu\nu} = \frac{1}{2} (\bar{\nabla}_\mu \alpha_\nu + \bar{\nabla}_\nu \alpha_\mu)$.

Recall that $G$ acts isometrically on $(\CP^1,\bar{g})$, and this action is orientation-preserving. Thus, the action of $G$ commutes with the operations $d$, $\bar{\star}$, and $\bar{\nabla}$, so for all $\phi \in G$,
\begin{gather}
    \alpha_E[\phi^*Y] = \phi^*\alpha_E[Y], \qquad \alpha_B[\phi^*Y] = \phi^*\alpha_B[Y], \label{eq:alpha-G-equivariant} \\
     T_E[\phi^*Y] = \phi^*T_E[Y], \qquad T_B[\phi^*Y] = \phi^*T_B[Y]. \label{eq:T-G-equivariant}
\end{gather}

Using these tensors, we define families of symmetric $2$-tensors on $M^\circ$ associated to a given $Y \in \bar{\calH}_j$. Here and below, pullbacks of scalar functions to $M^\circ$ are implicit. Define
\begin{alignat}{2}
    U_1[Y] &= c_k Y (\epsilon^0)^2,
    \qquad & U_2[Y] &= \sqrt{2}\, c_k Y
        \epsilon^0 \odot \epsilon^1, \\
    U_3[Y] &= c_k Y (\epsilon^1)^2,
    \qquad & U_4[Y] &= \frac{c_k Y}{\sqrt{2}} g_H.
\end{alignat}
If $j \geq 1$, also define
\begin{alignat}{2}
    V_1[Y] &= \sqrt{2}\, c_k
        \epsilon^0 \odot \alpha_E[Y]^H,
    \qquad & V_2[Y] &= \sqrt{2}\, c_k
        \epsilon^0 \odot \alpha_B[Y]^H, \\
    V_3[Y] &= \sqrt{2}\, c_k
        \epsilon^1 \odot \alpha_E[Y]^H,
    \qquad & V_4[Y] &= \sqrt{2}\, c_k
        \epsilon^1 \odot \alpha_B[Y]^H.
\end{alignat}
If $j \geq 2$, also define
\begin{align}
    W_1[Y] = c_k T_E[Y]^H, \qquad W_2[Y] = c_k T_B[Y]^H.
\end{align}

To facilitate later computations, we record some well-known identities for the above tensors. These follow from commuting covariant derivatives, applying Bochner formulas, and using standard differential-geometric identities. Thus we omit the proof.
\begin{lemma} \label{lem:base-harmonic-identities}
    The following identities involving $\alpha_P$ hold for $j \geq 1$, and those involving $T_P$ hold for $j \geq 2$.
    For $Y,Z \in \bar{\calH}_j$ and $P,Q \in \{E,B\}$, we have
    \begin{align}
        \inner{\alpha_P[Y]}{\alpha_Q[Z]}_{L^2(\CP^1,\bar{g})} &= \delta_{PQ} \inner{Y}{Z}_{L^2(\CP^1,\bar{g})}, \label{eq:alpha-orthog} \\
        \inner{T_P[Y]}{T_Q[Z]}_{L^2(\CP^1,\bar{g})} &= \delta_{PQ} \inner{Y}{Z}_{L^2(\CP^1,\bar{g})}. \label{eq:T-orthog}
    \end{align}
    Moreover, abbreviating $\alpha_P := \alpha_P[Y]$ and $T_P = T_P[Y]$,
    \begin{gather}
        -\bar{\Delta} \alpha_P = (\lambda-4) \alpha_P, \qquad -\bar{\Delta} T_P = (\lambda-16) T_P \\
        \div_{\bar{g}} \alpha_E = -\sqrt{\lambda} \, Y, \qquad \div_{\bar{g}} \alpha_B = 0, \qquad \div_{\bar{g}} T_P = -\kappa \alpha_P, \label{eq:div-identities} \\
        \operatorname{sym}\bar{\nabla} \alpha_E = \kappa T_E - \frac{\sqrt{\lambda}}{2} Y \bar{g}, \qquad \operatorname{sym}\bar{\nabla} \alpha_B = \kappa T_B. \label{eq:sym-alpha-identities} 
    \end{gather}
\end{lemma}

Collect the tensors defined above into the ordered family
\begin{align}
    (\tau_1[Y], \ldots, \tau_{n_j}[Y]) = (U_1[Y], \ldots, U_4[Y], V_1[Y], \ldots, V_4[Y], W_1[Y], W_2[Y]),
\end{align}
omitting absent groups if $j=0$ or $j=1$. Thus, $n_j = 4,8,10$ if $j = 0,1, \geq 2$ respectively.

\begin{proposition} \label{prop:spherical-sector-basis}
    Fix $s>0$, and let $\{Y_m\}_{m=-j}^j$ be an orthonormal basis of $\bar{\calH}_j$. Then
    \begin{align} \label{eq:tau-family}
        \{ \tau_n[Y_m](s,\cdot) : 1 \leq n \leq n_j, \, -j \leq m \leq j \}
    \end{align}
    is an orthonormal basis of the $(j,0)$-sector in $L^2(\scrE_s, d\omega)$, where $\scrE = S^2 T^*_{\C} M^\circ$. Moreover, each $\tau_n[Y_m]$ is parallel along $\del_s$.
\end{proposition}
\begin{proof}
    The action of $G$ on $(M^\circ, g)$ preserves the tensors $\epsilon^0$, $\epsilon^1$, and $g_H$. The Hopf map $\pi$ is also $G$-equivariant, so normalized horizontal lifting $(\cdot)^H$ commutes with the $G$-action. Together with \eqref{eq:alpha-G-equivariant} and \eqref{eq:T-G-equivariant}, it follows that each map $\bar{\calH}_j \to C^\infty(\scrE_s)$, $Y \mapsto \tau_n[Y](s,\cdot)$ is $G$-equivariant. Since $\bar{\calH}_j \cong_G V_{j,0}$, Schur's lemma and Proposition \ref{prop:sector-peter-weyl} imply that each map has image contained in the $(j,0)$-sector of $\scrE_s$.

    Since the Hopf fibers on $L_{2,k}$ have length $2\pi/k$ with respect to the round metric of curvature $1$, every smooth function $\phi$ on $\CP^1$ satisfies
    \begin{align} \label{eq:hopf-fiber-integration}
        c_k^2 \int_{L_{2,k}} \pi^*\phi \, d\omega = \int_{\CP^1} \phi \, \dvol_{\bar{g}}.
    \end{align}
    Combining this with the orthogonality of tensor types in \eqref{eq:neutral-base-decomposition}, and the identities \eqref{eq:alpha-orthog} and \eqref{eq:T-orthog}, we can check case by case that at any fixed $s>0$,
    \begin{align} \label{eq:160861}
        \int_{L_{2,k}} \inner{\tau_{n'}[Y](s,\cdot)}{\tau_n[Z](s,\cdot)} \, d\omega = \delta_{nn'} \inner{Y}{Z}_{L^2(\CP^1,\bar{g})},
    \end{align}
    for all $Y,Z \in \bar{\calH}_j$ and $1 \leq n,n' \leq n_j$.
    
    Meanwhile, Proposition \ref{prop:sector-peter-weyl} gives that the $(j,0)$-sector of $\scrE_s$ has dimension $\dim (V_{j,0} \otimes \calM_{j,0}(F)) = (2j+1)n_j$, where $F$ is the fiber of $\scrE_s$ above the basepoint $o$. Here we have used that $\dim V_{j,0} = 2j+1$, while \eqref{eq:symmetric-tensor-dims} and \eqref{eq:symmetric-tensor-multiplicity} imply $\dim \calM_{j,0}(F) = n_j$. Together with \eqref{eq:160861}, it follows that \eqref{eq:tau-family} defines an orthonormal basis of the $(j,0)$-sector of $\scrE_s$.

    Finally, since $\epsilon^0$, $\epsilon^1$, $g_H$, and the normalized lifts $(\cdot)^H$ of $s$-independent base tensors are parallel along $\del_s$, the same holds true for every $\tau_n[Y_m]$.
\end{proof}

Next, we show that to establish the stability inequality in a given $(j,0)$-sector, it suffices to consider tensors built out of one fixed $Y \in \bar{\calH}_j$. For fixed $s>0$, let $h$ and $k$ be smooth sections of $\scrE_s = S^2 T^*_{\C}M^\circ|_{\Sigma_s}$. Let $g_s$ be the induced metric on $\Sigma_s$. Choose any local $g_s$-orthonormal frame $(F_1,F_2,F_3)$ on $\Sigma_s$, and define the $G$-invariant pairing
\begin{align}
    \fka_s(h,k) := \int_{L_{2,k}} \left( \sum_{\alpha=1}^3 \inner{\nabla_{F_\alpha} h}{\nabla_{F_\alpha} k} - 2\inner{\Rm(h)}{k} \right) \, d\omega.
\end{align}
Here $\nabla$ and $\Rm$ are the ambient connection and curvature action, evaluated at $\Sigma_s$. The expression does not depend on the choice of $g_s$-orthonormal frame.

Since the maps $Y \mapsto \tau_n[Y]$ are $G$-equivariant, for every $n,n'$ the sesquilinear form on $\bar{\calH}_j$ defined by $(Y,Z) \mapsto \fka_s(\tau_{n'}[Y], \tau_n[Z])$ is $G$-invariant. Since $\bar{\calH}_j \cong V_{j,0}$ is irreducible, Schur's lemma implies that there is a Hermitian matrix $\sfC_j(s)$ such that
\begin{align} \label{eq:C-matrix-defining}
    \fka_s(\tau_{n'}[Y], \tau_n[Z]) = (\sfC_j(s))_{nn'} \inner{Y}{Z}_{L^2(\CP^1,\bar{g})}.
\end{align}
In particular, the entries of $\sfC_j(s)$ may be recovered by the formula
\begin{align} \label{eq:C-matrix-explicit}
    (\sfC_j(s))_{nn'} = \int_{L_{2,k}} \left( \sum_{\alpha=1}^3 \inner{\nabla_{F_\alpha} \tau_{n'}[Y]}{\nabla_{F_\alpha} \tau_n[Y]} - 2\inner{\Rm(\tau_{n'}[Y])}{\tau_n[Y]} \right) \, d\omega
\end{align}
for any choice of $Y \in \bar{\calH}_j$ such that $\norm{Y}_{L^2(\CP^1,\bar{g})} = 1$. We may then approach proving Theorem \ref{thm:neutral-coercivity} as follows.

\begin{corollary} \label{cor:Cj-positivity-implies-stability}
    Fix an integer $j \geq 0$. Suppose there exists $\mu > 0$ such that for every $X \in C^\infty_c((0,\infty); \C^{n_j})$, we have
    \begin{align} \label{eq:radial-coercivity-target}
        \int_0^\infty w\left( |X'|^2 + X^* \sfC_j X \right) \, ds \geq \mu \int_0^\infty \frac{w}{b^2} |X|^2 \, ds.
    \end{align}
    Then every smooth symmetric 2-tensor $h$ compactly supported in $M^\circ$ and belonging to the $(j,0)$-sector satisfies
    \begin{align}
        \calE[h] \geq \mu \int_M \frac{|h|^2}{b^2} e^f.
    \end{align}
\end{corollary}
\begin{proof}
    Choose an orthonormal basis $\{Y_m\}_{m=-j}^j$ of $\bar{\calH}_j$. By Proposition \ref{prop:spherical-sector-basis}, every smooth tensor $h$ in the $(j,0)$-sector has a unique expansion
    \begin{align} \label{eq:h-expansion-neutral}
        h(s,\cdot) = \sum_{m=-j}^j \sum_{n=1}^{n_j} X_{m,n}(s) \tau_n[Y_m].
    \end{align}
    Write $X_m = (X_{m,1}, \ldots, X_{m,n_j})^T$. If $h$ is compactly supported in $M^\circ$, then each $X_m$ is smooth and compactly supported in $(0,\infty)$. As the basis from Proposition \ref{prop:spherical-sector-basis} is orthonormal and parallel along $\del_s$, \eqref{eq:h-expansion-neutral} implies
    \begin{align}
        \int_{L_{2,k}} |h(s,\cdot)|^2 \, d\omega = \sum_{m=-j}^j |X_m(s)|^2, \qquad \int_{L_{2,k}} |\nabla_{\del_s} h(s,\cdot)|^2 \, d\omega = \sum_{m=-j}^j |X_m'(s)|^2.
    \end{align}
    Also, \eqref{eq:C-matrix-defining}, \eqref{eq:C-matrix-explicit}, and \eqref{eq:h-expansion-neutral} give $\mathfrak{a}_s(h,h) = \sum_{m=-j}^j X_m^* \sfC_j(s) X_m$. Thus,
    \begin{align}
        \calE[h] &= \int_0^\infty w \int_{L_{2,k}} \left( |\nabla h|^2 - 2\inner{\Rm(h)}{h} \right) \, d\omega \, ds \\
        &= \int_0^\infty w \left( \int_{L_{2,k}} |\nabla_{\del_s}h|^2 \, d\omega + \mathfrak{a}_s(h,h) \right) \, ds \\
        &= \sum_{m=-j}^j \int_0^\infty w\left( |X_m'|^2 + X_m^* \sfC_j X_m \right) \, ds. \label{eq:energy-formula-radial}
    \end{align}
    Using the assumption \eqref{eq:radial-coercivity-target}, we get
    \begin{align}
        \calE[h] \geq \mu \sum_{m=-j}^j \int_0^\infty \frac{w}{b^2} |X_m|^2 \, ds = \mu \int_M \frac{|h|^2}{b^2} e^f.
    \end{align}
\end{proof}

\subsection{Calculation of $\sfC_j$} \label{subsec:Cj-calculation}

In this subsection, we calculate the matrix $\sfC_j(s)$. We fix $Y \in \bar{\calH}_j$ with $\norm{Y}_{L^2(\CP^1,\bar{g})} = 1$, so that $\sfC_j(s)$ can be computed using \eqref{eq:C-matrix-explicit}. Abbreviate $\tau_n := \tau_n[Y]$ and accordingly $U_n := U_n[Y]$, $V_n := V_n[Y]$, $W_n := W_n[Y]$. The result is as follows.

\begin{lemma} \label{lem:spherical-angular-matrix}
    Put $\lambda=4j(j+1)$. If $j\geq2$, then the matrix $\sfC_j$ in the ordered basis
    \[
        (U_1,U_2,U_3,U_4,V_1,V_2,V_3,V_4,W_1,W_2)
    \]
    is given by
    \begin{align} \label{eq:spherical-angular-matrix}
        \sfC_j=
        \begin{pmatrix}
            \dfrac{\lambda}{b^2}I_4+\sfC_U
                &\dfrac{2\sqrt{\lambda}}b \sfM&0\\[2mm]
            \dfrac{2\sqrt{\lambda}}b \sfM^T
                &\dfrac{\lambda-4}{b^2}I_4+\sfC_V
                &\dfrac{2\sqrt{\lambda-8}}b \sfN\\[2mm]
            0&\dfrac{2\sqrt{\lambda-8}}b \sfN^T
                &\dfrac{\lambda-8}{b^2}I_2
        \end{pmatrix},
    \end{align}
    where
    \begin{gather}
        \sfC_U=
        \begin{pmatrix}
            2A^2+4B^2
                &0
                &-2(A^2+K_1)
                &-2\sqrt2(B^2+K_2)\\
            0
                &4A^2+2B^2+2D^2+2K_1
                &0&0\\
            -2(A^2+K_1)
                &0
                &2A^2+4D^2
                &-2\sqrt2(D^2+H_{12})\\
            -2\sqrt2(B^2+K_2)
                &0
                &-2\sqrt2(D^2+H_{12})
                &2B^2+2D^2-2H_{23}
        \end{pmatrix}, \\
        \sfC_V=
        \begin{pmatrix}
            c_1 & 0 & 0 & c_3 \\
            0 & c_1 & -c_3 & 0 \\
            0 & -c_3 & c_2 & 0 \\
            c_3 & 0 & 0 & c_2 
        \end{pmatrix},
        \quad
        \sfM=
        \begin{pmatrix}
            -\sqrt2B&0&0&0\\
            0&D&-B&0\\
            0&0&0&\sqrt2D\\
            B&0&0&-D
        \end{pmatrix},
        \quad
        \sfN=
        \begin{pmatrix}
            -B&0\\
            0&-B\\
            0&D\\
            -D&0
        \end{pmatrix},
    \end{gather}
    and
    \begin{align}
        c_1 &= A^2 + 5B^2 + 2D^2 + 2K_2, \\
        c_2 &= A^2 + B^2 + 6D^2 + 2H_{12}, \\
        c_3 &= 2D(A+B)-6M_2.
    \end{align}
    The quantities $A,B,D,K_1,K_2,H_{12},H_{23},M_2$ are defined in \S\ref{subsec:local-frame}.
    If $j=1$, omit the rows and columns corresponding to $(W_1,W_2)$, and take $\lambda=8$. If $j=0$, also omit the rows and columns corresponding to $(V_1,\ldots,V_4)$, and take $\lambda=0$. Thus $\sfC_0=\sfC_U$.
\end{lemma}

To calculate $\sfC_j(s)$, we will separately calculate the Hermitian matrices
\begin{align}
    (\sfG_j(s))_{nn'} &:= \int_{L_{2,k}} \sum_{\alpha=1}^3 \inner{\nabla_{F_\alpha} \tau_{n'}}{\nabla_{F_\alpha} \tau_n} \, d\omega, \label{eq:sfG-defn} \\
    (\sfK_j(s))_{nn'} &:= \int_{L_{2,k}} \inner{\Rm(\tau_{n'})}{\tau_n} \, d\omega, \label{eq:sfK-defn}
\end{align}
so that $\sfC_j(s) = \sfG_j(s) - 2\sfK_j(s)$. Let us first introduce the local frame used in the calculations.

Fix $\bar{p} \in \CP^1$ and choose a local oriented $\bar{g}$-orthonormal coframe $\bar{\epsilon}^2,\bar{\epsilon}^3$, with dual frame $\bar{E}_2,\bar{E}_3$, such that $\bar{\nabla} \bar{\epsilon}^\mu|_{\bar{p}} = 0$ for $\mu=2,3$. Upstairs in $M^\circ$, define the local orthonormal coframe
\begin{align}
    \vartheta^0 = \epsilon^0 = ds, \qquad \vartheta^1 = \epsilon^1 = a\sigma, \qquad \vartheta^\mu = (\bar{\epsilon}^\mu)^H = b \, \pi^*\bar{\epsilon}^\mu \quad (\mu=2,3).
\end{align}
Fix $p \in M^\circ$ above $\bar{p}$. Rotating $(\bar{\epsilon}^2,\bar{\epsilon}^3)$ if necessary, we may assume that $\vartheta^i|_p = \epsilon^i|_p$ for all $i$, where $(\epsilon^i)$ is the coframe introduced in \eqref{eq:appleton-epsilon-frame}.

Let $\tilde{E}_\mu$ be the horizontal lift of $\bar{E}_\mu$. The dual orthonormal frame is
\begin{align}
    F_0 = E_0 = \del_s, \qquad F_1 = E_1, \qquad F_\mu = \frac{1}{b} \tilde{E}_\mu \quad (\mu=2,3).
\end{align}
Recall the functions $A=a'/a$, $B=b'/b$, and $D=a/b^2$. We compute $[F_1,F_2] = [F_1,F_3] = 0$, $[F_0,F_1] = -AF_1$, and $[F_0,F_\mu] = -BF_\mu$ for $\mu=2,3$. Also, using \eqref{eq:dsigma-base}, we compute $[\tilde{E}_2,\tilde{E}_3] = \widetilde{[\bar{E}_2,\bar{E}_3]} + 2a E_1 = 2aE_1$, where the last equality holds at $p$. Applying the Koszul formula, we obtain the following connection table at the point $p$:
\begin{equation}\label{eq:base-lift-connection}
    \begin{array}{c|cccc}
        &\vartheta^0&\vartheta^1&\vartheta^2&\vartheta^3\\ \hline
        \nabla_{F_0}
        &0&0&0&0\\[2mm]
        \nabla_{F_1}
        &A\vartheta^1&-A\vartheta^0
        &-D\vartheta^3&D\vartheta^2\\[2mm]
        \nabla_{F_2}
        &B\vartheta^2&-D\vartheta^3
        &-B\vartheta^0&D\vartheta^1\\[2mm]
        \nabla_{F_3}
        &B\vartheta^3&D\vartheta^2
        &-D\vartheta^1&-B\vartheta^0
    \end{array}.
\end{equation}
Since $ds$ and $\sigma$ are Hopf-invariant, and $\vartheta^2,\vartheta^3$ are defined by pullback by $\pi$, it follows that the coframe $(\vartheta^i)$ is Hopf-invariant. As a consequence, suppose $\alpha$ is a Hopf-invariant $1$-form and $h$ is a Hopf-invariant symmetric $2$-tensor (which is true if they belong to a $(j,0)$-sector of the respective bundle). If they are written locally as
\begin{align}
    \alpha = \sum_{i=0}^3 \alpha_i \vartheta^i, \qquad h = \sum_{i,j=0}^3 h_{ij} \vartheta^i \otimes \vartheta^j,
\end{align}
then $F_1(\alpha_i) = 0$ and $F_1(h_{ij}) = 0$. Finally, if $\phi$ is a function on $\CP^1$, then $F_\mu(\pi^*\phi) = b^{-1} \pi^*(\bar{E}_\mu \phi)$ for $\mu=2,3$.

\subsubsection{Calculation of $\sfG_j$}

We calculate the orbital derivatives of $U_1,V_1,W_1$, and then calculate the entries of $\sfG_j$ at indices $(1,1), (5,5), (9,9), (1,5), (1,9), (5,9)$. The other entries follow similarly.

For these calculations, suppose $j \geq 2$, and abbreviate $\alpha = \alpha_E[Y]$ and $T = T_E[Y]$. Thus,
\begin{align}
    U_1 = c_k Y (\vartheta^0)^2, \qquad V_1 = \sqrt{2} \, c_k \vartheta^0 \odot \alpha^H, \qquad W_1 = c_k T^H.
\end{align}
All calculations below take place at $p$. Pullbacks of scalar functions to $M^\circ$ are implicit. Write
\begin{align}
    Y_\mu = \bar{E}_\mu Y, \qquad \alpha = \alpha_2 \bar{\epsilon}^2 + \alpha_3 \bar{\epsilon}^3, \qquad T_\mu = T(\bar{E}_\mu, \cdot) \quad (\mu=2,3).
\end{align}

We start with $U_1$. The connection table \eqref{eq:base-lift-connection} and the product rule give
\begin{align}
    \nabla_{F_1} U_1 &= 2c_k AY \, \vartheta^0 \odot \vartheta^1, \label{eq:U1-deriv-1} \\
    \nabla_{F_\mu} U_1 &= c_k \left( \frac{Y_\mu}{b} (\vartheta^0)^2 + 2BY \, \vartheta^0 \odot \vartheta^\mu \right), \quad \mu=2,3. \label{eq:U1-deriv-2}
\end{align}
For $V_1$, first write $\alpha^H = b \, \pi^*(\alpha_2\bar{\epsilon}^2 + \alpha_3 \bar{\epsilon}^3) = \alpha_2 \vartheta^2 + \alpha_3 \vartheta^3$, where we have identified $\alpha_\mu$ with its pullback. Since $\alpha^H$ is Hopf-invariant, $F_1(\alpha_\mu) = 0$ and $F_\nu(\alpha_\mu) = b^{-1} \bar{E}_\nu(\alpha_\mu)$ for $\mu, \nu \in \{2,3\}$. Thus, we calculate
\begin{align}
    \nabla_{F_1} \alpha^H &= D(\alpha_3 \vartheta^2 - \alpha_2 \vartheta^3), \\
    \nabla_{F_2} \alpha^H &= \frac{1}{b}(\bar{\nabla}_{\bar{E}_2}\alpha)^H - B \alpha_2 \vartheta^0 + D \alpha_3 \vartheta^1, \\
    \nabla_{F_3} \alpha^H &= \frac{1}{b} (\bar{\nabla}_{\bar{E}_3} \alpha)^H - B\alpha_3 \vartheta^0 - D\alpha_2 \vartheta^1.
\end{align}
It follows that
\begin{align}
    \nabla_{F_1} V_1 &= \sqrt{2} c_k \left[ A \vartheta^1 \odot \alpha^H + D \vartheta^0 \odot (\alpha_3 \vartheta^2 - \alpha_2 \vartheta^3) \right], \label{eq:V1-deriv-1} \\
    \nabla_{F_2} V_1 &= \sqrt{2} c_k \left[ B\vartheta^2 \odot \alpha^H + \frac{1}{b} \vartheta^0 \odot (\bar{\nabla}_{\bar{E}_2}\alpha)^H - B \alpha_2 (\vartheta^0)^2 + D \alpha_3 \vartheta^0 \odot \vartheta^1 \right], \label{eq:V1-deriv-2} \\
    \nabla_{F_3} V_1 &= \sqrt{2} c_k \left[ B \vartheta^3 \odot \alpha^H + \frac{1}{b} \vartheta^0 \odot (\bar{\nabla}_{\bar{E}_3} \alpha)^H - B \alpha_3 (\vartheta^0)^2 - D \alpha_2 \vartheta^0 \odot \vartheta^1 \right]. \label{eq:V1-deriv-3}
\end{align}
For $W_1$, since $T$ is $\bar{g}$-trace-free and hence $T^H$ is horizontal and $g$-trace-free, we can write
\begin{align}
    T^H = T_{22} \left( (\vartheta^2)^2 - (\vartheta^3)^2 \right) + 2T_{23} \, \vartheta^2 \odot \vartheta^3
\end{align}
where $T_{22}$ and $T_{23}$ are Hopf-invariant functions. Differentiating gives
\begin{align}
    \nabla_{F_1} W_1 &= 2c_k D \left[ T_{23} \left( (\vartheta^2)^2 - (\vartheta^3)^2 \right) - 2 T_{22} \, \vartheta^2 \odot \vartheta^3 \right], \label{eq:W1-deriv-1} \\
    \nabla_{F_2} W_1 &= c_k \left[ \frac{1}{b} (\bar{\nabla}_{\bar{E}_2} T)^H - 2B \vartheta^0 \odot T_2^H + 2D \vartheta^1 \odot T_3^H \right], \label{eq:W1-deriv-2} \\
    \nabla_{F_3} W_1 &= c_k \left[ \frac{1}{b}(\bar{\nabla}_{\bar{E}_3} T)^H - 2B \vartheta^0 \odot T_3^H - 2D \vartheta^1 \odot T_2^H \right]. \label{eq:W1-deriv-3}
\end{align}

Using these derivatives, we may calculate the entries $(\sfG_j)_{11}$, $(\sfG_j)_{55}$, $(\sfG_j)_{99}$, $(\sfG_j)_{15}$, $(\sfG_j)_{19}$, $(\sfG_j)_{59}$. Since $\sfG_j$ is Hermitian, this also gives entries below the diagonal. We first obtain the diagonal entries. For a $1$-form $\eta$ on $\CP^1$, we calculate
\begin{align} \label{eq:mixed-norm-squared}
    |\vartheta^0 \odot \eta^H|^2 = |\vartheta^1 \odot \eta^H|^2 = \frac{1}{2}|\eta|_{\bar{g}}^2, \qquad \sum_{\mu=2}^3 |\vartheta^\mu \odot \eta^H|^2 = \frac{3}{2} |\eta|_{\bar{g}}^2.
\end{align}
Using this with the displayed equations between \eqref{eq:U1-deriv-1} and \eqref{eq:W1-deriv-3}, we calculate
\begin{align}
    c_k^{-2} \sum_{\alpha=1}^3 |\nabla_{F_\alpha} U_1|^2 &= \frac{|dY|^2}{b^2} + (2A^2 + 4B^2)|Y|^2, \\
    c_k^{-2} \sum_{\alpha=1}^3 |\nabla_{F_\alpha} V_1|^2 &= \frac{|\bar{\nabla}\alpha|^2}{b^2} + (A^2 + 5B^2 + 2D^2) |\alpha|^2, \\
    c_k^{-2} \sum_{\alpha=1}^3 |\nabla_{F_\alpha} W_1|^2 &= \frac{|\bar{\nabla} T|^2}{b^2} + (2B^2 + 6D^2) |T|^2.
\end{align}
Lemma \ref{lem:base-harmonic-identities} and integration by parts give
\begin{gather}
    \norm{Y}_{L^2(\CP^1,\bar{g})}^2 = \norm{\alpha}_{L^2(\CP^1,\bar{g})}^2 = \norm{T}_{L^2(\CP^1,\bar{g})}^2 = 1, \\
    \norm{dY}_{L^2(\CP^1,\bar{g})}^2 = \lambda, \qquad \norm{\bar{\nabla}\alpha}_{L^2(\CP^1,\bar{g})}^2 = \lambda-4, \qquad \norm{\bar{\nabla} T}_{L^2(\CP^1,\bar{g})}^2 = \lambda-16.
\end{gather}
Integrating their pullbacks over $L_{2,k}$ with respect to $d\omega$ acquires an extra $c_k^{-2}$ factor, see \eqref{eq:hopf-fiber-integration}. Thus, by the definition \eqref{eq:sfG-defn} of $\sfG_j$, and the preceding calculations, we obtain
\begin{align}
    (\sfG_j)_{11} &= \frac{\lambda}{b^2} + 2A^2 + 4B^2, \\
    (\sfG_j)_{55} &= \frac{\lambda-4}{b^2} + A^2 + 5B^2 + 2D^2, \\
    (\sfG_j)_{99} &= \frac{\lambda-16}{b^2} + 2B^2 + 6D^2.
\end{align}

We next calculate the mixed entries. For the pair $(U_1,V_1)$, the above calculations give
\begin{align}
    c_k^{-2} \inner{\nabla_{F_1} V_1}{\nabla_{F_1} U_1} &= 0, \\
    c_k^{-2} \inner{\nabla_{F_\mu} V_1}{\nabla_{F_\mu} U_1} &= \frac{\sqrt{2}B}{b} \left[ -\alpha_\mu \overline{Y_\mu} + (\bar{\nabla}_{\bar{E}_\mu} \alpha)(\bar{E}_\mu) \overline{Y} \right] \quad (\mu=2,3).
\end{align}
Summing and integrating, then using \eqref{eq:div-identities} and $dY = \alpha\sqrt{\lambda}$, we obtain
\begin{align}
    (\sfG_j)_{15} &= \int_{G/H} \sum_{\alpha=1}^3 \inner{\nabla_{F_\alpha} V_1}{\nabla_{F_\alpha} U_1} \, d\omega \\
    &= \frac{\sqrt{2} B}{b} \left( \inner{\div_{\bar{g}} \alpha}{Y}_{L^2(\CP^1,\bar{g})} - \inner{\alpha}{dY}_{L^2(\CP^1,\bar{g})} \right) = -\frac{2\sqrt{2\lambda} \, B}{b}.
\end{align}
For the pair $(U_1, W_1)$, we have $\inner{\nabla_{F_1} W_1}{\nabla_{F_1} U_1} = 0$ and
\begin{align}
    c_k^{-2} \sum_{\mu=2}^3 \inner{\nabla_{F_\mu} W_1}{\nabla_{F_\mu} U_1} &= -2B^2 \bar{Y} \sum_{\mu=2}^3 T(\bar{E}_\mu, \bar{E}_\mu) = 0,
\end{align}
since $T$ is $\bar{g}$-trace-free. Thus $$(\sfG_j)_{19} = 0.$$
For the pair $(V_1,W_1)$, we have $\inner{\nabla_{F_1} W_1}{\nabla_{F_1} V_1} = 0$, while for $\mu=2,3$,
\begin{align}
    c_k^{-2} \inner{\nabla_{F_\mu} W_1}{\nabla_{F_\mu} V_1} &= \sqrt{2} \inner{\frac{1}{b} (\bar{\nabla}_{\bar{E}_\mu} T)^H - 2B \vartheta^0 \odot T_\mu^H}{B \vartheta^\mu \odot \alpha^H + \frac{1}{b} \vartheta^0 \odot (\bar{\nabla}_{\bar{E}_\mu} \alpha)^H} \\
    &= \frac{\sqrt{2}B}{b} \left( \inner{\bar{\nabla}_{\bar{E}_\mu} T}{\bar{\epsilon}^\mu \odot \alpha} - \inner{T_\mu}{\bar{\nabla}_{\bar{E}_\mu} \alpha} \right).
\end{align}
Summing over $\mu$, we get
\begin{align}
    \sum_{\mu=2}^3 \inner{\bar{\nabla}_{\bar{E}_\mu} T}{\bar{\epsilon}^\mu \odot \alpha} &= \inner{\div_{\bar{g}} T}{\alpha}, \\
    \sum_{\mu=2}^3 \inner{T_\mu}{\bar{\nabla}_{\bar{E}_\mu} \alpha} &= \inner{T}{\operatorname{sym}\bar{\nabla}\alpha}.
\end{align}
Consequently, using \eqref{eq:sym-alpha-identities} and \eqref{eq:div-identities} together with the trace-freeness of $T$,
\begin{align}
    (\sfG_j)_{59} &= \frac{\sqrt{2}B}{b} \left( \inner{\div_{\bar{g}} T}{\alpha}_{L^2(\CP^1,\bar{g})} - \inner{T}{\operatorname{sym}\bar{\nabla}\alpha}_{L^2(\CP^1,\bar{g})} \right) = -\frac{2\sqrt{\lambda-8} \, B}{b}.
\end{align}

The remaining entries are obtained by the same procedure. Collecting the calculations, for $j \geq 2$ we obtain
\begin{align} \label{eq:Gj-matrix}
    \sfG_j=
    \begin{pmatrix}
        G_U & \dfrac{2\sqrt{\lambda}}b \sfM & 0\\[2mm]
        \dfrac{2\sqrt{\lambda}}b \sfM^T
            &G_V&\dfrac{2\sqrt{\lambda-8}}b \sfN\\[2mm]
        0&\dfrac{2\sqrt{\lambda-8}}b \sfN^T&G_W
    \end{pmatrix},
\end{align}
where $\sfM, \sfN$ are as in Lemma \ref{lem:spherical-angular-matrix}, and the diagonal blocks are
\begin{gather}
    G_U
    =\frac{\lambda}{b^2}I_4+
    \begin{pmatrix}
        2A^2+4B^2&0&-2A^2&-2\sqrt2 B^2\\
        0&4A^2+2B^2+2D^2&0&0\\
        -2A^2&0&2A^2+4D^2&-2\sqrt2 D^2\\
        -2\sqrt2 B^2&0&-2\sqrt2 D^2&2B^2+2D^2
    \end{pmatrix}, \\
    G_V
    =\frac{\lambda-4}{b^2}I_4+
    \begin{pmatrix}
        g_1&0&0&g_3\\
        0&g_1&-g_3&0\\
        0&-g_3&g_2&0\\
        g_3&0&0&g_2
    \end{pmatrix}, \qquad
    G_W
    =\left(
        \frac{\lambda-16}{b^2}+2B^2+6D^2
        \right)I_2, \\
    g_1=A^2+5B^2+2D^2,\qquad
    g_2=A^2+B^2+6D^2,\qquad
    g_3=2D(A+B).
\end{gather}
If $j=1$, omit the rows and columns corresponding to $(W_1,W_2)$ and take $\lambda=8$. If $j=0$, also omit the rows and columns corresponding to $(V_1,\ldots,V_4)$ and take $\lambda=0$, so that $\sfG_j = G_U$.

\subsubsection{Calculation of $\sfK_j$}

Recall that $\Rm(h)_{ij} = \Rm_{kijl} h^{kl}$. Also, recall that the coframe $(\vartheta^i)$ agrees with $(\epsilon^i)$ at the point $p \in M^\circ$. Thus, the curvature formulas \eqref{eq:frame-curvature} apply in the $\vartheta$-frame at $p$. We calculate at this point
\begin{align}
    \Rm((\epsilon^0)^2) &= K_1(\epsilon^1)^2 + K_2 g_H, \\
    \Rm(\epsilon^0 \odot \epsilon^1) &= -K_1 \epsilon^0 \odot \epsilon^1, \\
    \Rm((\epsilon^1)^2) &= K_1(\epsilon^0)^2 + H_{12} g_H, \\
    \Rm(g_H) &= 2K_2(\epsilon^0)^2 + 2H_{12}(\epsilon^1)^2 + H_{23} g_H.
\end{align}
Applying the normalizations defining the $U_n$, we obtain
\begin{align}
    \Rm(U_1) &= K_1 U_3 + \sqrt{2} K_2 U_4, \\
    \Rm(U_2) &= -K_1 U_2, \\
    \Rm(U_3) &= K_1 U_1 + \sqrt{2} H_{12} U_4, \\
    \Rm(U_4) &= \sqrt{2} K_2 U_1 + \sqrt{2} H_{12} U_3 + H_{23} U_4.
\end{align}
Next, we compute
\begin{align}
    \Rm(\epsilon^0 \odot \vartheta^2) &= -K_2 \epsilon^0 \odot \vartheta^2 + 3M_2 \epsilon^1 \odot \vartheta^3, \\
    \Rm(\epsilon^0 \odot \vartheta^3) &= -K_2 \epsilon^0 \odot \vartheta^3 - 3M_2 \epsilon^1 \odot \vartheta^2, \\
    \Rm(\epsilon^1 \odot \vartheta^2) &= -H_{12} \epsilon^1 \odot \vartheta^2 - 3M_2 \epsilon^0 \odot \vartheta^3, \\
    \Rm(\epsilon^1 \odot \vartheta^3) &= -H_{12} \epsilon^1 \odot \vartheta^3 + 3M_2 \epsilon^0 \odot \vartheta^2.
\end{align}
If $\alpha$ is a $1$-form on $\CP^1$, write $\alpha^H = \alpha_2 \vartheta^2 + \alpha_3 \vartheta^3$. Using the above calculations, and that $\bar{\star} \bar{\epsilon}^2 = \bar{\epsilon}^3$, $\bar{\star} \bar{\epsilon}^3 = -\bar{\epsilon}^2$ (as $(\bar{\epsilon}^2,\bar{\epsilon}^3)$ is positively oriented by construction), it follows that
\begin{align}
    \Rm(\epsilon^0 \odot \alpha^H) &= -K_2 \epsilon^0 \odot \alpha^H + 3M_2 \epsilon^1 \odot (\bar{\star}\alpha)^H, \\
    \Rm(\epsilon^1 \odot \alpha^H) &= -H_{12} \epsilon^1 \odot \alpha^H - 3M_2 \epsilon^0 \odot (\bar{\star}\alpha)^H.
\end{align}
Since $\bar{\star} \alpha_E = \alpha_B$ and $\bar{\star} \alpha_B = -\alpha_E$, applying the above to $V_1,\ldots,V_4$ gives
\begin{align}
    \Rm(V_1) &= -K_2 V_1 + 3M_2 V_4, \\
    \Rm(V_2) &= -K_2 V_2 - 3M_2 V_3, \\
    \Rm(V_3) &= -H_{12} V_3 - 3M_2 V_2, \\
    \Rm(V_4) &= -H_{12} V_4 + 3M_2 V_1.
\end{align}
If $T$ is a $\bar{g}$-trace-free symmetric 2-tensor on $\CP^1$, we can write $T^H = T_{22} \left( (\vartheta^2)^2 - (\vartheta^3)^2 \right) + 2T_{23} \vartheta^2 \odot \vartheta^3$, where $T_{22}$ and $T_{23}$ are functions. Then we calculate
\begin{align}
    \Rm(T^H) = -H_{23} T^H.
\end{align}
Applying this to the tensors $W_n$ gives
\begin{align}
    \Rm(W_1) = -H_{23} W_1, \qquad \Rm(W_2) = -H_{23} W_2.
\end{align}

These formulas determine the matrix $\sfK_j$. Indeed, writing
\begin{align}
    \Rm(\tau_{n'}) = \sum_{n''=1}^{n_j} k_{n''n'}(s) \tau_{n''},
\end{align}
the definition \eqref{eq:sfK-defn} of $\sfK_j$ and the orthonormality of the $\tau_n$'s (Proposition \ref{prop:spherical-sector-basis}) give
\begin{align}
    (\sfK_j)_{nn'} &= \int_{G/H} \inner{\Rm(\tau_{n'})}{\tau_n} \, d\omega = \sum_{n''=1}^{n_j} k_{n''n'}(s) \int_{G/H} \inner{\tau_{n''}}{\tau_n} \, d\omega = k_{nn'}(s).
\end{align}
The above calculations therefore yield, for $j \geq 2$,
\begin{align}\label{eq:Kj-matrix}
    \sfK_j=\operatorname{diag}\left(
    \begin{pmatrix}
        0&0&K_1&\sqrt2K_2\\
        0&-K_1&0&0\\
        K_1&0&0&\sqrt2H_{12}\\
        \sqrt2K_2&0&\sqrt2H_{12}&H_{23}
    \end{pmatrix},
    \begin{pmatrix}
        -K_2&0&0&3M_2\\
        0&-K_2&-3M_2&0\\
        0&-3M_2&-H_{12}&0\\
        3M_2&0&0&-H_{12}
    \end{pmatrix},
    -H_{23}I_2
    \right).
\end{align}
If $j=1$, omit the last block; if $j=0$, omit the last two blocks.

\begin{proof}[Proof of Lemma \ref{lem:spherical-angular-matrix}]
    This follows from \eqref{eq:Kj-matrix}, the earlier calculation \eqref{eq:Gj-matrix} of $\sfG_j$, and the formula $\sfC_j = \sfG_j - 2\sfK_j$. The only additional simplification occurs in the $W$-block of $\sfC_j$, where we have used $H_{23} = 4/b^2 - 3D^2 - B^2$ to simplify that block to $(\lambda-8)b^{-2} I_2$.
\end{proof}

\subsection{Reduction to a pointwise positivity problem} \label{subsec:reduction-pointwise}

The goal now is to prove the stability inequality in $(j,0)$-sectors by establishing the sufficient criterion in Corollary \ref{cor:Cj-positivity-implies-stability}. Discarding the $|X'|^2$ term in \eqref{eq:radial-coercivity-target} would work if $\sfC_j \succeq \delta b^{-2} I$ uniformly in $j$ and $s$. However, this already fails for $j=0$: take the constant vector $X_g = (1,0,1,\sqrt{2})^T$, which corresponds to the tensor
\begin{align}
    h_g = U_1[Y] + U_3[Y] + \sqrt{2} U_4[Y] = c_k Y ((\epsilon^0)^2 + (\epsilon^1)^2 + g_H) = c_k Y g,
\end{align}
where $Y$ is a constant function on $\CP^1$. Since $\nabla g = 0$ and $\Rm(g) = \Ric$, we have
\begin{align}
    X_g^* \sfC_0(s) X_g &= \mathfrak{a}_s(h_g,h_g) = \int_{L_{2,k}} \left( \sum_{\alpha=1}^3 |\nabla_{F_\alpha} h_g|^2 - 2\inner{\Rm(h_g)}{h_g} \right) \, d\omega \\
    &= -2R(s) \int_{L_{2,k}} c_k^2 |Y|^2 \, d\omega < 0,
\end{align}
and hence $\sfC_0$ cannot be $\succeq \delta b^{-2} I$. Instead, we will rebalance the left-hand side of \eqref{eq:radial-coercivity-target} as follows, then prove positivity of the remainder matrix $\calR_j$.

\begin{lemma} \label{lem:spherical-radial-completion}
    Define
    \begin{align}
        Q=\frac{a}{b}, \qquad \theta = \frac{Q(1-Q^2)}{b^2 Q'}, \qquad \psi = f' + 2B + \theta, \qquad \gamma = \frac{f'}{2} + 3B.
    \end{align}
    Set
    \begin{gather}
        \Psi_U=
        \begin{pmatrix}
            \frac{\psi+\gamma}{2}&0
                &\frac{\psi-\gamma}{2}&-\sqrt2B\\
            0&\gamma&0&0\\
            \frac{\psi-\gamma}{2}&0
                &\frac{\psi+\gamma}{2}&-\sqrt2B\\
            -\sqrt2B&0&-\sqrt2B&\gamma
        \end{pmatrix},
        \qquad
        \Psi_j=\operatorname{diag}
            (\Psi_U,\gamma I_4,\gamma I_2), \\
        \mathbf r_U=
        \begin{pmatrix}
            \theta/\sqrt2\\0\\\theta/\sqrt2\\2B/3
        \end{pmatrix},
        \qquad
        \mathbf r_j= \begin{pmatrix} \bfr_U \\ 0_{4 \times 1} \\ 0_{2 \times 1} \end{pmatrix},
    \end{gather}
    where absent groups in $\Psi_j$ and $\bfr_j$ are omitted if $j=0$ or $j=1$.
    Then every $X \in C^\infty_c((0,\infty); \C^{n_j})$, written as $X = (X_U,X_V,X_W)$ (absent groups omitted) with $X_U = (u_1,u_2,u_3,u_4)^T$, satisfies
    \begin{align}\label{eq:spherical-radial-completion}
        \int_0^\infty
            w\left(|X'|^2+X^*\sfC_jX\right)\,ds
        ={}&
        \int_0^\infty w|X'+\Psi_jX|^2\,ds
        +3\int_0^\infty w
            \left|
                \frac{\theta}{\sqrt2}(u_1+u_3)
                +\frac23Bu_4
            \right|^2\,ds
        \notag\\
        &+\int_0^\infty\frac{w}{b^2}
            X^*\calR_jX\,ds,
    \end{align}
    where
    \begin{align}\label{eq:spherical-remainder-defn}
        \calR_j
        :=b^2\left(
            \sfC_j+\Psi_j'
            +\frac{w'}w\Psi_j-\Psi_j^2
            -3\mathbf r_j\mathbf r_j^*
        \right).
    \end{align}
\end{lemma}
\begin{proof}
    Since $\Psi_j$ is Hermitian, we have $(X^*\Psi_j X)' = 2\Re(X'^*\Psi_j X) + X^*\Psi_j' X$. Thus, expanding the square and integrating by parts gives
    \begin{align}
        \int_0^\infty w|X' + \Psi_j X|^2 \, ds &= \int_0^\infty w \left[ |X'|^2 + X^*(\Psi_j^2 - \Psi_j') X + (X^* \Psi_j X)' \right] \, ds \\
        &= \int_0^\infty w\left[ |X'|^2 + X^*\left( \Psi_j^2 - \Psi_j' - \frac{w'}{w} \Psi_j \right) X \right] \, ds.
    \end{align}
    The boundary term vanishes as $X$ has compact support in $(0,\infty)$. Rearranging and using \eqref{eq:spherical-remainder-defn}, we get
    \begin{align}
        \int_0^\infty w\left( |X'|^2 + X^* \sfC_j X \right) \, ds ={} & \int_0^\infty w|X' + \Psi_j X|^2 \, ds + 3 \int_0^\infty w|\bfr_j^*X|^2 \, ds \\
        &+ \int_0^\infty \frac{w}{b^2} X^* \calR_j X \, ds.
    \end{align}
    Finally, $\bfr_j^*X = \frac{\theta}{\sqrt{2}}(u_1+u_3) + \frac{2}{3} Bu_4$, proving the claimed identity.
\end{proof}

The matrix $\sfC_j$ was computed in Lemma \ref{lem:spherical-angular-matrix}. We now calculate the other parts of $\calR_j$ in terms of the same quantities used to express $\sfC_j$. This requires an application of soliton identities.
\begin{lemma} \label{lem:Psi-calculation}
    Let $\Psi_j$ and $\bfr_j$ be defined as in Lemma \ref{lem:spherical-radial-completion}. Then
    \begin{equation}\label{eq:spherical-radial-correction}
        b^2\left(
            \Psi_j'+\frac{w'}w\Psi_j-\Psi_j^2
            -3\mathbf r_j\mathbf r_j^*
        \right) =
        \operatorname{diag}\left(
        \begin{pmatrix}
            \frac{d+\eta}{2}&0&\frac{d-\eta}{2}
                &\frac{e}{\sqrt2}\\
            0&\eta&0&0\\
            \frac{d-\eta}{2}&0&\frac{d+\eta}{2}
                &\frac{e}{\sqrt2}\\
            \frac{e}{\sqrt2}&0&\frac{e}{\sqrt2}
                &\eta-\frac{16}{3}b^2B^2
        \end{pmatrix},
        \eta I_4,\eta I_2
        \right),
    \end{equation}
    with absent groups omitted. The quantities $d,e,\eta$ are defined by
    \begin{align}
        d &= 2 + b^2 (2K_1-6B^2-2D^2), \\
        e &= b^2(2K_2 + Bf' + 8B^2 - 2AB), \\
        \eta &= b^2 \left( \frac{1}{2}K_1 - 2K_2 - 6B^2 + \frac{(f')^2}{4} + \frac{1}{2} Af' + Bf' + 3AB \right).
    \end{align}
\end{lemma}
\begin{proof}
    Define the following orthonormal basis of $\C^4$:
    \begin{equation} \label{eq:rotated-ONB-C4}
        \bfp = \frac{1}{\sqrt{2}}(1,0,1,0)^T, \quad \bfe_2 = (0,1,0,0)^T, \quad \bfq = \frac{1}{\sqrt{2}} (1,0,-1,0)^T, \quad \bfe_4 = (0,0,0,1)^T.
    \end{equation}
    The subspace $\Span\{\bfp,\bfe_4\}$ contains $\bfr_U$ and is preserved by $\Psi_U$. On this subspace, $\Psi_U$ and $\bfr_U$ have the form
    \begin{align}
        \Psi_U|_{\Span\{\bfp,\bfe_4\}} = \begin{pmatrix} \psi & -2B \\ -2B & \gamma \end{pmatrix}, \qquad \bfr_U = \begin{pmatrix} \theta \\ 2B/3 \end{pmatrix}.
    \end{align}
    A calculation gives that on this subspace,
    \begin{align}
        &b^2\left( \Psi_U' + \frac{w'}{w} \Psi_U - \Psi_U^2 - 3\bfr_U \bfr_U^* \right) \\
        &= b^2 \begin{pmatrix}
            \psi' + \frac{w'}{w} \psi - \psi^2 - 4B^2 - 3\theta^2 & -2B' + 2B\left(\psi + \gamma - \frac{w'}{w} - \theta \right) \\
            -2B' + 2B\left( \psi + \gamma - \frac{w'}{w} - \theta \right) & \gamma' + \frac{w'}{w} \gamma - \gamma^2 - \frac{16}{3}B^2
        \end{pmatrix}.
    \end{align}
    We claim that this is equal to
    \begin{align}
        \begin{pmatrix}
            d & e \\ e & \eta - \frac{16}{3} b^2 B^2
        \end{pmatrix}.
    \end{align}
    Indeed, this will imply, by comparison of the bottom right entries, that $\eta = b^2(\gamma' + \frac{w'}{w}\gamma - \gamma^2)$. Since $\Psi_U$ acts as $\gamma I$ on the orthogonal complement of $\Span\{\bfp,\bfe_4\}$ (this is easily verified), we will have, in the basis \eqref{eq:rotated-ONB-C4},
    \begin{align}
        b^2\left(
            \Psi_U'+\frac{w'}w\Psi_U-\Psi_U^2
            -3 \bfr_U \bfr_U^*
        \right)
        =
        \begin{pmatrix}
            d & 0 & 0 & e \\
            0 & \eta & 0 & 0 \\
            0 & 0 & \eta & 0 \\
            e & 0 & 0 & \eta-\frac{16}{3} b^2 B^2
        \end{pmatrix}.
    \end{align}
    Changing back to the standard basis gives the first block in \eqref{eq:spherical-radial-correction}. The other blocks are trivially $\eta I$ since $\Psi_j$ acts as $\gamma I$ on the corresponding subspaces, and $\eta = b^2 (\gamma' + \frac{w'}{w}\gamma- \gamma^2)$.

    Thus, to prove the lemma, it remains to show that
    \begin{align}
        d &= b^2 \left( \psi' + \frac{w'}{w}\psi - \psi^2 - 4B^2 - 3\theta^2 \right), \label{eq:claimed-d-identity} \\
        e &= b^2 \left( -2B' + 2B \left( \psi + \gamma - \frac{w'}{w} - \theta \right) \right), \label{eq:claimed-e-identity} \\
        \eta &= b^2 \left( \gamma' + \frac{w'}{w}\gamma - \gamma^2 \right). \label{eq:claimed-sigma-identity}
    \end{align}
    To prove these identities, first use \eqref{eq:appleton-Q-ODE} with $m=2$ to get $\frac{Q''}{Q'} = -f' - 3B - 4\theta$. Noting that $\theta$ is smooth and positive because $0<Q<1$ and $Q'>0$ on $(0,\infty)$, it follows that
    \begin{align}
        \frac{\theta'}{\theta} &= (\log\theta)' = (\log Q)' + (\log(1-Q^2))' - 2(\log b)' - (\log Q')' \\
        &= A - B - \frac{2QQ'}{1-Q^2} - 2B + f' + 3B + 4\theta.
    \end{align}
    Multiplying both sides by $\theta$ gives $\theta' = (A+f')\theta + 4\theta^2 - 2Q^2/b^2$. Meanwhile, \eqref{eq:appleton-soliton-f} gives
    \begin{align} \label{eq:soliton-consequence-1}
        f'' + 2B' = K_1 - 2B^2.
    \end{align}
    Thus,
    \begin{align} \label{eq:psi-prime}
        \psi' = f'' + 2B' + \theta' = K_1 - 2B^2 + (A+f')\theta + 4\theta^2 - \frac{2Q^2}{b^2}.
    \end{align}
    Since $\frac{w'}{w} = f'+A+2B$ and $\psi = f'+2B+\theta$, we have $\frac{w'}{w}\psi - \psi^2 = (A-\theta)(f'+2B+\theta)$. Together with \eqref{eq:psi-prime}, and the fact that $D^2 = Q^2/b^2$, we get
    \begin{align}
        \psi' + \frac{w'}{w} \psi - \psi^2 + 2B^2 + 2D^2 - 2K_1 = -K_1 + Af' + 2AB + 2(A-B)\theta + 3\theta^2.
    \end{align}
    Now, $(A-B)\theta = \frac{Q'}{Q} \frac{Q(1-Q^2)}{b^2 Q'} = \frac{1-Q^2}{b^2}$, and \eqref{eq:appleton-soliton-a} implies $K_1 = Af' + 2AB - 2Q^2/b^2$. Plugging these into the right-hand side above and rearranging, we arrive at the identity \eqref{eq:claimed-d-identity}. Next, note that $B' = -K_2 - B^2$ and $\psi+\gamma-\frac{w'}{w}-\theta = \frac{f'}{2}+3B-A$. Thus,
    \begin{align}
        -2B' + 2B\left( \psi+\gamma - \frac{w'}{w} - \theta \right) &= 2K_2 + 2B^2 + Bf' + 6B^2 - 2AB = \frac{e}{b^2},
    \end{align}
    which proves \eqref{eq:claimed-e-identity}. Lastly, using \eqref{eq:soliton-consequence-1}, we get
    \begin{align}
        \gamma' &= \frac{1}{2}(f'' + 2B') + 2B' = \frac{1}{2}K_1 - 2K_2 - 3B^2.
    \end{align}
    Moreover,
    \begin{align}
        \gamma \left( \frac{w'}{w}-\gamma \right) &= \left( \frac{f'}{2} + 3B \right) \left( \frac{f'}{2} + A - B \right) = \frac{(f')^2}{4} + \frac{1}{2} Af' + Bf' + 3AB - 3B^2.
    \end{align}
    Adding these expressions gives \eqref{eq:claimed-sigma-identity}, completing the proof.
\end{proof}

Combining Lemmas \ref{lem:Psi-calculation} and \ref{lem:spherical-angular-matrix}, we obtain the entries of $\calR_j$ \eqref{eq:spherical-remainder-defn}. We now decompose the blocks of $\calR_j$ and introduce two matrices $\sfA, \sfB$ whose entries are expressed in the variables $x,v,y,Z$ from \S\ref{subsec:new-variables} and \eqref{eq:Z-epsilon-defn}. The definitions together with \eqref{eq:background-algebraic-identities} give
\begin{align} \label{eq:spherical-background-conversion}
    bA = \frac{a'}{\sqrt{x}} = \frac{y-xZ}{2\sqrt{x}}, \qquad bB = v\sqrt{x}, \qquad bD = \sqrt{x}, \qquad bf' = \sqrt{x}(Z-2v).
\end{align}
Setting $\lambda = 4j(j+1)$ and using these identities, the result is
\begin{align}\label{eq:spherical-remainder-matrix}
    \calR_j=
    \begin{pmatrix}
        \sfA+\lambda I_4
            &2b\sqrt{\lambda}\, \sfM&0\\
        2b\sqrt{\lambda}\, \sfM^T
            &\sfB+4b^2 \sfM^T \sfM+(\lambda-8)I_4
            &2b\sqrt{\lambda-8}\, \sfN\\
        0&2b\sqrt{\lambda-8}\, \sfN^T
            &(\eta+\lambda-8)I_2
    \end{pmatrix}.
\end{align}
Here $\sfM,\sfN$ are the matrices from Lemma \ref{lem:spherical-angular-matrix}. If $j=1$, omit the $W$ rows and columns and take $\lambda=8$; if $j=0$, we have $\calR_0 = \sfA$.

The scalar $\eta$ is defined in Lemma \ref{lem:Psi-calculation}. After substituting \eqref{eq:scaled-curvature-coefficients} and \eqref{eq:background-algebraic-identities}, it is given by
\begin{align}\label{eq:spherical-sigma}
    \eta = 8 + \frac{yZ}{2} - x\left( \frac{Z^2}{4} + 2vZ + 5 + 5v^2 \right).
\end{align}
The matrices $\sfA,\sfB$ are $4 \times 4$ real symmetric with possibly nonzero upper-triangular entries
\begin{align}\label{eq:spherical-AU-entries}
    \sfA_{11}
    &=\frac{y(y-xZ)}{2x}
      +\frac{\eta}{2}+1-3x+v^2x,\\
    \sfA_{22}
    &=\frac{y(y-xZ)}{x}
      +\eta-2x+2v^2x,\\
    \sfA_{33}
    &=\frac{y(y-xZ)}{2x}
      +\frac{\eta}{2}+1+x-3v^2x,\\
    \sfA_{44}
    &=\eta+8x-8-\frac43v^2x,
      \label{eq:A44}\\
    \sfA_{13}
    &=-\frac{y(y-xZ)}{2x}
      -\frac{\eta}{2}+1+x-3v^2x,\\
    \sfA_{14}
    &=\frac{1}{\sqrt2}
      \left(8-4x-2vy+vxZ+4v^2x\right),\\
    \sfA_{34}
    &=\frac{1}{\sqrt2}
      \left(-8-4x+2vy+vxZ+4v^2x\right),\\[1mm]
    \sfB_{11}
    &=\frac{y^2}{4x}+vy-vxZ
      +4+x-14v^2x,\\
    \sfB_{22}
    &=\frac{y^2}{4x}+vy-vxZ
      +4-3x-2v^2x,\\
    \sfB_{33}
    &=\frac{y^2}{4x}-vy-vxZ
      +12+3x-8v^2x,\\
    \sfB_{44}
    &=\frac{y^2}{4x}-vy-vxZ
      +12-9x-4v^2x,\\
    \sfB_{14}
    &=-2y+2xZ+12vx,\\
    \sfB_{23}
    &=2y-2xZ-4vx.
\end{align}

\subsection{Proof of Theorem \ref{thm:neutral-coercivity}} \label{subsec:neutral-proof}

In this subsection, we prove that $\calR_j \succeq \delta I$ for all integers $j \geq 0$ and a uniform $\delta > 0$, then use this fact to prove Theorem \ref{thm:neutral-coercivity}. The matrix $\calR_j$ is explicitly given at the end of \S\ref{subsec:reduction-pointwise}. To estimate it, we first estimate $\sfA$ and $\sfB$ as follows.
\begin{lemma} \label{lem:spherical-background-bounds}
    For every $s>0$, we have
    \begin{align}
        \sfA(s) \succeq \frac{1}{64} I_4, \qquad \sfB(s) \succeq \frac{1}{64} I_4.
    \end{align}
\end{lemma}
\begin{proof}
    We will establish the inequalities $\sfA, \sfB \succeq \frac{1}{64} I_4$ throughout the region
    \begin{align}
        \calD := \left\{ (x,v,y,Z) \in \R^4 :
        \begin{array}{l}
            0<x<1,\quad 0<v\leq2/3,\quad y\geq6,\\
            vZ\geq2,\quad y-xZ\geq2xv,\\
            yv\leq4-x+xv^2
        \end{array} \right\}.
    \end{align}
    This suffices to prove the lemma, since for every $s>0$, the variables $x,v,y,Z$ on the Appleton soliton lie in $\calD$. Indeed, \eqref{eq:4d-basic-bounds} gives $0<x<1$, $0<v<2/k \leq 2/3$, $y > 2k \geq 6$, and $y-xZ = 2a' = 2z+2xv > 2xv$, whereas Lemma \ref{lem:sharp-background-bounds} gives $vZ > 2$ and $yv = 4-x+xv^2-\xi < 4-x+xv^2$.

    Holding $(x,v,y)$ fixed, use the explicit entries of $\sfA$ and $\sfB$ listed earlier to compute
    \begin{align}
        \del_Z^2 \sfA = -\frac{x}{2} (I_4 - \bfp \bfp^T) \preceq 0, \qquad \partial_Z^2 \sfB = 0,
    \end{align}
    where $\bfp = \frac{1}{\sqrt{2}}(1,0,1,0)^T$. Thus, $\sfA$ and $\sfB$ are both concave in $Z$ with respect to the order $\preceq$. It therefore suffices to establish that $\sfA, \sfB \succeq \frac{1}{64} I_4$ at the smallest and largest $Z$ in $\calD$ for each fixed $(x,v,y)$, i.e. at
    \begin{align} \label{eq:Z-endpoints}
        Z_L = \frac{2}{v}, \qquad Z_R = \frac{y}{x}-2v.
    \end{align}
    To verify these inequalities, define $W=vy$. Whenever $(x,v,y,Z) \in \calD$, the triple $(x,v,W)$ belongs to
    \begin{align}
        \calD_* := \left\{ (x,v,W) \in \R^3 :
            \begin{array}{l}
                0<x\leq1,\quad 0<v\leq2/3,\\
                6v\leq W,\quad
                2x(1+v^2)\leq W\leq4-x+xv^2
            \end{array} \right\}.
    \end{align}
    Indeed, the bounds $2/v \leq Z \leq y/x - 2v$ imply $2x \leq yv-2xv^2 = W-2xv^2$, and hence $W \geq 2x(1+v^2)$.

    Let $\sfA_L, \sfA_R, \sfB_L, \sfB_R$ be the matrices obtained by substituting the lower and upper bounds for $Z$ \eqref{eq:Z-endpoints}, i.e.
    \begin{gather}
        \sfA_L = \sfA|_{Z=2/v}, \qquad \sfA_R = \sfA|_{Z=\frac{W}{vx}-2v}, \\
        \sfB_L = \sfB|_{Z=2/v}, \qquad \sfB_R = \sfB|_{Z=\frac{W}{vx}-2v}.
    \end{gather}
    The upper bound for $Z$ comes from substituting $y=W/v$ into $\frac{y}{x}-2v$. All entries of these four matrices are (real-valued) rational functions of $x,v,W$. We thus need to check that all four matrices $\succeq \frac{1}{64} I_4$ in $\calD_*$.

    For a matrix $H$ and an index set $I \subset \{1,2,3,4\}$, let $H^I$ be the square matrix obtained by taking the $i$-th row and column for each $i \in I$. We verify that the following 16 inequalities hold everywhere on $\calD_*$, with $\delta = 1/64$:
    \begin{align}
        \det(\sfA_L^{\{2\}}-\delta I)&>0,
        &\det(\sfA_R^{\{2\}}-\delta I)&>0,\\
        \det(\sfA_L^{\{1\}}-\delta I)&>0,
        &\det(\sfA_R^{\{1\}}-\delta I)&>0,\\
        \det(\sfA_L^{\{1,3\}}-\delta I)&>0,
        &\det(\sfA_R^{\{1,3\}}-\delta I)&>0,\\
        \det(\sfA_L^{\{1,3,4\}}-\delta I)&>0,
        &\det(\sfA_R^{\{1,3,4\}}-\delta I)&>0,\\[1mm]
        \det(\sfB_L^{\{1\}}-\delta I)&>0,
        &\det(\sfB_R^{\{1\}}-\delta I)&>0,\\
        \det(\sfB_L^{\{1,4\}}-\delta I)&>0,
        &\det(\sfB_R^{\{1,4\}}-\delta I)&>0,\\
        \det(\sfB_L^{\{2\}}-\delta I)&>0,
        &\det(\sfB_R^{\{2\}}-\delta I)&>0,\\
        \det(\sfB_L^{\{2,3\}}-\delta I)&>0,
        &\det(\sfB_R^{\{2,3\}}-\delta I)&>0.
    \end{align}
    This is done in the supplementary Wolfram Mathematica notebook accompanying this paper. There, we construct $\sfA_L, \sfA_R, \sfB_L, \sfB_R$ from their definitions and compute the listed determinants. Each determinant simplifies to a rational function $p/q$ of $(x,v,W)$. It is obvious from these expressions that $q>0$ throughout $\calD_*$, while a cylindrical algebraic decomposition is used to show that $p>0$ throughout $\calD_*$. This is an exact symbolic procedure completed in finitely many steps, so there are no subtleties arising from approximation.

    As shown at the end of \S\ref{subsec:reduction-pointwise}, the matrix $\sfA$ preserves the coordinate subspaces indexed by $\{2\}$ and $\{1,3,4\}$, and the matrix $\sfB$ preserves the coordinate subspaces indexed by $\{1,4\}$ and $\{2,3\}$. It then follows from these inequalities, and Sylvester's criterion, that $\sfA_L, \sfA_R, \sfB_L, \sfB_R \succeq \delta I_4$ throughout $\calD_*$.
\end{proof}

We conclude this section with the proof of Theorem \ref{thm:neutral-coercivity}.

\begin{proof}[Proof of Theorem \ref{thm:neutral-coercivity}]
    Let $\delta = 1/64$. We first show that $\calR_j \succeq \delta I$ for every $j \geq 0$. By Lemma \ref{lem:spherical-background-bounds}, $\sfA \succeq \delta I_4$ and $\sfB \succeq \delta I_4$. For $j=0$, we have $\calR_0 = \sfA$, so the desired bound is immediate.
    
    Suppose now that $j \geq 1$, and let $\lambda = 4j(j+1)$. The formula \eqref{eq:A44} for $\sfA_{44}$ gives
    \begin{align}
        \eta = \sfA_{44} + 8(1-x) + \frac{4}{3}v^2 x.
    \end{align}
    Since $\sfA_{44} \geq \delta$ and $0<x<1$, this gives $\eta \geq \delta$. Let $X = (X_U,X_V,X_W) \in \C^{n_j}$, omitting $X_W$ when $j=1$. Using \eqref{eq:spherical-remainder-matrix} and the preceding bounds, we obtain, if $j \geq 2$,
    \begin{align}
        X^*(\calR_j - \delta I)X \geq{}& \lambda|X_U|^2 + 4b\sqrt{\lambda} \, \Re\inner{\sfM X_V}{X_U} + 4b^2 |\sfM X_V|^2 \\
        & + (\lambda-8)(|X_V|^2 + |X_W|^2) + 4b\sqrt{\lambda-8} \Re\inner{\sfN X_W}{X_V} \\
        ={} & \left| \sqrt{\lambda} X_U + 2b \sfM X_V \right|^2 \\
        & + (\lambda-8)(|X_V|^2 + |X_W|^2) + 4b\sqrt{\lambda-8} \, \Re\inner{\sfN X_W}{X_V}. \label{eq:spherical-final-square}
    \end{align}
    For $j=1$, the same calculation applies with $X_W$ absent and $\lambda = 8$, leaving only the nonnegative square on the right. Thus $\calR_1 \succeq \delta I$.

    For $j \geq 2$, the formula for $\sfN$ from Lemma \ref{lem:spherical-angular-matrix} gives $\sfN^T \sfN = (B^2 + D^2) I_2$. Using that $0<x<1$ and $0<v<2/3$, we also have $b^2 (B^2 + D^2) = x(1+v^2) < 13/9$. Moreover, $\lambda-8 = 4j(j+1)-8 \geq 16$. Thus, we can bound the operator norm of $\sfN$ as follows:
    \begin{align}
        2b\sqrt{\lambda-8} \, \norm{\sfN}_{\mathrm{op}} = 2\sqrt{x(1+v^2)(\lambda-8)} \leq \lambda-8.
    \end{align}
    By the Cauchy--Schwarz and absorbing inequalities, we therefore have
    \begin{align}
        4b\sqrt{\lambda-8} \, \Re\inner{\sfN X_W}{X_V} &\geq -(\lambda-8) (|X_V|^2 + |X_W|^2).
    \end{align}
    Substituting this into \eqref{eq:spherical-final-square} shows that $\calR_j \succeq \delta I$.

    Thus, by Lemma \ref{lem:spherical-radial-completion}, for any integer $j \geq 0$ and $X \in C^\infty_c((0,\infty); \C^{n_j})$, we have
    \begin{align}
        \int_0^\infty w\left( |X'|^2 + X^*\sfC_j X \right) \, ds \geq \int_0^\infty \frac{w}{b^2} X^* \calR_j X \, ds \geq \delta \int_0^\infty \frac{w}{b^2} |X|^2 \, ds.
    \end{align}
    Combining this with Corollary \ref{cor:Cj-positivity-implies-stability} proves Theorem \ref{thm:neutral-coercivity}.
\end{proof}

\section{Conclusion: strict linear stability} \label{sec:conclusion}

\begin{proof}[Proof of Theorem \ref{thm:appleton-stability}]
    Fix $k \geq 3$, and let $(M,g,f)$ be the four-dimensional Appleton soliton on $\O_{\CP^1}(-k)$. By Theorems \ref{thm:charged-coercivity} and \ref{thm:neutral-coercivity}, every smooth symmetric $2$-tensor $h$ compactly supported in $M^\circ$, i.e. away from the zero section, and belonging to a single admissible sector, satisfies
    \begin{align}
        \calE[h] \geq \frac{1}{64} \int_M \frac{|h|^2}{b^2} e^f.
    \end{align}
    Proposition \ref{prop:sector-stability-reduction} then implies that $(M,g,f)$ is strictly linearly stable.
\end{proof}

\bibliography{Refs}

\begin{bibdiv}
\begin{biblist}

\bib{appleton}{article}{
      author={Appleton, A.},
       title={A family of non-collapsed steady {R}icci solitons in even dimensions greater or equal to four},
        date={2017},
     journal={Preprint, arXiv:1708.00161},
}

\bib{bamler-kleiner-m1}{article}{
      author={Bamler, R.~H.},
      author={Kleiner, B.},
       title={On the multiplicity one conjecture for mean curvature flows of surfaces},
        date={2023},
     journal={Preprint, arXiv:2312.02106},
}

\bib{boucetta-sphere}{article}{
      author={Boucetta, M.},
       title={Spectra and symmetric eigentensors of the {L}ichnerowicz {L}aplacian on {$S^n$}},
        date={2009},
     journal={Osaka J. Math.},
      volume={46},
      number={1},
       pages={235\ndash 254},
}

\bib{cao-hamilton-ilmanen}{article}{
      author={Cao, H.-D.},
      author={Hamilton, R.},
      author={Ilmanen, T.},
       title={Gaussian densities and stability for some {R}icci solitons},
        date={2004},
     journal={Preprint, arXiv:math/0404165},
}

\bib{cao-he}{article}{
      author={Cao, H.-D.},
      author={He, C.},
       title={Linear stability of {P}erelman's $\nu$-entropy on symmetric spaces of compact type},
        date={2015},
     journal={J. Reine Angew. Math.},
      volume={709},
       pages={229\ndash 246},
}

\bib{ccms}{article}{
      author={Chodosh, O.},
      author={Choi, K.},
      author={Mantoulidis, C.},
      author={Schulze, F.},
       title={Mean curvature flow with generic initial data},
        date={2024},
     journal={Invent. Math.},
      volume={237},
       pages={121\ndash 220},
}

\bib{cms}{article}{
      author={Chodosh, O.},
      author={Choi, K.},
      author={Schulze, F.},
       title={Mean curvature flow with generic initial data {II}},
        date={2023},
     journal={Preprint, arXiv:2302.08409},
}

\bib{cm-generic}{article}{
      author={Colding, T.~H.},
      author={Minicozzi~II, W.~P.},
       title={Generic mean curvature flow {I}; generic singularities},
        date={2012},
     journal={Ann. of Math.},
      volume={175},
      number={2},
       pages={755\ndash 833},
}

\bib{deng-wang}{article}{
      author={Deng, Y.},
      author={Wang, G.},
       title={Classification of cohomogeneity one steady gradient {R}icci solitons in dimension four and five},
        date={2026},
     journal={Preprint, arXiv:2611.08206},
}

\bib{deruelle-stability}{article}{
      author={Deruelle, A.},
       title={Stability of non compact steady and expanding gradient {R}icci solitons},
        date={2015},
     journal={Calc. Var. Partial Differential Equations},
      volume={54},
       pages={2367\ndash 2405},
}

\bib{deruelle-kroncke}{article}{
      author={Deruelle, A.},
      author={Kr\"oncke, K.},
       title={Stability of {ALE} {R}icci-flat manifolds under {R}icci flow},
        date={2021},
     journal={J. Geom. Anal.},
      volume={31},
       pages={2829\ndash 2870},
}

\bib{deruelle-ozuch}{article}{
      author={Deruelle, A.},
      author={Ozuch, T.},
       title={Dynamical (in)stability of {R}icci-flat {ALE} metrics along the {R}icci flow},
        date={2023},
     journal={Calc. Var. Partial Differential Equations},
      volume={62},
      number={84},
}

\bib{ivey-analytic}{article}{
      author={Ivey, T.~A.},
       title={Local existence of {R}icci solitons},
        date={1996},
     journal={Manuscripta Math.},
      volume={91},
       pages={151\ndash 162},
}

\bib{kim-ozuch}{article}{
      author={Kim, D.},
      author={Ozuch, T.},
       title={Blowdowns along {R}icci flow in dimension four},
        date={2026},
     journal={Preprint, arXiv:2611.00139},
}

\bib{knopf-sesum-cpn}{article}{
      author={Knopf, D.},
      author={Šešum, N.},
       title={Dynamical instability of {$\mathbb{CP}^N$} under {R}icci flow},
        date={2019},
     journal={J. Geom. Anal.},
      volume={29},
       pages={902\ndash 916},
}

\bib{kroncke-cpn}{article}{
      author={Kr\"oncke, K.},
       title={Stability of {E}instein metrics under {R}icci flow},
        date={2020},
     journal={Comm. Anal. Geom.},
      volume={28},
      number={2},
       pages={351\ndash 394},
}

\bib{lavoyer-thompson}{article}{
      author={Lavoyer, L.},
      author={Thompson, A.},
       title={Linear stability of steady and expanding {K}\"ahler-{R}icci solitons},
        date={2026},
     journal={Preprint, arXiv:2609.21472},
}

\bib{symmetry-law}{article}{
      author={Law, M.~B.},
       title={On steady and expanding {R}icci solitons with asymptotic symmetries},
        date={2026},
     journal={Preprint, arXiv:2505.20576},
}

\bib{naff-ozuch-kahler}{article}{
      author={Naff, K.},
      author={Ozuch, T.},
       title={Linear stability and instability of {K}\"ahler {R}icci solitons},
        date={2025},
     journal={Preprint, arXiv:2511.15885},
}

\bib{naff-ozuch-fik}{article}{
      author={Naff, K.},
      author={Ozuch, T.},
       title={Linear stability of the blowdown {R}icci shrinker in 4{D}},
        date={2025},
     journal={Preprint, arXiv:2509.05470},
}

\bib{ribeiro-tran}{article}{
      author={Ribeiro~Jr., E.},
      author={Tran, H.},
       title={Steady conformally {K}\"ahler {R}icci solitons with isoparametric potential},
        date={2026},
     journal={Preprint, arXiv:2610.40994},
}

\bib{wallach}{book}{
      author={Wallach, N.~R.},
       title={Harmonic analysis on homogeneous spaces},
   publisher={Dover Publications},
        date={2018},
}

\end{biblist}
\end{bibdiv}
\end{document}